\documentclass[aos,nameyear,dvips]{arximspdf}
\usepackage{graphicx}

\doi{10.1214/10-AOS795}
\volume{38}
\issue{6}
\pubyear{2010}
\firstpage{3487}
\lastpage{3566}

\makeatletter
\renewcommand{\epsilon}{\varepsilon}
\DeclareMathAlphabet\mathcaligr{OMS}{cmsy}{m}{n}
\newcommand{\cal}{\mathcaligr}

\newcommand{\eps}{\epsilon}

\newcommand{\myExp}{\mathbf{E}}
\newcommand{\iid}{\stackrel{i.i.d.}{\backsim}}
\newcommand{\tendsto}{\rightarrow}
\newcommand{\me}{\mathrm{e}}
\newcommand{\id}{\mathrm{Id}}
\newcommand{\gO}{\mathrm{O}}
\newcommand{\lo}{\mathrm{o}}

\newtheorem{theorem}{Theorem}[section]
\newtheorem{lemma}{Lemma}[section]
\newtheorem{lemmaNumThm}[theorem]{Lemma}
\newtheorem{fact}{Fact}[section]
\newproclaim{factNumThm}[theorem]{Fact}
\newtheorem{corollaryNumThm}[theorem]{Corollary}
\newtheorem{corollary}[theorem]{Corollary}
\newtheorem{propositionNumThm}[theorem]{Proposition}

\newcommand{\opnorm}[1]{|\!|\!|#1|\!|\!|_2}

\newcommand{\trsp}{'}
\newcommand{\Wishart}{\cal{W}}
\newcommand{\Poisson}{\mathrm{Po}}
\newcommand{\WeakCv}{\Longrightarrow}
\newcommand{\equalInLaw}{\stackrel{{\cal L}}{=}}
\newcommand{\SigmaHat}{\hat{\Sigma}}
\newcommand{\ebold}{\mathbf{e}}
\newcommand{\muHat}{\hat{\mu}}
\newcommand{\muTilde}{\widetilde{\mu}}
\newcommand{\wOpt}{w_{\mathrm{optimal}}}
\newcommand{\fTheo}{f_{\mathrm{theo}}}
\newcommand{\fEmp}{f_{\mathrm{emp}}}
\newcommand{\wTheo}{w_{\mathrm{theo}}}
\newcommand{\wEmp}{w_{\mathrm{emp}}}
\newcommand{\wOracle}{w_{\mathrm{oracle}}}
\newcommand{\myVector}[2]{\pmatrix{ #1 \cr #2 }}
\newcommand{\VHat}{\widehat{V}}
\newcommand{\omegahat}{\widehat{\omega}}
\newcommand{\MHat}{\widehat{M}}
\newcommand{\NHat}{\widehat{N}}
\newcommand{\MTilde}{\widetilde{M}}
\newcommand{\XTilde}{\widetilde{X}}
\newcommand{\YTilde}{\widetilde{Y}}

\renewcommand{\citep}[1]{[\citeauthor{#1} (\citeyear{#1})]}
\renewcommand{\cite}{\citet}
\newcommand{\eqref}[1]{\textup{(\ref{#1})}}
\newcommand{\limitScaling}{\mathfrak{s}}
\makeatother

\begin{document}
\begin{frontmatter}

\title{High-dimensionality effects in the Markowitz problem and other
quadratic programs with linear~constraints:
Risk underestimation\thanksref{T1}}
\runtitle{High-dimensional quadratic programs}
\thankstext{T1}{Supported by the France--Berkeley Fund, a Sloan
research Fellowship and NSF Grants DMS-06-05169 and DMS-08-47647 (CAREER).}

\begin{aug}
\author{\fnms{Noureddine} \snm{El Karoui}\corref{}\ead[label=e1]{nkaroui@stat.berkeley.edu}}
\runauthor{N. El Karoui}
\affiliation{University of California, Berkeley}
\address{Department of Statistics\\ 367 Evans Hall\\University of
California, Berkeley\\Berkeley, California 94720-3860\\USA\\
\printead{e1}} 
\end{aug}

\received{\smonth{8} \syear{2009}}
\revised{\smonth{1} \syear{2010}}

%
\begin{abstract}
We first study the properties of solutions of quadratic programs with
linear equality constraints whose parameters are estimated from data in
the high-dimensional setting where $p$, the number of variables in the
problem, is of the same order of magnitude as $n$, the number of
observations used to estimate the parameters. The Markowitz problem in
Finance is a subcase of our study. Assuming normality and independence
of the observations we relate the efficient frontier computed
empirically to the ``true'' efficient frontier. Our computations show
that there is a separation of the errors induced by estimating the mean
of the observations and estimating the covariance matrix. In
particular, the price paid for estimating the covariance matrix is an
underestimation of the variance by a factor roughly equal to $1-p/n$.
Therefore the risk of the optimal population solution is underestimated
when we estimate it by solving a similar quadratic program with
estimated parameters.

We also characterize the statistical behavior of linear functionals of
the empirical optimal vector and show that they are biased estimators
of the corresponding population quantities.

We investigate the robustness of our Gaussian results by extending the
study to certain elliptical models and models where our $n$
observations are correlated (in ``time''). We show a lack of robustness
of the Gaussian results, but are still able to get results concerning
first order properties of the quantities of interest, even in the case
of relatively heavy-tailed data (we require two moments). Risk
underestimation is still present in the elliptical case and more
pronounced than in the Gaussian case.

We discuss properties of the nonparametric and parametric bootstrap in
this context. We show several results, including the interesting fact
that standard applications of the bootstrap generally yield
inconsistent estimates of bias.

We propose some strategies to correct these problems and practically
validate them in some simulations. Throughout this paper, we will
assume that $p$, $n$ and $n-p$ tend to infinity, and $p<n$.

Finally, we extend our study to the case of problems with more general
linear constraints, including, in particular, inequality constraints.
\end{abstract}

%

\end{frontmatter}

\section{Introduction}
Many statistical estimation problems are now formulated, implicitly or
explicitly, as solutions of certain optimization problems. Naturally,
the parameters of these problems tend to be estimated from data and it
is therefore important that we understand the relationship between the
solutions of two types of optimization problems: those which use the
population parameters and those which use the estimated parameters.
This question is particularly relevant in high-dimensional inference
where one suspects that the differences between the two solutions might
be considerable. The aim of this paper is to contribute to this
understanding by focusing on quadratic programs with linear
constraints. An important example of such a program where our questions
are very natural is the celebrated Markowitz optimization problem in
Finance which will serve as a supporting example throughout the paper.

The Markowitz problem \citep{Markowitz1952} is a classic portfolio
optimization problem in Finance, where investors choose to invest
according to the following framework: one picks assets in such a way
that the portfolio guarantees a certain level of expected returns but
minimizes the ``risk'' associated with them. In the standard framework,
this risk is measured the variance of the portfolio.

Markowitz's paper was highly influential and much work has followed. It
is now part of the standard textbook literature on these issues [\cite
{Ruppert06FinanceBook}, \cite{CampbellLoMacKinlay}].
Let us recall the setup of the Markowitz problem.
\begin{itemize}
\item We have the opportunity to invest in $p$ assets, $A_1,\ldots,A_p$.
\item In the ideal situation, the mean returns are known and
represented by a $p$-dimensional vector, $\mu$.
\item Also, the covariance between the returns is known; we denote it
by $\Sigma$.
\item We want to create a portfolio, with guaranteed mean return $\mu
_P$, and minimize its risk, as measured by variance.
\item The question is how should items be weighted in portfolio? What
are weights~$w$?
\end{itemize}
We note that $\Sigma$ is positive semi-definite and hence is in
particular symmetric.
In the ideal (or population) solution, the covariance and the mean are
known. The mathematical formulation is then the following simple
quadratic program. We wish to find the weights $w$ that solve the
following problem:
\[
\cases{
\min\frac{1}{2} w\trsp\Sigma w, \cr
w\trsp\mu= \mu_P ,\cr\displaystyle
w\trsp\ebold=1.
}
\]
Here $\ebold$ is a $p$-dimensional vector with 1 in every entry. If
$\Sigma$ is invertible, the solution is known explicitly (see Section
\ref{Sec:GeneralQP}). If we call $\wOpt$ the solution of this
problem, the curve $\wOpt\trsp\Sigma\wOpt$, seen as a function of
$\mu_P$, is called the \textit{efficient frontier}.

Of course, in practice, we do not know $\mu$ and $\Sigma$ and we need
to estimate them. An interesting question is therefore to know what
happens in the Markowitz problem when we replace population quantities
by corresponding estimators.

Naturally, we can ask a similar question for general quadratic programs
with linear constraints [see below or \cite{boydvandenberghe04} for a
definition], the Markowitz problem being a particular instance of such
a problem. This paper provides an answer to these questions under
certain distributional assumptions on the data. Hence our paper is
really about the impact of estimation error on certain high-dimensional
$M$-estimation problems.

It has been observed by many that there are problems in practice when
replacing population quantities by standard estimators [see \cite
{LaiXing08FinanceBook}, Section~3.5], and alternatives have been
proposed. A famous one is the Black--Litterman model [\cite
{BlackLittermanGS90}, \cite{MeucciBook} and, e.g., \cite
{MeucciEnhancingBlackLitterman08}]. Adjustments to the standard
estimators have also been proposed:  \cite
{LedoitWolf04}, partly motivated by portfolio optimization problems,
proposed to ``shrink'' the sample covariance matrix toward another
positive definite matrix (often the identity matrix properly scaled),
while   \cite{MichaudBook98} proposed to use the bootstrap and
to average bootstrap weights to find better-behaved weights for the
portfolio. As noted in \cite{LaiXing08FinanceBook}, there is a dearth
of theoretical studies regarding, in particular, the behavior of
bootstrap estimators.

An aspect of the problem that is of particular interest to us is the
study of large-dimensional portfolios (or quadratic programs with
linear constraints). To make matters clear, we focus on a portfolio
with $p=100$ assets. If we use a year of daily data to estimate $\Sigma
$, the covariance between the daily returns of the assets, we have
$n\simeq250$ observations at our disposal. In modern statistical
parlance, we are therefore in a ``large $n$, large $p$'' setting, and we
know from random matrix theory that $\SigmaHat$ the sample covariance
matrix is a poor estimator of $\Sigma$, especially when it comes to
spectral properties of~$\Sigma$. There is now a developing statistical
literature on properties of sample covariance matrices when $n$ and $p$
are both large, and it is now understood that, though $\SigmaHat$ is
unbiased for $\Sigma$, the eigenvalues and eigenvectors of $\SigmaHat
$ behave very differently from those of $\Sigma$. We refer the
interested reader to \cite{imj}, \citeauthor{nekGencov} (\citeyear{nekGencov,nekSparseMatrices,nekCorrEllipD}),
\cite{bickellevinaThresh07}, \cite
{RothmanSpice08}  for a partial introduction to
these problems. We wish with this study to make clear that the ``large
$n$, large~$p$'' character of the problem has an important impact of the
empirical solution of the problem. By contrast, standard but thorough
discussions of these problems \citep{MeucciBook} give only a cursory
treatment of dimensionality issues (e.g., one page out of a whole book).

Another interesting aspect of this problem is that the high-dimensional
setting does not allow, by contrast to the classical ``small $p$, large
$n$'' setting, a perturbative approach to go through. In the ``small
$p$, large $n$'' setting, the paper \cite{JobsonKorkie80} is concerned,
in the Gaussian case, with issues similar to the ones we will be investigating.

The ``large $n$, large $p$'' setting is the one with which random matrix
theory is concerned---and the high-dimensional Markowitz problem has
therefore been of interest to random matrix theorists for some time
now. We note in particular the paper \cite{lalouxetalbis}, where a
random matrix-inspired (shrinkage) approach to improved estimation of
the sample covariance matrix is proposed in the context of the
Markowitz problem.

Let us now remind the reader of some basic facts of random matrix
theory that suggest that serious problems may arise if one solves
naively the high-dimensional Markowitz problem or other quadratic
programs with linear equality constraints. A key result in random
matrix theory is the Mar\v{c}enko--Pastur equation \citep{mp67}
which characterizes
the limiting distribution of the eigenvalues of the sample covariance
matrix and relates it to the spectral distribution of the population
covariance matrix. We give only in this introduction its simplest form
and refer the reader to \cite{mp67}, \cite{wachter78}, \cite
{silverstein95}, \cite{bai99} and, for example, \cite{nekCorrEllipD}
for a more thorough introduction and very recent developments, as well
as potential geometric and statistical limitations of the models
usually considered in random matrix theory.

\begin{figure}

\includegraphics{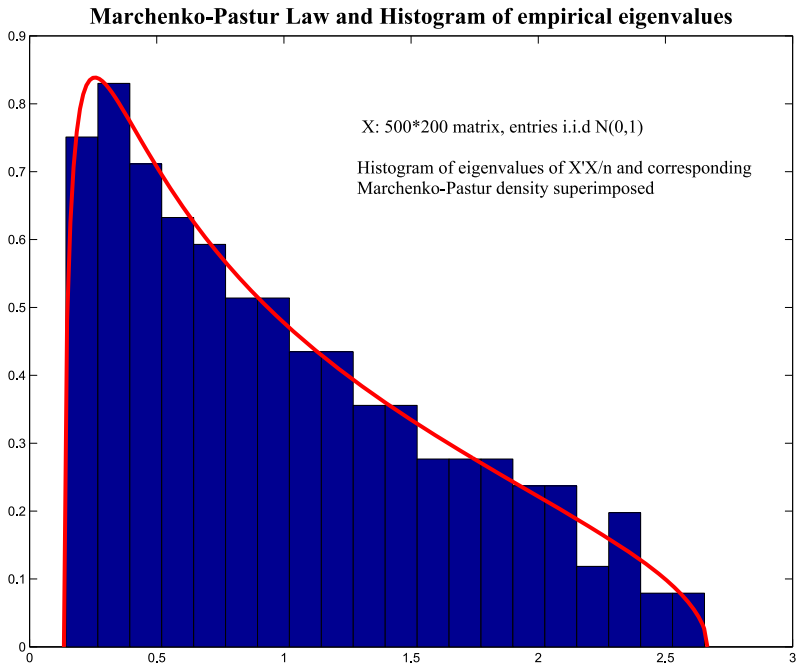}

\caption{Illustration of Mar\v{c}enko--Pastur-law,
$n=500$, $p=200$. The red curve is the density of the Mar\v
{c}enko--Pastur-law for $\rho
=2/5$. The simulation was done with i.i.d. Gaussian data. The histogram
is the histogram of eigenvalues of $X\trsp X/n$.}\label
{fig:IllustrationMapaLaw}
\end{figure}

In the simplest setting, we consider data $\{X_i\}_{i=1}^n$, which are
$p$-dimensional. In a financial context, these vectors would be vectors
of (log)-returns of assets, the portfolio consisting of $p$ assets. To
simplify the exposition, let us assume that the $X_i$'s are i.i.d. with
distribution ${\cal N}(0,\id_p)$. We call $X$ the $n\times p$ matrix
whose $i$th row is the vector $X_i$. Let us consider the sample
covariance matrix
\[
\SigmaHat=\frac{1}{n-1}(X-\bar{X})\trsp(X-\bar{X}) ,
\]
where $\bar{X}$ is a matrix whose rows are all equal to the column
mean of $X$. Now let us call $F_p$ the spectral distribution of
$\SigmaHat$, that is, the probability distribution that puts mass
$1/p$ at each of the $p$ eigenvalues of $\SigmaHat$. A graphical
representation of this probability distribution is naturally the
histogram of eigenvalues of $\SigmaHat$. A~consequence of the main
result of the very profound paper \cite{mp67} is that $F_p$, though a
random measure, is asymptotically nonrandom, and its limit, in the
sense of weak convergence of distributions, $F$ has a density (when
$p<n$) that can be computed. $F$~depends on $\rho=\lim_{n\tendsto
\infty} p/n$ in the following manner: if $p<n$, the density of $F$ is
\[
f_{\rho}(x)=\frac{1}{2\pi\rho}\frac{\sqrt
{(y_+-x)(x-y_-)}}{x}1_{y_-\leq x \leq y_+} ,
\]
where $y_+=(1+\sqrt{\rho})^2$ and $y_-=(1-\sqrt{\rho})^2$. Figure
\ref{fig:IllustrationMapaLaw} presents a graphical illustration of
this result.

%
%

What is striking about this result is that it implies that the largest
eigenvalue of $\Sigma$, $\lambda_1$, will be overestimated by $l_1$
the largest eigenvalue of $\SigmaHat$. Also, the smallest eigenvalue
of $\Sigma$, $\lambda_p$, will be underestimated by the smallest
eigenvalue of $\SigmaHat$, $l_p$. As a matter of fact, in the model
described above, $\Sigma$ has all its eigenvalues equal to 1, so
$\lambda_1(\Sigma)=\lambda_p(\Sigma)=1$, while $l_1$ will
asymptotically be larger or equal to $(1+\sqrt{\rho})^2$ and $l_p$
smaller or equal to $(1-\sqrt{\rho})^2$ (in the Gaussian case and
several others, $l_1$ and $l_p$ converge to those limits). We note that
the result of \cite{mp67} is not limited to the case where $\Sigma$
is identity, as presented here, but holds for general covariance~$\Sigma$ ($F_p$ has of course a different limit then).

Perhaps more concretely, let us consider a projection of the data along
a vector~$v$, with $\| v\|_2=1$, where $\| v\|_2$ is the
Euclidian norm of $v$. Here it is clear that, if $X\sim{\cal N}(0,\id
_p)$, $\operatorname{var}(v\trsp X)=1$, for all $v$, since $v\trsp X
\sim{\cal
N}(0,1)$. However, if we do not know $\Sigma$ and estimate it by
$\SigmaHat$, a naive (and wrong) reasoning suggests that we can find
direction of lower variance than 1, namely those corresponding to
eigenvectors of $\SigmaHat$ associated with eigenvalues that are less
than 1. In particular, if $v_p$ is the eigenvector associated with
$l_p$, the smallest eigenvalue of $\SigmaHat$, by naively estimating,
for $X$ independent of $\{X_i\}_{i=1}^n$, the variance in the direction
of $v_p$, $\operatorname{var}(v_p\trsp X)$, by the empirical version
$v_p\trsp
\SigmaHat v_p$, one would commit a severe mistake: the variance in any
direction is 1, but it would be estimated by something roughly equal to
$(1-\sqrt{p/n})^2$ in the direction of $v_p$.

In a portfolio optimization context, this suggests that by using
standard estimators, such as the sample covariance matrix, when solving
the high-dimensional Markowitz problem, one might underestimate the
variance of certain portfolios (or ``optimal'' vectors of weights). As a
matter of fact, in the previous toy example, thinking (wrongly) that
there is low variance in the direction $v_p$, one might (numerically)
``load'' this direction more than warranted, given that the true
variance is the same in all directions.

This simple argument suggests that severe problems might arise in the
high-dimensional Markowitz problem and other quadratic programs with
linear constraints, and in particular, risk might be underestimated.
While this heuristic argument is probably clear to specialists of
random matrix theory, the problem had not been investigated at a
mathematical level of rigor in that literature before this paper was
submitted [the paper \cite{BaiLiuWongMarko09} has appeared while this
paper was being refereed. It is concerned with different models than
the ones we will be investigating and our results do not overlap]. It
has received some attention at a physical level of rigor [see, e.g.,
\cite{PafkaKondor03}, where the authors treat only the Gaussian case,
and do not investigate the effect of the mean, which as we show below
creates problems of its own].
In this paper, we propose a theoretical analysis of the problem in a
Gaussian and elliptical framework for general quadratic programs with
linear constraints, one of them involving the parameter~$\mu$. Our
results and contributions are several-fold. We relate the empirical
efficient frontier to the theoretical efficient frontier that is key to
the Markowitz theory, in a variety of theoretical settings. We show
that the empirical frontier generally yields an underestimation of the
risk of the portfolio and that Gaussian analysis gives an
over-optimistic view of this problem. We show that the expected returns
of the naive ``optimal'' portfolio are poorly estimated by $\mu_P$. We
argue that the bootstrap will not solve the problems we are pointing
out here. Beside new formulas, we also provide robust estimators of the
various quantities we are interested in.

The paper is divided into four main parts and a conclusion. In
Section~\ref{Sec:GeneralQP}, to make the paper self-contained, we
discuss the solution of
quadratic problems with linear equality constraints---a focus of this
paper. In Section~\ref{sec:GaussianCase}, we study the impact of
parameter estimation on the
solution of these problems when the observed data is i.i.d. Gaussian
and obtain some exact distributional results for fixed $p$ and $n$. In
Section~\ref{sec:EllipticalCase}, we obtain results in the case where
the data is elliptically
distributed. This allows us also to understand the impact of
correlation between observations in the Gaussian case and to get
information about the behavior of the nonparametric bootstrap. In
Section~\ref{sec:ComparisonGaussianElliptical}, we apply the results
of Section~\ref{sec:EllipticalCase} to the quadratic programs
at hand and compare the elliptical and the Gaussian cases. We show,
among other things, that the Gaussian results are not robust in the
class of elliptical distribution. In particular, two models may yield
the same $\mu$ and $\Sigma$ but can have very different empirical
behavior. In Section~\ref{sec:ComparisonGaussianElliptical}, we also
propose various schemes to correct the
problems we highlight (see pages \pageref{fig:ReturnsCorrectionSmall},
\pageref{fig:ReturnsCorrection} and \pageref{fig:FrontiersCorrection}
for pictures) and study more general problems with linear constraints
(see Section~\ref{Subsec:InequalityConstraints}). The conclusion
summarizes our findings and the \hyperref[appm]{Appendix} contains various facts and
proofs that did not naturally flow in the main text or were better
highlighted by being stated separately.

Several times in the paper $\SigmaHat^{-1}$ and $\Sigma^{-1}$ will
appear. Unless otherwise noted, when taking the inverse of a population
matrix, we implicitly assume that it exists. The question of existence
of inverse of sample covariance matrices is well understood in the
statistics literature. Because our models will have a component with a
continuous distribution, there are essentially no existence problems
(unless we explicitly mention and treat them) as proofs similar to
standard ones found in textbooks [e.g., \cite{anderson03}] would show.
Hence, we do not belabor this point any further in the rest of the
paper as our focus is on things other than rather well-understood
technical details, and the paper is already a bit long.

Finally, let us mention that while the Finance motivation for our study
is important to us, we treat the problem in this paper as a
high-dimensional $M$-estimation question (which we think has practical
relevance). We will not introduce particular modelization assumptions
which might be relevant for practitioners of Finance but might make the
paper less relevant in other fields. A companion paper \citep
{nekMarkoRealizedRisk} deals with more ``financial'' issues and the
important question of the realized risk of portfolios that are
``plug-in'' solutions of the Markowitz problem.

\section{Quadratic programs with linear equality constraints}\label
{Sec:GeneralQP}
We discuss here the properties of the solution of quadratic programs
with linear equality constraints as they lay the foundations for our
analysis of similar problems involving estimated parameters (and of
problems with inequality constraints). We included this section for the
convenience of the reader to make the paper as self-contained as possible.

The problem we want to solve is the following:
\renewcommand{\theequation}{QP-eqc}
\begin{equation}\label{eq:GeneralQP}
\cases{\displaystyle
\min_{w \in\mathbb{R}^p} \frac{1}{2} w\trsp\Sigma w, \cr
\displaystyle
w\trsp v_i = u_i ,\qquad1\leq i \leq k .
}
\end{equation}
Here $\Sigma$ is a positive definite matrix of size $p\times p$, $v_i
\in\mathbb{R}^p$ and $u_i\in\mathbb{R}$. We have the following theorem:
\begin{theorem}\label{thm:SolnGeneralQP}
Let us call $V$ the $p\times k$ matrix whose $i$th column is $v_i$, $U$
the~$k$ dimensional vector whose $i$th entry is $u_i$ and $M$ the
$k\times k$ matrix
\[
M=V\trsp\Sigma^{-1} V .
\]
We assume that the $v_i$'s are such that $M$ is invertible.
The solution of the quadratic program with linear equality constraints
\eqref{eq:GeneralQP} is achieved for
\[
\wOpt=\Sigma^{-1}V M^{-1} U,
\]
and we have
\[
\wOpt\trsp\Sigma\wOpt=U\trsp M^{-1} U .
\]
\end{theorem}

\begin{pf}
Let us call $\lambda$ a $k$ dimensional vector of Lagrange
multipliers. The Lagrangian function is, in matrix notation,
\[
L(w,\lambda)=\frac{w\trsp\Sigma w}{2}-\lambda\trsp(V\trsp w-U) .
\]
This is clearly a (strictly) convex function in $w$, since $\Sigma$ is
positive definite by assumption. We have
\[
\frac{\partial L}{\partial w}=\Sigma w - V\lambda.
\]
So $\wOpt=\Sigma^{-1}V\lambda$. Now we know that $U=V\trsp\wOpt$.
So $U=V\trsp\Sigma^{-1} V \lambda=M\lambda$. Therefore,
\[
\wOpt=\Sigma^{-1}V M^{-1} U .
\]
We deduce immediately that
\[
\wOpt\trsp\Sigma\wOpt=U\trsp M^{-1} U .
\]
\upqed
\end{pf}

We now turn to another result which will prove to be useful later. It
gives a compact representation of linear combinations of the weights
of
the optimal solution, and we will rely heavily on it in particular in
the case of Gaussian data.\vadjust{\goodbreak}
\begin{lemmaNumThm}\label{lemma:repLinearCombOptimalWeights}
Let us consider $\wOpt$ the solution of the optimization problem
\eqref{eq:GeneralQP}. Let $\gamma$ be a vector in $\mathbb{R}^p$.
Let us call ${\cal M}$ the $(k+1)\times(k+1)$ matrix that is written
in block form
\[
{\cal M}= \pmatrix{
V\trsp\Sigma^{-1} V & V\trsp\Sigma^{-1} \gamma\cr
\gamma\trsp\Sigma^{-1} V & \gamma\trsp\Sigma^{-1}\gamma
}
.
\]
\renewcommand{\theequation}{\arabic{equation}}\setcounter
{equation}{0}Assume that ${\cal M}$ is invertible. Then
\begin{equation}\label{eq:RepLinCombiWeights}
\gamma\trsp\wOpt=-\frac{1}{{\cal M}^{-1}_{k+1,k+1}} (U
\trsp0 ){\cal M}^{-1} \myVector{0_k}{1}.
\end{equation}

\end{lemmaNumThm}

\begin{pf}
The proof is a consequence of the results discussed in the \hyperref[appm]{Appen-}
\hyperref[appm]{dix}
concerning inverses of partitioned matrices [see Section \ref
{subsec:ClassicalLinearAlgebra} and equation \eqref
{eq:Rep21BlockInInverse} there]. Let us write
\[
{\cal M}= \pmatrix{
{\cal M}_{11} & {\cal M}_{12}\cr
{\cal M}_{21} & {\cal M}_{22}
} ,
\]
where ${\cal M}_{11}$ is $k\times k$, ${\cal M}_{12}$ is naturally
$k\times1$ and ${\cal M}_{22}$ is a scalar. With the same block
notation, we have
\[
{\cal M}^{-1}= \pmatrix{
{\cal M}^{11} & {\cal M}^{12}\cr
{\cal M}^{21} & {\cal M}^{22}
}
.
\]
Then, we know [see equation \eqref{eq:Rep21BlockInInverse}] that $
{\cal M}^{12}=-{\cal M}_{11}^{-1} {\cal M}_{12} {\cal M}^{22} ,$
but since ${\cal M}^{22}$ is a scalar, equal to ${\cal
M}^{-1}(k+1,k+1)$, we have
\[
{\cal M}_{11}^{-1} {\cal M}_{12}=-{\cal M}^{12}/{\cal M}^{22} .
\]
Now ${\cal M}_{11}^{-1} {\cal M}_{12}=(V\trsp\Sigma^{-1} V)^{-1}
V\trsp\Sigma^{-1} \gamma$, so $U\trsp{\cal M}_{11}^{-1} {\cal
M}_{12}=\wOpt\trsp\gamma$. Hence,
\[
\wOpt\trsp\gamma=-\frac{1}{{\cal M}^{22}}(U\trsp0) {\cal
M}^{-1} \myVector{0_k}{1} .
\]
\upqed
\end{pf}

We note that here $({\cal M}^{22})^{-1}=\gamma\trsp\Sigma^{-1}\gamma
-\gamma\trsp\Sigma^{-1}V M^{-1} V\trsp\Sigma^{-1} \gamma$, as an
application of equation \eqref{eq:Rep22BlockInInverse} clearly shows.

\section{QP with equality constraints: Impact of parameter estimation
in the Gaussian case} \label{sec:GaussianCase}
From now on, we will assume that we are in the high-dimensional setting
where $p$ and $n$ go to infinity. Our study will be divided into two.
We will first consider the Gaussian setting (in this section) and then
study an elliptical distribution setting (in Section \ref
{sec:EllipticalCase}). (We note that for the Markowitz problem, the
assumption of Gaussianity would be satisfied if we worked under
Black--Scholes diffusion assumptions for our assets and were
considering log-returns as our observations.) Interestingly, we will
show that the results are not robust against the assumption of
Gaussianity, which is not (so) surprising in light of recent random
matrix results [see \cite{nekCorrEllipD}]. We will also show that
understanding the elliptical setting allows us to understand the impact
of correlation between observations and to discuss bootstrap-related
ideas. In particular, we will see that various problems arise with the
bootstrap in high-dimension and that the results change when one deals
with observations that are correlated (in time) or not.

We also address similar questions concerning inequality constrained
problems in Section~\ref{Subsec:InequalityConstraints}.

Before we proceed, we need to set up some notations: we call $\ebold$
the $p$-dimensional vector whose entries are all equal to 1. We call
$V$, as above, the matrix containing all of our constraint vectors,
which we may have to estimate (for instance, if $v_i=\mu$ for a
certain $i$). We call $\widehat{V}$ the matrix of estimated constraint vectors.

The template question for all our investigations will be the following
(Marko\-witz) question: what can be said of the statistical properties of
the solution of
\[
\cases{\displaystyle
\min_{w \in\mathbb{R}^p} w\trsp\SigmaHat w, \cr\displaystyle
w\trsp\muHat= \mu_P ,\cr\displaystyle
w\trsp\ebold=1
}
\]
compared to the solution of the population version
\[
\cases{\displaystyle
\min_{w \in\mathbb{R}^p} w\trsp\Sigma w, \cr\displaystyle
w\trsp\mu= \mu_P ,\cr\displaystyle
w\trsp\ebold=1?
}
\]

We will solve the problem at a much greater degree of generality, by
considering first quadratic programs with linear equality constraints
(see Section~\ref{Subsec:InequalityConstraints} for inequality
constraints) and comparing the solutions of
\renewcommand{\theequation}{QP-eqc-Emp}
\begin{equation}\label{eq:QP-eqc-Emp}
\cases{\displaystyle
\min_{w \in\mathbb{R}^p} w\trsp\SigmaHat w, \cr\displaystyle
w\trsp v_i = u_i , \qquad1\leq i \leq k-1 ,\cr\displaystyle
w\trsp\muHat= u_k
}
\end{equation}
and
\renewcommand{\theequation}{QP-eqc-Pop}
\begin{equation}\label{eq:QP-eqc-Pop}
\cases{\displaystyle
\min_{w \in\mathbb{R}^p} w\trsp\Sigma w, \cr\displaystyle
w\trsp v_i = u_i ,\qquad 1\leq i \leq k-1 ,\cr\displaystyle
w\trsp\mu= u_k.
}
\end{equation}
Here $\SigmaHat$ and $\muHat$ will be estimated from the data. We
call $\wEmp$ the vector that yields a solution of problem \eqref
{eq:QP-eqc-Emp} and $\wTheo$ the vector that yields a solution of
problem \eqref{eq:QP-eqc-Pop}.

We call $\VHat$ the $p\times k$ matrix containing $\{v_i\}
_{i=1}^{k-1}$ and $\muHat$, and $V$ its population counterpart, which
contains $\{v_i\}_{i=1}^{k-1}$ and $\mu$. We assume that $\{v_i\}
_{i=1}^{k-1}$ are deterministic and known (just like the vector $\ebold
$ in the Markowitz problem).\vadjust{\goodbreak} In our analysis, $k$ will be held fixed.
(The $k$th column of $\VHat$ will contain $\muHat$ in general or our
estimator of $\mu$.)

As should be clear from Theorem~\ref{thm:SolnGeneralQP}, the
properties of the entries of the matrix $\widehat{V}\trsp\SigmaHat
^{-1}\widehat{V}$ as compared to those of the matrix $V\trsp\Sigma
^{-1}V$ will be key to our understanding of this question. In what
follows, we assume that the vectors $\hat{v}_i$ are either
deterministic or equal to $\muHat$. The extension to linear
combinations of a deterministic vector and $\muHat$ is
straightforward. We also note that in the Gaussian case, we could just
assume that the $\hat{v}_i$ are (deterministic) functions of $\muHat$
(because $\muHat$ and $\SigmaHat$ are independent in this case). On
the other hand, the vector $U$ is assumed to be deterministic.

Before we proceed, let us mention that after our study was completed,
we learned of similar results (restricted to the Markowitz case and not
dealing with general quadratic programs with linear equality
constraints) by \cite{KanAndSmith08}. We stress the fact that our work
was independent of theirs and is more general which is why it is
included in the paper.

\subsection{Efficient frontier problems}\label{subsec:GaussianCase}
We first study questions concerning the efficient frontier and then
turn to information we can get about linear functionals of the
empirical weights.
\begin{theorem}\label{thm:GalQPGaussianCase}
Let us assume that we observe data $X_i\iid{\cal N}(\mu,\Sigma)$,
for $i=1,\ldots,n$. Here $\Sigma$ is $p\times p$ and $p<n$. Suppose
we estimate $\Sigma$ with the sample covariance matrix $\SigmaHat$,
and $\mu$ with the sample mean $\muHat$. Suppose we wish to solve the problem
\renewcommand{\theequation}{QP-eqc-Pop}
\begin{equation}
\cases{\displaystyle
\min_{w \in\mathbb{R}^p} w\trsp\Sigma w , \cr\displaystyle
w\trsp v_j = u_j ,\qquad1\leq j \leq k .
}
\end{equation}
where $u_j$ are deterministic, $v_j$ are deterministic and given for
$j<k$ and $v_k=\mu$.
Assume that we use as a proxy for the previous problem the empirical
version with plugged-in parameters.
Let us consider the solution of the problem
\renewcommand{\theequation}{QP-eqc-Emp}
\begin{equation}
\cases{\displaystyle
\min_{w \in\mathbb{R}^p} w\trsp\SigmaHat w, \cr\displaystyle
w\trsp\hat{v}_j = u_j ,\qquad1\leq j \leq k .
}
\end{equation}
Now $\hat{v}_j=v_j$ for $j<k$ and $\hat{v}_k=f(\muHat)$, for a given
deterministic function $f$.
Let us call $\wEmp$ the corresponding ``weight'' vector. The plug-in
estimate of $w\trsp\Sigma w $ is $\wEmp\trsp\SigmaHat\wEmp$.
Let us call $\wOracle$ the optimal solution of the quadratic program
obtained under the assumption that $\Sigma$ is given, but $\mu$ is
not and is estimated by $f(\muHat)$. Finally, we assume that $n-1-p+k>0$.

Then we have
\renewcommand{\theequation}{\arabic{equation}}\setcounter{equation}{1}
\begin{equation}\label{eq:OracleRepRiskUnderEstimationGeneralQP}
\wEmp\trsp\SigmaHat\wEmp= \wOracle\trsp\Sigma\wOracle\frac
{\chi^2_{n-1-p+k}}{n-1} ,
\end{equation}
where $\wOracle\trsp\Sigma\wOracle$ is random (because $\muHat$
is) but is statistically independent of $\chi^2_{n-1-p+k}$. Also,
\[
\wOracle\trsp\Sigma\wOracle=U\trsp(\VHat\trsp\Sigma
^{-1}\VHat)^{-1} U .
\]
\end{theorem}

The previous theorem means that the cost of not knowing the covariance
matrix and estimating it is the apparition of the $\frac{\chi
^2_{n-1-p+k}}{n-1}$. In the high-dimensional setting when $p$ and $n$
are of the same order of magnitude and $n-p$ is large, this terms is
approximately $1-(p-k)/(n-1)$. Hence, the theorem quantifies the random
matrix intuition that having to estimate the high-dimensional
covariance matrix at stake here leads to risk \textit
{underestimation}, by the factor $1-(p-k)/(n-1)$. In other words, using
plug-in procedures leads to over-optimistic conclusions in this situation.

We also note that the previous theorem shows that, in the Gaussian
setting under study here, the effect of estimating the mean and the
covariance on the solution of the quadratic program are ``separable'':
the effect of the mean estimation is in the oracle term, while the
effect of estimating the covariance is in the $\chi^2_{n-p-1+k}/(n-1)$
term. To show risk underestimation, it will therefore be necessary to
relate $\wOracle\trsp\Sigma\wOracle$ to $\wTheo\trsp\Sigma\wTheo
$. We do it in Proposition \ref
{prop:RiskUnderEstimationByMeanEstimation} but first give a proof of
Theorem~\ref{thm:GalQPGaussianCase}.

\begin{pf*}{Proof of Theorem~\ref{thm:GalQPGaussianCase}}
The crux of the proof is the following result, which is well known by
statisticians, concerning (essentially) blocks of the inverse of a
Wishart matrix: if $S\sim\Wishart_p(\Sigma,m)$, that is, $S$ is a
$p\times p$ Wishart matrix with $m$ degree of freedoms and covariance
$\Sigma$, and $A$ is $p\times k$, deterministic matrix, then, when $m>p$,
\[
(A\trsp S^{-1} A)^{-1}\sim\Wishart_k\bigl((A\trsp\Sigma
^{-1}A)^{-1},m-p+k\bigr) .
\]
We refer to \citeauthor{EatonMultivariateStatBook83} [(\citeyear{EatonMultivariateStatBook83}), Proposition 8.9, page
312] for a proof, and to \citeauthor{mardiakentbibby} [(\citeyear{mardiakentbibby}), pages 70--73] for
related results.

Another important remark is the well-known fact that, in the situation
we are considering, $\muHat$ is ${\cal N}(\mu, \Sigma/n)$ and
independent of $\SigmaHat$. Finally, it is also well known that if
$S\sim\Wishart_p(\Sigma,m)$ and $U$ is a $p$-dimensional
deterministic vector, then $U\trsp S U=U\trsp\Sigma U \chi^2_m$.

Now $\SigmaHat\sim\Wishart_p(\Sigma,n-1)/(n-1)$. Therefore, since
$\widehat{V}$ is a function of $\muHat$, we have, by independence of
$\muHat$ and $\SigmaHat$,
\[
(\widehat{V}\trsp\SigmaHat^{-1} \widehat{V})^{-1} |
\muHat\sim\Wishart_k\bigl((\VHat\trsp\Sigma^{-1}\VHat
)^{-1},n-1-p+k\bigr)/(n-1) .
\]
Therefore,
\[
\frac{U\trsp(\widehat{V}\trsp\SigmaHat^{-1} \widehat
{V})^{-1} U }{U\trsp(\VHat\trsp\Sigma^{-1}\VHat)^{-1} U}\bigg|
\muHat\sim\frac{\chi^2_{n-p-1+k}}{n-1} .
\]
Because the right-hand side does not depend on $\muHat$, we have
established the independence of
\[
\frac{U\trsp(\widehat{V}\trsp\SigmaHat^{-1} \widehat{V})^{-1} U
}{U\trsp(\VHat\trsp\Sigma^{-1}\VHat)^{-1} U} \quad\mbox
{and}\quad\frac
{\chi^2_{n-p-1+k}}{n-1} .
\]
Hence, we conclude that
\[
U\trsp(\widehat{V}\trsp\SigmaHat^{-1} \widehat{V})^{-1} U = U\trsp
(\VHat\trsp\Sigma^{-1}\VHat)^{-1} U \frac{\chi^2_{n-p-1+k}}{n-1} ,
\]
and the two terms are independent. Now the term $U\trsp(\VHat\trsp
\Sigma^{-1}\VHat)^{-1} U$ is the estimate we would get for the
solution of problem \eqref{eq:QP-eqc-Pop}, if $\Sigma$ were known and
$\mu$ were estimated by $f(\muHat)$. In other words, it is the
``oracle'' solution described above.
\end{pf*}

\subsubsection{Some remarks on the oracle solution}
Theorem~\ref{thm:GalQPGaussianCase} sheds light on the separate
effects of mean and covariance estimation on the problem considered
above. To understand further the problem of risk estimation, we need to
better understand the role the estimation of the mean might play. This
is what we do now.

\begin{propositionNumThm}\label{prop:RiskUnderEstimationByMeanEstimation}
Suppose that the last column of $\VHat$ is $\muHat$. Let us call
$V_{-k}$ the $p\times k-1$ dimensional matrix whose $j$th column is
$v_j$, which are known deterministic vectors. Suppose that $M=V\trsp
\Sigma^{-1}V=\gO(1)$. Suppose further that $\lambda_k(V\trsp\Sigma
^{-1} V)\gg n^{-1/2}$, where $\lambda_k(S)$ is the smallest eigenvalue
of the $k\times k$ matrix~$S$.

Further, call $M=V\trsp\Sigma^{-1} V \in\mathbb{R}^{k\times k}$ and
call $e_i$ the canonical basis vectors in $\mathbb{R}^k$. Finally,
call $\alpha=\chi^2_p/n$.

Then, when $p/n\tendsto\rho\in(0,1)$, asymptotically,
\[
\wOracle\trsp\Sigma\wOracle= \wTheo\trsp\Sigma\wTheo-\alpha
\frac{(U\trsp M^{-1} e_k)^2}{1+\alpha e_k\trsp M^{-1} e_k}+\lo_P
(\wTheo\trsp\Sigma\wTheo) .
\]
\end{propositionNumThm}

Let us discuss a little bit this result before we provide a proof. In
the asymptotics we have in mind and are considering, $p/n\tendsto\rho
\in(0,1)$ and therefore $\alpha\simeq p/n+\gO(n^{-1/2})$. So if
$\delta_n=(U\trsp M^{-1}e_k)^2/(1+p/n e_k\trsp M^{-1}e_k)$, when the
above analysis applies, the impact of the estimation of $\mu$ by
$\muHat$ will be risk underestimation, just as is the case for the
case of the covariance matrix. Here, we can also quantify the impact of
this estimation of $\mu$ by $\muHat$: it leads to risk
underestimation by the amount $\alpha\delta_n$.

\begin{pf*}{Proof of Proposition~\ref
{prop:RiskUnderEstimationByMeanEstimation}}
Let us write $\muHat=\mu+e$, where $e\sim{\cal N}(0,\Sigma/n)$.
Clearly, $e=n^{-1/2}\Sigma^{1/2}Z$, where $Z$ is ${\cal N}(0,\id_p)$.
We have, using block notations,
\[
\VHat\trsp\Sigma^{-1} \VHat=V\trsp\Sigma^{-1} V+
\pmatrix{ 0 & 0 \cr0 & e\trsp\Sigma^{-1} e
}
+
\pmatrix{ 0 & V_{-k}\trsp\Sigma^{-1} e\cr e\trsp\Sigma^{-1}
V_{-k}& 2 \mu\trsp\Sigma^{-1} e
}
.
\]
Replacing $e$ by its value, we have $\mu\trsp\Sigma^{-1}e \sim{\cal
N}(0,\mu\trsp\Sigma^{-1} \mu/n)$. By the same token, we can also
get that
\[
V_{-k}\trsp\Sigma^{-1} e = \frac{1}{\sqrt{n}} V_{-k}\trsp\Sigma
^{-1/2} Z\sim{\cal N} \biggl(0,\frac{V_{-k}\trsp\Sigma^{-1}
V_{-k}}{n} \biggr) .
\]
Our assumption that $V\trsp\Sigma^{-1}V=\gO(1)$ implies that $\mu
\trsp\Sigma^{-1} \mu=\gO(1)$ and $V_{-k}\trsp\times\Sigma^{-1}
V_{-k}=\gO(1)$. Therefore,
\[
\pmatrix{ 0 & V_{-k}\trsp\Sigma^{-1} e\cr e\trsp\Sigma^{-1}
V_{-k}& 2 \mu\trsp\Sigma^{-1} e
}
=\gO_P \biggl(\frac{1}{\sqrt{n}} \biggr) .
\]
Hence, since $e\trsp\Sigma^{-1} e =Z\trsp Z /n=\alpha$,
\[
\VHat\trsp\Sigma^{-1} \VHat=V\trsp\Sigma^{-1} V+\alpha e_k
e_k\trsp+\gO_P(n^{-1/2}) .
\]
Our assumptions guarantee that $\lambda_k(V\trsp\Sigma^{-1} V)\gg
n^{-1/2}$, and therefore\break $\lambda_k(V\trsp\Sigma^{-1}\times  V+\alpha
e_ke_k\trsp)\gg n^{-1/2}$.
In other respects, let $A$ be a matrix such that $\lambda_p(A)\gg
n^{-1/2}$ and $E$ be a matrix such that $E=\gO(n^{-1/2})$. Recall that
for symmetric matrices, $\lambda_p(A+E)\geq\lambda_p(A)+\lambda
_p(E)$ [see, e.g., Weyl's theorem, \cite{hornjohnson94}, page 185]. So
in this situation, $(A+E)^{-1}=\lo(n^{1/2})$. Let us now consider the
implications of this remark on the difference of $(A+E)^{-1}$ and $A^{-1}$.
We claim that $(A+E)^{-1}=A^{-1}+\lo(A^{-1})$. By the first resolvent
identity, $(A+E)^{-1}=A^{-1}-(A+E)^{-1}E A^{-1}$; our previous remark
implies that $\sigma_1[(A+E)^{-1}E]=\lo(1)$ and the result follows.
Applying the results of this discussion to $A=V\trsp\Sigma^{-1}
V+\alpha e_k e_k\trsp$ and $A+E=\VHat\trsp\Sigma^{-1} \VHat$, we have
\[
\VHat\trsp\Sigma^{-1} \VHat=(V\trsp\Sigma^{-1} V+\alpha e_k
e_k\trsp)^{-1}+\lo_P\bigl((V\trsp\Sigma^{-1} V+\alpha e_k e_k\trsp
)^{-1}\bigr) .
\]
We can now use well-known results concerning inverses of rank-1
perturbation of matrices, namely
\[
(V\trsp\Sigma^{-1} V+\alpha e_k e_k\trsp)^{-1}=(M+\alpha e_ke_k\trsp
)^{-1}=M^{-1}-\alpha\frac{M^{-1}e_ke_k\trsp M^{-1}}{1+\alpha e_k\trsp
M^{-1} e_k} .
\]
This allows us to conclude that
\[
U\trsp(\VHat\trsp\Sigma^{-1} \VHat)^{-1} U= U\trsp M^{-1}U -\alpha
\frac{(U\trsp M^{-1}e_k)^2}{1+\alpha e_k\trsp M^{-1}e_k}+\lo_P(U\trsp
M^{-1}U) .
\]
This is the result announced in the theorem and the proof is complete.
\end{pf*}

We can now combine the results of Theorem~\ref{thm:GalQPGaussianCase}
and Proposition~\ref{prop:RiskUnderEstimationByMeanEstimation} to
obtain the following corollary.
\begin{corollaryNumThm}\label{coro:GaussianCaseFinalApproxEffFrontier}
We assume that the assumptions of Theorem~\ref{thm:GalQPGaussianCase}
and Proposition~\ref{prop:RiskUnderEstimationByMeanEstimation} hold
and that $p/n$ has a finite\vadjust{\goodbreak} nonzero limit, as $n\tendsto\infty$, and
$n-p$ tends to infinity. Then we have
\renewcommand{\theequation}{\arabic{equation}}\setcounter{equation}{2}
\begin{eqnarray}\label{eq:RiskUnderEstimationGalQPGaussianCase}
\wEmp\trsp\SigmaHat\wEmp&=& \biggl(1-\frac{p-k}{n-1} \biggr)
\biggl(\wTheo\trsp\Sigma\wTheo-\frac{p}{n} \frac{(U\trsp M^{-1}
e_k)^2}{1+ ({p}/{n}) e_k\trsp M^{-1} e_k} \biggr)
\nonumber
\\[-8pt]
\\[-8pt]
&&{}+\lo_P (\wTheo
\trsp\Sigma\wTheo\vee n^{-1/2} ),\nonumber
\end{eqnarray}
where $M$ is the population quantity $M=V\trsp\Sigma^{-1}V$.
\end{corollaryNumThm}

The corollary shows that the effects of both covariance and mean
estimation are to underestimate the risk, and the empirical frontier is
asymptotically deterministic.

\subsection{On the optimal weights}
Our matrix characterization of the empirical optimal weights (Lemma
\ref{lemma:repLinearCombOptimalWeights}) allows us to give a precise
characterization of the statistical properties of linear functionals of
these weights. We give here some exact results, concerning
distributions and expectations of those functionals. A longer
discussion, including robustness and more detailed bias issues can be
found in Section~\ref{sec:ComparisonGaussianElliptical}.

\begin{propositionNumThm}\label{prop:exactResultsWeightsGaussianCase}
Assume that the assumptions of Theorem~\ref{thm:GalQPGaussianCase}
hold and in particular $X_i$ are i.i.d. ${\cal N}(\mu,\Sigma_p)$. Let
$\gamma$ be a fixed $n$-dimensional vector. Let us call $\VHat
_{\gamma}=(\VHat\gamma)$ the $p\times(k+1)$ matrix whose first
$k$ columns are those of $\VHat$. Let $\NHat_{\gamma}=(\VHat
_{\gamma}\trsp\Sigma^{-1} \VHat_{\gamma})^{-1}$ and $W_{\gamma}$
be a $(k+1)\times(k+1)$ matrix with distribution $\Wishart
_{k+1}(\NHat_{\gamma},n-p+k)$ (conditional on $\muHat$). Then,
\[
\gamma\trsp\wEmp| \muHat\equalInLaw-\frac{\sum
_{i=1}^k u_i W_{\gamma}(i,k+1)}{W_{\gamma}(k+1,k+1)} .
\]
In particular,
\[
\mathbf{E}(\gamma\trsp\wEmp|\muHat)=-\frac{\sum_{i=1}^k u_i
\NHat
_{\gamma}(i,k+1)}{\NHat_{\gamma}(k+1,k+1)} .
\]
\end{propositionNumThm}

We note, somewhat heuristically, that when $\mu$ is estimated by
$\muHat$, since $\muHat\sim{\cal N}(\mu, \Sigma/n)$, $\muHat\trsp
\Sigma^{-1}\muHat\simeq\mu\trsp\Sigma^{-1}\mu+ p/n$, when $p$,
$n$ and $n-p$ are all large (we refer again to Section \ref
{sec:ComparisonGaussianElliptical} for a more precise statement). Hence
$\NHat_{\gamma}$ is a not a consistent estimator of $N_{\gamma
}=(V_{\gamma}\trsp\Sigma^{-1} V_{\gamma})^{-1}$. As we will see in
Section~\ref{Sec:GaussEllip:Subsec:WeightsPortfolio} and as can be
expected from the previous proposition, this will also imply bias for
linear combinations of empirical optimal weights. We will show in
particular that returns are overestimated when using $\muHat$ as an
estimator for~$\mu$.

Another interesting aspect of the previous proposition is that it
allows us to understand the fluctuation behavior of $\gamma\trsp\wEmp
$ when $n-p+k$ is large: as a matter of fact, the limiting fluctuation
behavior of the entries of a (fixed-dimensional) Wishart matrix with
large number of degrees of freedom\vadjust{\goodbreak} is well known [see, e.g., \cite
{anderson03}, Theorem 3.4.4, page 87] and the $\delta$-method can be
applied to get the information---conditional on $\muHat$.

For instance, if we assume that, conditional on $\muHat$, the matrix
$\NHat_{\gamma}$ converges to a matrix $N^0_{\gamma}$, which
possibly depends on $\muHat$, we see that calling $\nu$ the last
column $W_{\gamma}/(n-p+k)$, $\nu$ is asymptotically normal (all
statements are conditional on $\muHat$), if $n-p+k$ goes to infinity
when $p$ and $n$ go to infinity. Furthermore we know the limiting
covariance of $\nu$ (after scaling by $\sqrt{n-p+k}$), using Theorem
3.4.4 in
\cite{anderson03}. Let us call it $\Gamma_0$ and let us call $\nu_0$
the limit of $\nu$---which we assume exists.

If we assume that $\nu_0(k+1)$ is not 0, Slutsky's lemma and the
$\delta$-method give us through simple computations that
\[
\sqrt{n-p+k} \biggl(\gamma\trsp\wEmp+\frac{\sum_{i=1}^k
u_i \nu_0(i)}{\nu_0(k+1)} \biggr) \Big| \muHat\WeakCv\frac
{1}{\nu_0(k+1)^2}{\cal N}(0,C\trsp\Gamma_0 C) ,
\]
where $C=\nu_0(k+1) {U\choose0}- ( {U\choose0}\trsp\nu
_0 ) e_{k+1}$.

We know the distribution of $\muHat$, so we could get (limiting)
unconditional results for $\gamma\trsp\wEmp$. This is not hard but a
bit tedious if we want explicit expressions, and because our focus is
mostly on first-order properties in this paper, we do not state the result.

\begin{pf*}{Proof of Proposition \ref
{prop:exactResultsWeightsGaussianCase}}
The proof follows from the representation we gave in Lemma \ref
{lemma:repLinearCombOptimalWeights}, that is,
\[
\gamma\trsp\wEmp=-\frac{1}{(\VHat_{\gamma}\trsp\SigmaHat^{-1}
\VHat_{\gamma})^{-1}(k+1,k+1)}(U\trsp0) (\VHat_{\gamma}\trsp
\SigmaHat^{-1} \VHat_{\gamma})^{-1} \myVector{0_k}{1} ,
\]
and the fact that, by the same arguments as before, conditional on
$\muHat$,
\[
(\VHat_{\gamma}\trsp\SigmaHat^{-1} \VHat_{\gamma
})^{-1} | \muHat\sim\Wishart_{k+1}\bigl((\VHat_{\gamma}\trsp
\Sigma^{-1}\VHat_{\gamma})^{-1},n-p+k\bigr)/(n-1) .
\]
We conclude that
\[
\gamma\trsp\wEmp| \muHat\equalInLaw-\frac{(U\trsp
0) W_{\gamma}{0_k\choose 1}}{W_{\gamma}(k+1,k+1)}=-\frac{\sum
_{i=1}^k u_i W_{\gamma}(i,k+1)}{W_{\gamma}(k+1,k+1)} .
\]
This shows the fist part of the proposition.

The second part follows from the following observation. Suppose the
matrix $P$ is $\Wishart_p(\id_p,K)$. If $\alpha$ and $\beta$ are
$n$-dimensional, orthogonal vectors, let us consider
\[
\frac{\alpha\trsp P \beta}{\beta\trsp P \beta} .
\]
We can, of course, write $P=\sum_{i=1}^K Y_iY_i\trsp$, where $Y_i$
are i.i.d. ${\cal N}(0,\id_p)$. In other respects, $Y_i\trsp\alpha$
and $Y_i\trsp\beta$ are clearly independent normal random variables,
since their covariance is $\alpha\trsp\beta=0$, and they are normal.
So
\[
\mathbf{E}\biggl(\frac{\alpha\trsp P \beta}{\beta\trsp P \beta
}\Big|\{
Y_i\trsp\beta\}_{i=1}^K\biggr)=0
\]
because the quantity whose expectation we are taking is a linear
combination of mean 0 independent normal random variables. Hence, also,
\[
\mathbf{E}\biggl(\frac{\alpha\trsp P \beta}{\beta\trsp P \beta}\biggr)=0 .
\]
Now, when $\alpha$ is not orthogonal to $\beta$, we write $\alpha
=\beta(\alpha\trsp\beta)/\| \beta\|_2^2+\delta$, where $\delta
$ is orthogonal to $\beta$. We immediately deduce that in general,
\[
\mathbf{E}\biggl(\frac{\alpha\trsp P \beta}{\beta\trsp P \beta}\biggr)=\frac
{\alpha\trsp\beta}{\| \beta\|_2^2}+\mathbf{E}\biggl(\frac{\delta\trsp
P \beta}{\beta\trsp P \beta}\biggr)=\frac{\alpha\trsp\beta}{\| \beta\|
_2^2} .
\]
Furthermore, when $P$ is $\Wishart_p(\Sigma,K)$, because we can write
$P=\Sigma^{1/2}P_0\Sigma^{1/2}$, where $P_0\sim\Wishart_p(\id
_p,K)$, we finally have
\[
\mathbf{E}\biggl(\frac{\alpha\trsp P \beta}{\beta\trsp P \beta}\biggr)=\frac
{\alpha\trsp\Sigma\beta}{\beta\trsp\Sigma\beta} .
\]

In the case of interest to us, we have $\alpha= {U\choose0}$,
$\beta=e_{k+1}$ and $\Sigma=\NHat_{\gamma}$. Applying the previous
formula gives us the second part of the proposition.
\end{pf*}

We now turn to the question of understanding the robustness properties
of the Gaussian results we just obtained. We will do so by studying the
same problems under more general distributional assumptions, and
specifically we will now assume that the observations are elliptically
distributed.

\section{Solutions of quadratic programs when the data is elliptically
distributed}\label{sec:EllipticalCase}
In Section~\ref{sec:GaussianCase}, we studied the properties of the
``plug-in'' solution of problem \eqref{eq:QP-eqc-Pop} under the
assumption that the data was normally distributed. While this allowed
us to shed light on the statistical properties of the solution of
problem \eqref{eq:QP-eqc-Emp}, it is naturally extremely important to
understand how robust the results are to our normality assumptions.

In this section, we will consider elliptical models, that is, models
such that the data can be expressed as
\[
X_i=\mu+\lambda_i \Sigma^{1/2} Y_i ,
\]
where $\lambda_i$ is a random variable and $Y_i$ are i.i.d. ${\cal
N}(0,\id_p)$ entries. $\lambda_i$ and $Y_i$ are assumed to be
independent, and to lift the indeterminacy between $\Sigma$ and
$\lambda$, we assume that $\mathbf{E}(\lambda_i^2)=1$. Under this
assumption, we clearly have $\operatorname{cov}(X_i)=\Sigma$. We note
that this is
not the standard definition of elliptical models, which generally
replaces $Y_i$ with a vector uniformly distributed on the sphere in
$\mathbb{R}^p$, but it captures the essence of the problem. We refer
the interested reader to \cite{anderson03} and \cite{FangKotzNg90}
for extensive discussions of elliptical distributions.

Our motivation for undertaking this study comes also from the fact that
for certain types of data, such as financial data, it is sometimes
argued that elliptical models are more reasonable than Gaussian ones,
for instance, because they can capture nontrivial tail dependence [see
\cite{FrahmJaekel05} where such models are advocated for
high-dimensional modelization of financial returns, \cite{MeucciBook}
for a discussion of their relevance for certain financial markets,
\cite{BiroliBouchaudPottersStudentEnsemble07} for modelization
considerations quite similar to \cite{FrahmJaekel05} and \cite
{EmbrechtsEtAlQRM} for a thorough discussion of tail dependence]. From
a theoretical standpoint, considering elliptical models will also help
in several other ways: the results will yield alternative proofs to
some of the results we obtained in the Gaussian case, they will allow
us to deal with some situations where the data $X_i$ are not
independent and they will also allow us to understand the properties of
the bootstrap.

We also want to point out that elliptical distributions allow us to not
fall into the geometric ``trap'' of standard random matrix models
highlighted in \cite{nekCorrEllipD}: the fact that data vectors drawn
from standard random matrix models are essentially assumed to be almost
orthogonal to one another and that their norm (after renormalization by
$1/\sqrt{p}$) is almost constant. In a sense, studying elliptical
models will allow us to understand what is the impact of the implicit
geometric assumptions made about the data when assuming normality. (We
purposely do so not under minimal assumptions but under assumptions
that capture the essence of the problem while allowing us to show in
the proofs the key stochastic phenomena at play.) This part of the
article can therefore be viewed as a continuation of the investigation
we started in \cite{nekCorrEllipD} where we showed a lack of
robustness of random matrix models (contradicting claims of
``universality'') by thoroughly investigating limiting spectral
distribution properties of high-dimensional covariance matrices when
the data is drawn according to elliptical models and generalizations.
We show here that the theoretical problems we highlighted in \cite
{nekCorrEllipD} have important practical consequences. [For more
references on elliptical models in a random matrix context, we refer
the reader to \cite{nekCorrEllipD} where an extended bibliography can
be found.]

We now turn to the problem of understanding the solution of problem
\eqref{eq:QP-eqc-Emp} in the setting where the data is elliptically
distributed. We will limit ourselves to the case where the matrix
$\VHat$ is full of known and deterministic vectors, except possibly
for the sample mean. In this section we restrict ourselves to
convergence in probability results. It is clear from Section \ref
{Sec:GeneralQP} that to tackle the problems we are considering we need
to understand at least three types of quantities: $v\trsp\SigmaHat
^{-1} v$ for a deterministic $v$ with unit norm, $\muHat\trsp
\SigmaHat^{-1} v$ and $\muHat\trsp\SigmaHat^{-1} \muHat$.

Here is a brief overview of our findings. When we consider elliptical
models, our results say that roughly speaking, under certain
assumptions given precisely later:
\begin{enumerate}[3.]
\item$\frac{v\trsp\SigmaHat^{-1} v}{v\trsp\Sigma^{-1} v} \tendsto
\limitScaling$, where $\limitScaling$ satisfies, if $G$ is the limit
law of the empirical distribution of the $\lambda_i^2$ and
$p/n\tendsto\rho\in(0,1)$, $\int\frac{dG(\tau)}{1+\tau\rho
\limitScaling}=1-\rho$.
\item If $\mu=0$, $\muHat\trsp\SigmaHat^{-1} \muHat\tendsto\rho
/(1-\rho)$.
\item If $\mu=0$, $\muHat\trsp\SigmaHat^{-1} v\tendsto0$.
\end{enumerate}
All these convergence results are to be understood in probability. They
naturally allow us---under certain conditions on the population
parameters---to conclude about the convergence in probability of the
matrix $\VHat\trsp\SigmaHat^{-1} \VHat$. The results mentioned
above are stated in all details in Theorems \ref
{thm:QuadFormsInverseEllipticalCase} and \ref
{thm:SummaryQuadFormsMuHatAndSigmaHat}.

In the situation where $\lambda_i$ are i.i.d., the results above hold
when $\lambda_i$ have a second moment and they do not put too much
mass near 0. This is interesting in practice because
it tells us that our results hold for heavy-tailed data, which are of
particular interest in some financial applications.

The bootstrap situation corresponds basically to $G$ being Poisson(1),
which we denote by $\Poisson(1)$. Also in the statement above for
$\muHat\trsp\SigmaHat^{-1}\muHat$, one should replace $\rho
/(1-\rho)$ by $\limitScaling-1$ in the bootstrap case. This is
explained in Theorem \ref
{thm:bootQuadFormsInverseEllipticalCaseNormalCase} and  Section
\ref{subsubsec:BootResQuadFormsMuHatSigmaHat}. Finally, in the case of
Gaussian data with ``temporal'' correlation, that is, when the data can
be written in matrix form $X=\ebold_n\mu\trsp+\Lambda Y \Sigma
^{1/2}$, where $\Lambda$ is not diagonal (and $\ebold_n$ is an
$n$-dimensional vector with only 1's in its entries), one should replace
$G$ by the limiting spectral distribution of $\Lambda\trsp\Lambda$.
The question of convergence of $\muHat\trsp\SigmaHat^{-1} \muHat$
is then more involved. We refer to Proposition \ref
{proposition:QuadFormsMuHatSigmaInverseCorrCase} for details about this
situation.

Though we are taking a fundamentally random matrix theoretic approach,
our presentation purposely avoids borrowing too many techniques from
random matrix theory in the hope of making clear(er) the phenomena that
yield the results we will obtain. A more general but considerably more
technically complicated (for non-specialists of random matrix theory)
approach is being developed in our study of a connected problem and
will appear in another paper.

This section is divided into four subsections. The first two are
devoted to the main technical issues arising in the study of the
problem when the data is elliptically distributed. The third discusses
the impact of correlation between observations when the data is
Gaussian, as it can be recast as a variant of elliptical problems. The
last subsection discusses questions related to the (nonparametric) bootstrap.

\subsection{\texorpdfstring{On quadratic forms of the type $v\trsp\SigmaHat^{-1} v$}
{On quadratic forms of the type v' Sigma -1 v}}

The focus of this subsection is on understanding statistics of the type
$v\trsp\SigmaHat^{-1} v$, where $v$ is a deterministic vector. We
will prove the following important theorem.
\begin{theorem}\label{thm:QuadFormsInverseEllipticalCase}
Suppose we observe $n$ observations $X_i$, where $X_i$ has the form
$X_i=\mu+\lambda_i \Sigma^{1/2} Y_i$, with $Y_i\iid{\cal N}(0,\id
_p)$ and $\{\lambda_i\}_{i=1}^n$ is independent of $\{Y_i\}_{i=1}^n$.
$\Sigma^{1/2}$ is deterministic and $\mathbf{E}(\lambda_i^2)=1$.

We call $\rho_n=p/n$ and assume that $\rho_n\tendsto\rho\in(0,1)$.

We use the notation $\tau_i=\lambda_i^2$ and assume that the
empirical distribution, $G_n$, of $\tau_i$ converges weakly in
probability to a deterministic limit $G$. We also assume that $\tau
_i\neq0$ for all $i$.

If $\tau_{(i)}$ is the $i$th largest $\tau_k$, we assume that we can
find a random variable $N \in\mathbb{N}$ and positive real numbers
$\eps_0$ and $C_0$ such that
\renewcommand{\theequation}{Assumption-BB}
\begin{equation}\label{eq:AssumptionSmallestSingularValueBoundedBelow}
\hspace*{80pt}\cases{\displaystyle
P(p/N<1-\eps_0)\tendsto1,  \qquad\mbox{as }n\tendsto\infty,\cr\displaystyle
P\bigl(\tau_{(N)}>C_0\bigr)\tendsto1 ,\cr\displaystyle
\exists\eta_0>0   \quad \mbox{such that}\quad   P(N/n>\eta_0)\tendsto1, \qquad\mbox{as } n\tendsto\infty .
}
\end{equation}

Under these assumptions, if $v$ is a (sequence of) deterministic vector,
\[
\frac{v\trsp\SigmaHat^{-1} v}{v\trsp\Sigma^{-1} v}\tendsto
\limitScaling\qquad\mbox{in probability} ,
\]
where $\limitScaling$ satisfies
\renewcommand{\theequation}{\arabic{equation}}\setcounter{equation}{3}
\begin{equation}\label{eq:defLimitScaling}
\int\frac{dG(\tau)}{1+\rho\tau\limitScaling}=1-\rho.
\end{equation}
\end{theorem}

A few comments are in order before we turn to the proof. First, the
assumption that $\lambda_i\neq0$ for all $i$ could be dispensed of,
as long as all assumptions stated above hold when $n$ is understood to
denote the number of nonzero $\lambda_i$'s. Second, \eqref
{eq:AssumptionSmallestSingularValueBoundedBelow} concerning $N$ and $C$
will generally hold as soon as $G$ does not put too much mass at 0, the
only problem-specific question remaining being how much mass is put at
0 by $G$ compared to $\rho$, the limit of $p/n$.

In particular, in the case where the $\tau_i$'s are i.i.d., if there
exists $C_0>0$ and $x_0>0$ such that $P_G(X>C_0)=x_0>0$, and if $G_n$
is the empirical distribution of the $\tau_i$'s, if $G_n\WeakCv G$, we
see, using, for example, Lemma 2.2 in \cite{vandervaart}, that
\[
\liminf_{n\tendsto\infty} P_{G_n}(X>C_0)=\frac{\operatorname
{Card}\{\tau
_i>C_0\}}{n}\geq P_G(X>C_0)=x_0 .
\]
So picking $N=(1-\delta)x_0 n$ will guarantee that we have, if
$G_n\WeakCv G$ in probability, $P(\tau_{(N)}>C_0)\tendsto1$ and, of
course, $P(N/n>\eta)\tendsto1$. Hence, in checking whether the
theorem applies, we just need to see whether $p/N$ stays bounded away
from 1.

In the simpler case when all the $|\lambda_i|$ are bounded away from
0, the conditions on $N$ and $C$ apply directly by taking $N=n$.
Finally, let us say that \eqref
{eq:AssumptionSmallestSingularValueBoundedBelow} is needed in the proof
to guarantee that the smallest eigenvalues of $\SigmaHat$ stay bounded
away from 0 with high-probability.\vadjust{\goodbreak}

We now briefly compare the Gaussian and elliptical cases. A simple
convexity argument [relying on the fact that $1/(1+x)$ is a convex
function of $x$ for $x\geq0$ and Jensen's inequality] shows that, if
$\mu_G$ is the mean of $G$,
\[
\limitScaling\geq\frac{1}{1-\rho}\frac{1}{\mu_G} .
\]
In the case of Gaussian data, $G=\delta_{1}$, that is, it is a point
mass at 1 and we have $\limitScaling=1/(1-\rho)$. In other respects,
for $X_i$ to have covariance $\Sigma$, we need $\mathbf{E}(\lambda
_i^2)=1$. When the $\lambda_i$'s are i.i.d., with $\lambda_i^2$ having
distribution $G$, $\mu_G=\mathbf{E}(\lambda_i^2)=1$, and we know that
$G_n\WeakCv G$ in probability. Therefore, in the class of elliptical
distributions considered here, risk underestimation, which is
essentially measured by $1/\limitScaling$ (see Theorem \ref
{thm:SolnGeneralQP} and Section \ref
{sec:ComparisonGaussianElliptical}) will be least severe in the
Gaussian case. In other words, the Gaussian results lead to
over-optimistic conclusions (in terms of proximity between sample and
population solutions of the quadratic programs we are considering)
within the class of elliptical distributions.

We go back to these questions in more detail in Section \ref
{sec:ComparisonGaussianElliptical} and now turn to the proof of Theorem
\ref{thm:QuadFormsInverseEllipticalCase}. The proof could be carried
out in at least two ways. We take one that is not standard but we feel
best explains the phenomenon that is occurring.\vspace*{-3pt}

\begin{pf*}{Proof of Theorem \ref
{thm:QuadFormsInverseEllipticalCase}}
The proof is easier to carry out when we write the problem in matrix
form. Because we focus on $\SigmaHat$, we can assume without loss of
generality (wlog) that $\mu=0$. Let us consider the $n\times p$ data
matrix $X$ whose $i$th row is $X_i$. Similarly, we denote by $Y$ the
$n\times p$ data matrix whose $i$th row is $Y_i$. Let us call $\Lambda
$ the diagonal matrix with $i$th diagonal entry $\lambda_i$ and $H=\id
_n-\ebold\ebold\trsp/n$, where $\ebold$ is an $n$-dimensional
vector whose entries are all equal to 1. Note that $H\trsp H=H$. With
these notations, we have, since we assume that $\mu=0$,
\[
X=\Lambda Y \Sigma^{1/2} .
\]
Therefore, $X-\bar{X}=HX$, and
\[
\SigmaHat=\frac{1}{n-1} (X-\bar{X})\trsp(X-\bar{X})= \frac
{1}{n-1} \Sigma^{1/2} Y\trsp\Lambda H\Lambda Y \Sigma^{1/2}.
\]
Let us call $L$ the matrix $L=\Lambda H \Lambda$. Note that $Y\trsp L
Y$ is a rank $p$ matrix with probability 1, if we assume that $p\leq
n-1$ (recall that all the entries of $\Lambda$ are nonzero). Hence,
$Y\trsp L Y$ is invertible with probability 1. Therefore,
\[
\SigmaHat^{-1}=\Sigma^{-1/2} \biggl(\frac{1}{n-1}Y\trsp L Y
\biggr)^{-1} \Sigma^{-1/2} .
\]
Finally, we have
\[
\frac{v\trsp\SigmaHat^{-1} v}{v\trsp\Sigma^{-1} v}=\nu\trsp
\biggl(\frac{1}{n-1}Y\trsp L Y \biggr)^{-1} \nu,
\]
where $\nu=\Sigma^{-1/2}v/\| \Sigma^{-1/2}v\|_2$ is a vector of
$\ell_2$ norm 1.\vadjust{\goodbreak}

We now make all of our statements conditional on $\Lambda$. Because of
the independence of $Y$ and $\Lambda$, we can therefore treat the
$\lambda_i$'s as if they were constant and the $Y_{i,j}$'s as i.i.d.
${\cal N}(0,1)$ random variables. $\Lambda$ is now assumed to be in
the set of matrices ${\cal L}_{\eps,\delta}$, defined just below, for
which we have control of the smallest eigenvalue of ${\cal S}=Y\trsp L
Y/(n-1)$. In the steps that follow that are conditional on $\Lambda$,
we therefore consider that we control the smallest eigenvalue of ${\cal
S}$. We note that if $\Lambda$ is in ${\cal L}_{\eps,\delta}$, $N$
is lower bounded. Because $N$ is a function of the $\lambda_i$'s and
hence of $\Lambda$, we write all the results conditionally on $\Lambda
$, but the reader should keep in mind that this conditioning constrains
also the possible values of $N$.

$\bullet$ \textit{The set ${\cal L}_{\eps,\delta}$}.

In Lemma~\ref{lemma:LowerBoundOnSmallestEig} in the \hyperref[appm]{Appendix}, we prove
the following result: when $\Lambda$ is such that $p/N<1-\eps$, if
$\mathfrak{C}_n=C_0 \frac{N-1}{n-1}$ [see \eqref
{eq:AssumptionSmallestSingularValueBoundedBelow} and Lemma~\ref
{lemma:LowerBoundOnSmallestEig} for definitions] and $\gamma_p$ is the
smallest eigenvalue of $Y\trsp L Y/n-1$, we have, if $P_{\Lambda}$
denotes probability conditional on $\Lambda$,
\[
P_{\Lambda} \bigl(\sqrt{\gamma_p}\leq\sqrt{\mathfrak{C}_n}
\bigl[\bigl(1-\sqrt{1-\eps}\bigr)-t \bigr] \bigr)\leq\exp\bigl(-(N-1)t^2
\bigr) .
\]

Let us call ${\cal L}_{\eps,\delta}$ the set of matrices $\Lambda$
such that $p/N<1-\eps$ and $C_0 (N-1)/(n-1)>\delta$. Under \eqref
{eq:AssumptionSmallestSingularValueBoundedBelow}, for a $\delta$
bounded away from 0 (e.g., $\delta=C_0 \eta_0/2$, since we need a
bound on $\liminf C_0 N/n$ that holds with probability going to 1),
$P(\Lambda\in{\cal L}_{\eps,\delta})\tendsto1$. In other respects,
if $\Lambda\in{\cal L}_{\eps,\delta}$,
\[
P_{\Lambda} \bigl(\sqrt{\gamma_p}\leq\sqrt{\delta} \bigl[\bigl(1-\sqrt
{1-\eps}\bigr)-t \bigr] \bigr)\leq\exp\bigl(-(n-1)\delta t^2/C_0
\bigr) .
\]

$\bullet$ \textit{Getting results conditionally on ${\Lambda}$}.

If $O$ is an orthogonal matrix, $O\trsp Y\trsp L Y O \equalInLaw Y\trsp
L Y$, because $Y$ is full of i.i.d. ${\cal N}(0,1)$ random variables
and is therefore invariant (in law) by left and right rotation.
Therefore the eigenvalues and eigenvectors of $Y\trsp L Y$ are
independent and its matrix of eigenvectors is uniformly (i.e., Haar)
distributed on the orthogonal group [see also \cite{ChikuseBook03},
page 40, equation (2.4.4)]. Let us write a spectral decomposition of
$Y\trsp L Y$
\[
{\cal S}=\frac{1}{n-1}Y\trsp L Y=\sum_{i=1}^p \gamma_i \upsilon_i
\upsilon_i\trsp.
\]
We know that a.s. $\gamma_i\neq0$ for all $i$, so
\[
\nu\trsp{\cal S}^{-1} \nu=\sum_{i=1}^p \frac{1}{\gamma_i} (\nu
\trsp\upsilon_i)^2 .
\]
We claim that
\[
\nu\trsp{\cal S}^{-1} \nu- \frac{1}{p}\sum_{i=1}^p \frac
{1}{\gamma_i} \Big| (\{\gamma_i\}_{i=1}^p, \Lambda)
\tendsto0.
\]
To see this, note that $\mathbf{E}((\nu\trsp\upsilon_i)^2)=\| \nu\|
_2^2/p=1/p$ because $\upsilon_i$ is uniformly distributed on the unit
sphere when $\Upsilon$ (the matrix containing the $\upsilon_i$)\vadjust{\goodbreak} is
Haar distributed on the orthogonal group. Hence, given the independence
between $\gamma_i$ and $\upsilon_i$,
\[
\mathbf{E}(\nu\trsp{\cal S}^{-1} \nu|\{\gamma_i \}
_{i=1}^n,\Lambda)=\frac{1}{p}\sum_{i=1}^p \frac{1}{\gamma_i} .
\]
Now let us call $w$ the vector with $w_i=(\nu\trsp\upsilon_i)^2$,
and $g$ the vector with $i$th entry $g_i=1/\gamma_i$. Clearly, since
$\nu\trsp{\cal S}^{-1} \nu=g \trsp w$, $\operatorname{var}(\nu
\trsp{\cal S}^{-1} \nu|\{\gamma_i\},\Lambda)=g \trsp
\operatorname{cov}(w)
g$. By symmetry it is clear that $\operatorname
{cov}(w)(i,i)=\operatorname{cov}(w)(1,1)$ and
$\operatorname{cov}(w)(i,j)=\operatorname{cov}(w)(1,2)$ if $i\neq j$.
Further, since the matrix
$\Upsilon$ containing the vectors $\upsilon_i$ is Haar distributed on
the orthogonal group, we can assume without loss of generality that
$\nu=e_1$ for all the computations at stake. As a matter of fact, if
$O_1$ is an orthogonal matrix such that $O_1\nu=e_1$, then $\nu\trsp
\upsilon_i=e_1\trsp O_1\upsilon_i=e_1\trsp\tilde{\upsilon}_i$
where the matrix $\widetilde{\Upsilon}=O_1 \Upsilon$ is again Haar
distributed on the orthogonal group.

So from now on, we assume (without loss of generality) that $\nu=e_1$,
and we therefore simply need to understand the correlation between
$(\upsilon_1(1))^2$ and $(\upsilon_2(1))^2$. Now, the first row of an
orthogonal matrix uniformly distributed on the orthogonal group is a
unit vector uniformly distributed on the unit sphere, because if $O$ is
Haar distributed, so is $O\trsp$. We now recall the fact that a vector
uniformly distributed on the unit sphere, $\upsilon$ can be generated
by drawing at random a ${\cal N}(0,\id_p)$ random vector and
normalizing it. In other words, if $Z\sim{\cal N}(0,\id_p)$,
$\upsilon=Z/\| Z\|_2$.

So our task has now been considerably simplified, and it consists in
understanding the covariance between 2 random variables, $r_1$ and
$r_2$ such that, if $Z_i$ are i.i.d. ${\cal N}(0,1)$,
\[
r_i=\frac{Z_i^2}{\sum_{i=1}^p Z_i^2} .
\]
Now, by symmetry, $\mathbf{E}(r_1 r_2)=\mathbf{E}(r_i r_j)$ for all
$i\neq j$ and
$p(p-1)\mathbf{E}(r_1 r_2)=\sum_{i\neq j} \mathbf{E}(r_i r_j)$. In
other words,
\[
p(p-1)\mathbf{E}(r_1 r_2)=\mathbf{E}\biggl(\sum_{i\neq j}\frac{Z_i^2
Z_j^2}{(\sum _{i=1}^p Z_i^2)^2}\biggr)=\mathbf{E}\biggl(\frac{\sum_{i,j}Z_i^2
Z_j^2}{(\sum_{i=1}^p Z_i^2)^2}-\sum_{i=1}^p \frac{Z_i^4}{(\sum
_{i=1}^p Z_i^2)^2}\biggr) .
\]
We can therefore conclude that
\[
p(p-1)\mathbf{E}(r_1 r_2)=1-p \mathbf{E}\biggl(\frac{Z_1^4}{(\sum_{i=1}^p
Z_i^2)^2}\biggr) .
\]
Hence, $\mathbf{E}(r_1 r_2)\leq1/(p(p-1))$. On the other hand,
\[
\mathbf{E}\biggl(\frac{Z_1^4}{(\sum_{i=1}^p Z_i^2)^2}\biggr)\leq\mathbf
{E}\biggl(\frac {Z_1^4}{(\sum_{i=2}^p Z_i^2)^2}\biggr)=\frac{3}{(p-3)(p-5)} ,
\]
since $\sum_{i=2}^p Z_i^2\sim\chi^2_{p-1}$, and $\mathbf{E}((\chi
^2_{p-1})^r)=2^r \Gamma((p-1)/2+r)/\Gamma((p-1)/2)$, for $r>-(p-1)/2$
[see, e.g., \cite{mardiakentbibby}, page 487].\vadjust{\goodbreak} Applying these results with
$r=-2$ yields the above result as soon as $p>5$, by using the fact that
$\Gamma(x+1)=x\Gamma(x)$. We therefore have
\[
1-\frac{3p}{(p-3)(p-5)}\leq p(p-1)\mathbf{E}(r_1 r_2) \leq1 .
\]
Since, for instance by symmetry, $\mathbf{E}(r_1)=1/p$, and
$1/(p(p-1))-1/p^2=(p^2(p-1))^{-1}$, we conclude that
\[
\frac{1}{p^2(p-1)}-\frac{3p}{p(p-1)(p-3)(p-5)}\leq\operatorname
{cov}(r_1,r_2)\leq\frac{1}{p^2(p-1)} .
\]
We have therefore established the fact that
\[
|\operatorname{cov}(r_1,r_2) |=\gO(p^{-3}) .
\]

On the other hand, since $\mathbf{E}(r_1^2)=\mathbf{E}(Z_1^4(\sum
_{i=1}^p Z_i^2)^{-2})$, we have
\[
0\leq\operatorname{var}(r_i)\leq\frac{3}{(p-3)(p-5)}-\frac
{1}{p^2} .
\]
Now using the (standard) fact that, for symmetric matrices $M$, if
$\sigma_1(M)$ is the largest singular value of $M$,
\[
\sigma_1(M)\leq\max_{i} \sum_{j} |m_{i,j}| ,
\]
[it can easily be proved using, for instance, Theorems 5.6.6 and 5.6.9
in \cite{hornjohnson94}, or Ger\v{s}gorin's theorem (Theorem 6.1.1 in
the same reference)] we have
\[
\sigma_1(\operatorname{cov}(r))\leq\biggl(\frac
{3}{(p-3)(p-5)}-\frac
{1}{p^2} \biggr)+\gO(p^{-2})=\gO(p^{-2}) .
\]
The first term in the previous bound comes from the contribution of the
diagonal and the second term is the sum over the $p-1$ off-diagonal
elements on a given row of the upper-bound we had on each such element,
that is, $Cp^{-3}$ for some $C$.

Let us now return to our initial question which was to show that the
conditional variance of interest to us was going to zero. Recall that
$g$ is a vector whose $i$th entry is $1/\gamma_i$. Since
\[
\operatorname{var}(\nu\trsp{\cal S}^{-1} \nu|\{\gamma_i\}
,\Lambda)=g \trsp\operatorname{cov}(w)g ,
\]
and $\operatorname{cov}(w)=\operatorname{cov}(r)$, we have, for $C$ a
constant, and if $\opnorm
{A}$ denotes the operator norm (or largest singular value) of the
matrix $A$,
\[
\operatorname{var}(\nu\trsp{\cal S}^{-1} \nu|\{\gamma_i\}
,\Lambda)\leq\opnorm{\operatorname{cov}(r)} \| g\|_2^2\leq
C \frac{\| g\|_2^2}{p^2}=C\frac{1}{p^2}\sum_{i=1}^p \frac{1}{\gamma
_i^2} .
\]
Now given the assumptions we made on $\Lambda$, according to the
arguments given at the beginning of this proof and Lemma \ref
{lemma:LowerBoundOnSmallestEig} in the \hyperref[appm]{Appendix},\vadjust{\goodbreak} $\gamma_i^2 \geq
\mathfrak{C}_n(1-\sqrt{p/(N-1)})^2/2$, where $\mathfrak{C}_n=C_0
(N-1)/(n-1)$, with high ($\{Y_i\}_{i=1}^n$)-probability. So we conclude
that all the $\gamma_i$'s are bounded away [uniformly for $\Lambda$ in
${\cal L}_{\eps,\delta}$ and with high ($\{Y_i\}
_{i=1}^n$)-probability] from 0, and when this is the case,
\[
\operatorname{var}(\nu\trsp{\cal S}^{-1} \nu|\{\gamma_i\}
,\Lambda)\tendsto0 .
\]
Therefore,
\[
\nu\trsp{\cal S}^{-1} \nu-\frac{1}{p}\sum_{i=1}^p \frac
{1}{\gamma_i}\Big |\{\gamma_i\}_{i=1}^p,\Lambda\tendsto0
\qquad\mbox{in probability}.
\]
Let us now show that this implies convergence in probability to 0
(conditional on $\Lambda$ only) of $Q_n=\nu\trsp{\cal S}^{-1} \nu-
\frac{1}{p}\sum_{i=1}^p \frac{1}{\gamma_i}$. Let us call $h_n=C
\frac{\| g\|_2^2}{p^2}=\frac{C}{p^2}\sum_{i=1}^p \frac{1}{\gamma
_i^2}$. For $\zeta_n$ to be determined later, we have
\[
P(|Q_n|>\eps|\Lambda)\leq P(|Q_n|>\eps\ \&\ h_n \leq\zeta
_n|\Lambda)+P(h_n>\zeta_n|\Lambda) .
\]
On the other hand,
\[
P(|Q_n|>\eps\ \&\ h_n \leq\zeta_n|\Lambda)=\mathbf{E}\bigl(\mathbf
{E}\bigl(1_{|Q_n|>\eps} 1_{h_n \leq\zeta_n}|\{g_i\},\Lambda\bigr)|\Lambda\bigr) .
\]
Because $h_n$ is a function of the $g_i$'s and $\operatorname
{var}(Q_n|\{\gamma_i\} _{i=1}^p,\Lambda)\leq h_n$,
\[
\mathbf{E}\bigl(1_{|Q_n|>\eps} 1_{h_n \leq\zeta_n}|\{g_i\},\Lambda\bigr)=1_{h_n
\leq\zeta_n}\mathbf{E}\bigl(1_{|Q_n|>\eps}|\{g_i\},\Lambda\bigr)\leq1_{h_n
\leq
\zeta_n} \frac{h_n}{\eps^2}\leq\frac{\zeta_n}{\eps^2} .
\]
But when $\Lambda\in{\cal L}_{\eps,\delta}$, under our assumptions
and their consequences on the $\gamma_i^2$'s mentioned above [i.e.,
$\gamma_i^2 \geq\mathfrak{C}_n(1-\sqrt{p/(N-1)})^2/2$ with high $\{
Y_i\}_{i=1}^n$ probability], we have $h_n|\Lambda=\gO_P(1/p)$, so
taking $\zeta_n=n^{-1/2}$, we have $P(h_n>\zeta_n|\Lambda)\tendsto
0$ and of course, $\zeta_n/\eps^2\tendsto0$. Hence, for any $\eps>0$,
\[
P(|Q_n|>\eps|\Lambda)\tendsto0 .
\]

Let us now turn to the question of identifying the limit.

$\bullet$ \textit{About $\frac{1}{p}\sum_{i=1}^p \frac
{1}{\gamma_i}$}.

The Stieltjes transform of the spectral distribution of $Y\trsp L
Y/(n-1)$ is
\[
s_p(z)=\frac{1}{p}\sum_{i=1}^p \frac{1}{\gamma_i -z} .
\]
The quantity $\frac{1}{p}\sum_{i=1}^p \frac{1}{\gamma_i}$ is
therefore $s_p(0)$ and we are interested in its limit, if it exists,
which would correspond to $\limitScaling$.

Recall the Mar\v{c}enko--Pastur equation, from \cite{mp67}, \cite
{wachter78} and
\cite{silverstein95}: if $Y$ is $n\times p$ has i.i.d. entries with
mean 0 and variance 1 and $L$ is positive semidefinite, has limiting
spectral distribution $G$ and is independent of $Y$, if $p/n\tendsto
\rho>0$, and if $m_p$ is the Stieltjes transform of the spectral
distribution of $Y\trsp L Y/p$, then $m_p(z)$ tends\vadjust{\goodbreak} (in probability) to
$m(z)$ for all $z$ in $\mathbb{C}^+$ and $m$ satisfies
\renewcommand{\theequation}{\arabic{equation}}\setcounter
{equation}{4}
\begin{equation}\label{eq:mapa}
-\frac{1}{m(z)}=z-\frac{1}{\rho}\int\frac{\tau \,dG(\tau)}{1+\tau
m(z)} .
\end{equation}
Note that, if $p/n=\rho_n$, we have
\[
\rho_n s_p(\rho_n z)=m_p(z) .
\]

Therefore, according to \cite{mp67}, \cite{wachter78} and \cite
{silverstein95}, we know that $s_p(z)$ converges for $z\in\mathbb
{C}^+$ to a nonrandom quantity $s(z)$, in probability. Note that $s$
satisfies, in light of equation \eqref{eq:mapa},
\[
-\frac{1}{s(z)}=z-\int\frac{\tau \,dG(\tau)}{1+\tau\rho s(z)} .
\]

Here, because we know using our assumptions (see the end of the proof)
that $\gamma_i$ are bounded away from 0 with probability going to 1,
we can also conclude that $s_p(0)\tendsto s(0)$ with probability going
to 1, because of the weak convergence (in probability) of spectral
distributions that pointwise convergence of Stieltjes transforms
implies (as a test function, we can use a function that coincides with
$1/x$ except in a interval near 0 where we are guaranteed that there
are no eigenvalues asymptotically). We also know that $s$ is continuous
(and actually analytic) at $0$ in this situation since the $s$ is the
Stieltjes transform of a measure who has support bounded away from 0.
So the previous equation holds for $z=0$, and we have
\[
-\frac{1}{s(0)}=-\int\frac{\tau \,dG(\tau)}{1+\tau\rho s(0)} .
\]
Multiplying both sides by $-\rho s(0)$, we get, after we recall that
$G$ is a probability measure,
\[
\rho=\int\frac{\rho\tau s(0)\,dG(\tau)}{1+\tau\rho s(0)}=\int
\biggl( 1-\frac{1}{1+\tau\rho s(0)} \biggr) \,dG(\tau)= 1- \int\frac
{1}{1+\tau\rho s(0)} \,dG(\tau) .
\]
Calling $s(0)=\limitScaling$, we have the result we announced,
conditionally on $\Lambda$. Now, here $G$ is the limiting spectral
distribution of $\Lambda H \Lambda$, but because this matrix is a rank
one perturbation of $\Lambda^2$, these two matrices have the same
limiting spectral distribution. This concludes this part of the proof.

$\bullet$ \textit{Getting results unconditionally on ${\Lambda}$}.

All the statements above were made conditional on $\Lambda$. If we can
show that our probability bounds and our characterization of the limit
hold uniformly in $\Lambda$, we will have an unconditional statement,
as we seek.

The fact that the limit does not depend on $\Lambda$ is essentially
obvious from its description: all that matters is the limiting spectral
distribution, which is the same for all $\Lambda$. Let us consider the
question of uniform probability bounds. All we need to do is show that
we control $P(h_n>\zeta_n|\Lambda)$ uniformly in $\Lambda$. At this
point, it is helpful to recall that $N$ can be viewed as a function of
$\Lambda$.\looseness=-1\vadjust{\goodbreak}

Recall also that if $\Lambda\in{\cal L}_{\eps,\delta}$,
\[
P_{\Lambda} \bigl(\sqrt{\gamma_p}\leq\sqrt{\delta} \bigl[\bigl(1-\sqrt
{1-\eps}\bigr)-t \bigr] \bigr)\leq\exp\bigl(-(n-1)\delta t^2/C_0
\bigr) .
\]
Hence, when $\Lambda\in{\cal L}_{\eps,\delta}$, if $\zeta
_n=n^{-1/2}$, $P(h_n>\zeta_n|\Lambda)\leq f_n(C_0,\eps,\delta)$,
where $f_n(C_0,\eps,\delta)$ tends to 0 as $n$ tends to infinity. In
other words, we have now established that if $\Lambda\in{\cal
L}_{\eps,\delta}$, and $Q_n=\nu\trsp{\cal S}^{-1} \nu- \frac
{1}{p}\sum_{i=1}^p \frac{1}{\gamma_i}$, for any $t>0$,
\[
P(|Q_n|>t|\Lambda)\leq\frac{\zeta_n}{t^2}+f_n(C_0,\eps,\delta) .
\]
Using the fact that $P(|Q_n|>t)\leq P (|Q_n|>t\ \&\ \{\Lambda\in
{\cal L}_{\eps,\delta}\} )+P (\Lambda\notin{\cal
L}_{\eps,\delta} )$, we conclude that $P(|Q_n|>t)\tendsto0$ as
$n$ tends to infinity for any $t>0$ and the proof is complete.
\end{pf*}

As a consequence of Theorem~\ref{thm:QuadFormsInverseEllipticalCase},
we have the following practically useful result.

\begin{lemmaNumThm}\label{lemma:DotProductsInvolvingSigmaInverse}
We assume that the assumptions of Theorem \ref
{thm:QuadFormsInverseEllipticalCase} hold and that $G$ is such that
$\limitScaling$ is not $\infty$.

Suppose that $v_1$ and $v_2$ are deterministic vectors such that
\[
\frac{v_1\trsp\Sigma^{-1} v_2}{(v_1+v_2)\trsp\Sigma^{-1}
(v_1+v_2)} \quad \mbox{and}\quad  \frac{v_1\trsp\Sigma^{-1}
v_2}{(v_1-v_2)\trsp\Sigma^{-1} (v_1-v_2)}
\]
are bounded away from 0.
Then under the assumptions of Theorem~\ref{thm:QuadFormsInverseEllipticalCase},
\[
\frac{v_1\trsp\SigmaHat^{-1} v_2}{v_1\trsp\Sigma^{-1} v_2}
\tendsto\limitScaling\qquad\mbox{in probability}.
\]

In other respects, suppose that $v_1\trsp\Sigma^{-1} v_2\tendsto0$,
while $v_1\trsp\Sigma^{-1} v_1$ and $v_2\trsp\Sigma^{-1} v_2$ stay
bounded away from $\infty$. Then, under the assumptions of Theorem
\ref{thm:QuadFormsInverseEllipticalCase},
\[
v_1\trsp\SigmaHat^{-1} v_2 \tendsto0 \qquad\mbox{in probability}.
\]
\end{lemmaNumThm}

\begin{pf}
The proof of the first part of the lemma is an immediate consequence of
Theorem~\ref{thm:QuadFormsInverseEllipticalCase}, after writing
\begin{eqnarray*}
2\frac{v_1\trsp\SigmaHat^{-1} v_2}{v_1\trsp\Sigma^{-1} v_2}&=&\frac
{(v_1+v_2)\trsp\SigmaHat^{-1}(v_1+v_2)}{(v_1+v_2)\trsp\Sigma^{-1}
(v_1+v_2)}\frac{(v_1+v_2)\trsp\Sigma^{-1} (v_1+v_2)}{v_1\trsp\Sigma
^{-1} v_2}\\
&&{}-\frac{(v_1-v_2)\trsp\SigmaHat^{-1}(v_1-v_2)}{(v_1-v_2)\trsp\Sigma
^{-1} (v_1-v_2)}\frac{(v_1-v_2)\trsp\Sigma^{-1} (v_1-v_2)}{v_1\trsp
\Sigma^{-1} v_2}.
\end{eqnarray*}
For the proof of the second part, we note that Theorem \ref
{thm:QuadFormsInverseEllipticalCase} implies that
\[
v\trsp\SigmaHat^{-1} v=\limitScaling v\trsp\Sigma^{-1} v +\lo
_P(v\trsp\Sigma^{-1} v) .
\]
Note that since for $i=1,2$, $v_i\trsp\Sigma^{-1} v_i$ is assumed to
stay bounded, the same is true of $(v_1+\eps v_2)\trsp\Sigma^{-1}
(v_1+\eps v_2)$, where $\eps=\pm1$. Now we write
\[
2 v_1\trsp\SigmaHat^{-1} v_2=(v_1+v_2)\trsp\SigmaHat^{-1}
(v_1+v_2)-(v_1-v_2)\trsp\SigmaHat^{-1} (v_1-v_2) .\vadjust{\goodbreak}
\]
Our previous remark and the assumption of boundedness of $v_i\trsp
\Sigma^{-1} v_i$ implies that, when $v_1\trsp\Sigma^{-1} v_2\tendsto0$,
\begin{eqnarray*}
2 v_1\trsp\SigmaHat^{-1} v_2&=&\limitScaling\bigl((v_1+v_2)\trsp\Sigma
^{-1} (v_1+v_2)-(v_1-v_2)\trsp\Sigma^{-1} (v_1-v_2)\bigr)+\lo_P(1)  \\
&=&\limitScaling2 v_1\trsp\Sigma^{-1} v_2+ \lo_P(1)=\lo_P(1) .
\end{eqnarray*}
\upqed
\end{pf}

\subsection{\texorpdfstring{On quadratic forms involving $\muHat$ and ${\SigmaHat^{-1}}$}
{On quadratic forms involving mu and Sigma -1}}

As is clear from the solutions of problems \eqref{eq:GeneralQP} and
\eqref{eq:QP-eqc-Emp}, when $\muHat$ appears in the matrix $\VHat$,
its influence on the solution of our quadratic program\vspace*{2pt} will manifest
itself in the form of quantities of the type $\muHat\trsp\SigmaHat
^{-1}\muHat$ and $v_i\trsp\SigmaHat^{-1}\muHat$. It is therefore
important that we get a good understanding of those quantities.

Compared to the Gaussian case, in the elliptical case, $\muHat$ is not
independent of $\SigmaHat$ anymore, which generates some
complications. They are fully addressed in Theorem \ref
{thm:SummaryQuadFormsMuHatAndSigmaHat}, but as a stepping stone to that
result (the main of this subsection), we need the following theorem,
which essentially takes care of the problem of understanding $\muHat
\trsp\SigmaHat^{-1}\muHat$ for the class of elliptical distributions
we consider when the population mean is 0.

\begin{theorem}\label{thm:simplequadformmuhatandsigmahat}
Suppose $Y$ is an $n\times p$ matrix whose rows are the vectors $Y_i$,
which are i.i.d. ${\cal N}(0,\id_p)$.

Suppose $\Lambda$ is a diagonal matrix whose $i$th entry is $\lambda
_i$, which is possibly random and is independent of $Y$. Call $\tau
_i=\lambda_i^2$. We assume that $\tau_i\neq0$ for all $i$ and
\renewcommand{\theequation}{Assumption-BLa}
\begin{equation}\label{eq:AssumptionControlFourthEmpiricalMomentLambda}
\frac{1}{n^2}\sum_{i=1}^n \lambda_i^4=\frac{1}{n^2} \sum_{i=1}^n
\tau_i^2 \tendsto0\qquad \mbox{in probability} .\hspace*{-50pt}
\end{equation}
If $\tau_{(i)}$ is the $i$th largest $\tau_k$, we assume that we can
find a random variable $N \in\mathbb{N}$ and positive real numbers
$\eps_0$ and $C_0$ such that
\renewcommand{\theequation}{Assumption-BB}
\begin{equation}
\cases{\displaystyle
P(p/N<1-\eps_0)\tendsto1,  \qquad\mbox{as } n\tendsto\infty,\cr\displaystyle
P\bigl(\tau_{(N)}>C_0\bigr)\tendsto1,\cr\displaystyle
\exists\eta_0>0  \quad  \mbox{such that}\quad  P(N/n>\eta_0)\tendsto1,\qquad \mbox{as } n\tendsto\infty.
}\hspace*{-100pt}
\end{equation}

Let us call $\rho_n=p/n$ and $\rho=\lim_{n\tendsto\infty} \rho
_n$. We assume that $\rho\in(0,1)$. We call
\[
Z_{n,p}=\frac{1}{n^2}\ebold\trsp\Lambda Y (Y\trsp\Lambda^2
Y/n)^{-1} Y\trsp\Lambda\ebold.
\]

Then we have
\[
Z_{n,p}\tendsto\rho\qquad \mbox{in probability. }
\]
If the $n\times p$ data matrix $\XTilde$ is written $\XTilde=\Lambda
Y \Sigma^{1/2}$, and if $\hat{m}=\Sigma^{1/2}Y\trsp\Lambda\ebold/n$
is the vector of column means of $\XTilde$, and if $\SigmaHat$ is the
sample covariance\vadjust{\goodbreak} matrix computed from $\XTilde$, we have
\[
\hat{m}\trsp\SigmaHat^{-1}\hat{m}\tendsto\kappa=\frac{\rho
}{1-\rho
} \qquad \mbox{in probability} .
\]
\end{theorem}

Some comments on this theorem are in order. First, $Z_{n,p}$ is
unchanged if we rescale all the $\lambda_i$'s by the same constant. So
it appears we could assume that they are all less than 1, for instance,
and dispense entirely with \eqref
{eq:AssumptionControlFourthEmpiricalMomentLambda}. However, that would
potentially violate the conditions of \eqref
{eq:AssumptionSmallestSingularValueBoundedBelow} which appear to
guarantee that $Z_{n,p}$ has variance going to zero. We also note that
because the $Y_i$'s have a continuous distribution and we know that all
the $\lambda_i$'s are different from 0, the existence of $Z_{n,p}$ is
guaranteed with probability 1.

Some practical clarifications are also in order concerning the condition
\[
\frac{1}{n^2}\sum_{i=1}^n \lambda_i^4=\frac{1}{n^2} \sum_{i=1}^n
\tau_i^2 \tendsto0 \qquad \mbox{in probability} .
\]
When the $\lambda_i$'s are i.i.d., this condition is satisfied (almost
surely and hence in probability) if for, instance, the $\lambda_i$'s
have finite second moment according to the Marcinkiewicz--Zygmund law
of large numbers [see \cite{ChowTeicherBook97}, page~125]. This is very
interesting from a practical standpoint as it basically means that we
only require our random variables $X_i$ to have a second moment for the
theorem to hold. We note that if there were no variance, the premises
of the problem would be essentially flawed (after all the quadratic
form we are optimizing involves a proxy for the population covariance,
and, in the absence of a second moment for the $\lambda_i$'s, the
population covariance would not exist), and hence we require minimal
conditions from the point of view of the practical problem at stake.

Finally, and remarkably, the limit of $Z_{n,p}$ does not depend on the
empirical distribution of the $\lambda_i$'s. In particular, in the
class of elliptical distributions (satisfying the assumptions of
Theorem~\ref{thm:simplequadformmuhatandsigmahat}), the limit of $\hat{m}
\trsp\SigmaHat^{-1}\hat{m}$ is always the same: $\kappa=\rho
/(1-\rho)$.

We now turn to proving Theorem \ref
{thm:simplequadformmuhatandsigmahat}. The proof will be facilitated by
the following lemma, which essentially gives us
$\mathbf{E}(Z_{n,p})$.\vspace*{-3pt}

\begin{lemmaNumThm}\label{lemma:MeanRandomProjection}
Let $Y$ be an $n\times p$ random matrix, with $n\geq p$ with, for
instance, independent rows, $Y_i$. Assume that $Y_i$ have symmetric
distributions, that is, $Y_i \equalInLaw-Y_i$. Let $\Lambda$ be an
$n\times n$ diagonal matrix with possibly random entries. Let
$P=\Lambda Y (Y\trsp\Lambda^2 Y)^{-1} Y\trsp\Lambda$ be a random
projection matrix. $Y$ is assumed to be independent of $\Lambda$ and
$Y$ and $\Lambda$ are assumed to be such that $P$ exists with
probability~1. Then,
\[
\mathbf{E}(\ebold\trsp P \ebold|\Lambda)=\mathbf{E}(\ebold\trsp P
\ebold)=p .
\]
\end{lemmaNumThm}

In particular, the result applies when $Y_i$ are normally distributed,
and $\Lambda$ is such that \eqref
{eq:AssumptionSmallestSingularValueBoundedBelow} holds, and $P$ is
defined with probability one.\vadjust{\goodbreak}
\begin{pf*}{Proof of   lemma~\ref{lemma:MeanRandomProjection}}
Let us note that $P=f_{\Lambda}(Y_1,\ldots,Y_n)$. Now, conditional on
$\Lambda$, $P\equalInLaw f_{\Lambda}(-Y_1,Y_2,\ldots,Y_n)=\tilde
{P}$. However, $\tilde{P}(1,j)=-P(1,j)$, if $j\neq1$. As a matter of fact,
\[
P(1,j)=\lambda_1 \lambda_j Y_1\trsp\Biggl(\sum_{i=1}^n \lambda_i^2 Y_i
Y_i\trsp\Biggr)^{-1} Y_j .
\]
Hence, conditional on $\Lambda$, $P(1,j)\equalInLaw-P(1,j)$. Now $P$
is an orthogonal projection matrix, $P=P\trsp$, so all its entries are
less than 1 in absolute value, the operator norm of $P$. In particular,
all the entries have an expectation. Since, if $j\neq1$, $P(1,j)$ has
a symmetric distribution (conditional on $\Lambda$), we conclude that
\[
\mathbf{E}(P(1,j)|\Lambda)=0 , \qquad\mbox{if } j\neq1 .
\]
Note that the same arguments would apply if $1$ were replaced by $i$,
so we really have
\[
\mathbf{E}(P(i,j)|\Lambda)=0 , \qquad\mbox{if } j\neq i .
\]
Therefore,
\[
\mathbf{E}(\ebold\trsp P \ebold|\Lambda)=\mathbf{E}(\mathrm
{trace}(P)|\Lambda )=p ,
\]
since $P$ has rank $p$ and is a projection matrix.

The same results hold when we take expectations over $\Lambda$ by
similar arguments.
\end{pf*}

To prove Theorem~\ref{thm:simplequadformmuhatandsigmahat}, all we have
to do (in light of Lemma~\ref{lemma:MeanRandomProjection}) is to show
that we control the variance of
\[
Z_{n,p}=\frac{1}{n} \ebold\trsp P \ebold.
\]
We are going to do this now by using rank 1 perturbation arguments, in
connection with the Efron--Stein inequality.
\begin{pf*}{Proof of Theorem \ref
{thm:simplequadformmuhatandsigmahat}} As before, we first work
conditionally on $\Lambda$. We assume until further notice that
$\Lambda\in{\cal L}_{\eps_0,\delta_0}$, a set of matrices which is
defined at the end of the proof, will have measure going to 1
asymptotically, and is such that all the technical issues appearing in
the proof can be taken care of. (The arguments are not circular.)

We will use the notation
\[
{\cal S}=\frac{1}{n}\sum_{k=1}^n \lambda_k^2 Y_k Y_k\trsp
\quad \mbox{and}\quad
{\cal S}_i={\cal S}-\frac{1}{n}\lambda_i^2 Y_i Y_i\trsp.
\]
Note that ${\cal S}_i$ is symmetric and positive semi-definite.
Naturally, in matrix form we can write ${\cal S}=(Y\trsp\Lambda^2
Y)/n$ and ${\cal S}_i=(Y\trsp\Lambda_i^2 Y)/n$, where $\Lambda_i^2$\vadjust{\goodbreak}
is the same matrix as $\Lambda$, except that $\Lambda_i(i,i)=0$.
Our aim is to approximate
\[
Z_{n,p}=\frac{\ebold\trsp\Lambda Y}{n} \biggl(\frac{Y\trsp\Lambda
^2 Y}{n} \biggr)^{-1}\frac{Y\trsp\Lambda\ebold}{n}=f(X_1,\ldots
,X_n) ,
\]
by a random variable involving only $(Y_1,\ldots
,Y_{i-1},Y_{i+1},\ldots,Y_n)$, that is, not involving $Y_i$. Using
classic matrix perturbation results [see \cite{hj}, page 19], we have
\[
{\cal S}^{-1}= \biggl({\cal S}_i+\frac{\lambda_i^2}{n}Y_iY_i\trsp
\biggr)^{-1}={\cal S}_i^{-1}-\frac{\lambda_i^2}{n}\frac{{\cal
S}_i^{-1}Y_iY_i\trsp{\cal S}_i^{-1}}{1+\lambda_i^2 (Y_i\trsp{\cal
S}_i^{-1}Y_i/n)} .
\]
Of course, if $e_i$ is the $i$th canonical basis vector in $\mathbb{R}^n$,
\[
W\triangleq\Lambda Y =\sum_{i=1}^n \lambda_i e_i Y_i\trsp\triangleq
W_i +\lambda_i e_i Y_i\trsp.
\]
Let us now call $q_i=Y_i\trsp{\cal S}_i^{-1} Y_i/n$ and $r_i=W_i{\cal
S}_i^{-1} Y_i$. We have
\renewcommand{\theequation}{\arabic{equation}}\setcounter{equation}{5}
\begin{equation}\label{eq:IntermedExpansionQuadFormMatrixForm}
\Lambda Y {\cal S}^{-1}=W_i{\cal S}_i^{-1}-\frac{\lambda_i^2}{n}
\frac{r_i Y_i\trsp{\cal S}_i^{-1}}{1+\lambda_i^2 q_i}+\lambda_i e_i
Y_i\trsp{\cal S}_i^{-1}-\lambda_i^3q_i \frac{e_iY_i\trsp{\cal
S}_i^{-1}}{1+\lambda_i^2 q_i} .
\end{equation}
Similarly,
\begin{eqnarray}\label{eq:keyExpansionProjectionMatrix}
\Lambda Y {\cal S}^{-1}Y\trsp\Lambda&=&W_i {\cal S}_i^{-1} W_i\trsp
-\frac{\lambda_i^2}{n}\frac{r_ir_i\trsp}{1+\lambda_i^2
q_i}+\lambda_i e_i r_i\trsp-\lambda_i^3 q_i \frac{e_i r_i\trsp
}{1+\lambda_i^2 q_i}\nonumber
\\[-8pt]
\\[-8pt]
&&{}+\lambda_i r_i e_i\trsp-\lambda_i^3 q_i \frac{r_ie_i\trsp
}{1+\lambda_i^2 q_i}+\lambda_i^2 n q_i e_i e_i\trsp-\lambda_i^4 n
q_i^2 \frac{e_ie_i\trsp}{1+\lambda_i^2 q_i}.\nonumber
\end{eqnarray}
This is, in some sense, the key expansion in this proof. Now let us
call $\muHat_i\trsp=\ebold\trsp W_i/n$ and $w_i=\ebold\trsp
r_i/n=\muHat_i\trsp{\cal S}_i^{-1}Y_i$. We have
\begin{eqnarray*}
Z_{n,p}&=&\muHat_i\trsp{\cal S}_i^{-1} \muHat_i-\frac{\lambda
_i^2}{n}\frac{w_i^2}{1+\lambda_i^2 q_i}+2\frac{\lambda_i}{n}w_i\\
&&{}-\frac{2}{n}\lambda_i^3 \frac{q_i w_i}{1+\lambda_i^2 q_i}+\frac
{\lambda_i^2}{n}q_i-\frac{\lambda_i^4}{n}\frac{q_i^2}{1+\lambda
_i^2 q_i}.
\end{eqnarray*}
Now let us call $Z_i=\muHat_i\trsp{\cal S}_i^{-1} \muHat_i$.
Clearly, $Z_i$ does not depend on $Y_i$. Now, it is easily verified that
\[
\biggl(2\lambda_i w_i +\lambda_i^2 q_i - \frac{\lambda
_i^4q_i^2}{1+\lambda_i^2 q_i}-\frac{\lambda_i^2 w_i^2}{1+\lambda
_i^2 q_i}-2 \frac{\lambda_i^3 q_i w_i}{1+\lambda_i^2 q_i}
\biggr)=1-\frac{(1-\lambda_i w_i)^2}{1+\lambda_i^2 q_i} .
\]
We finally conclude that
\begin{equation}\label{eq:ApproxZbyZi}
Z_{n,p}=Z_i+\frac{1}{n} \biggl(1-\frac{(1-\lambda_i w_i)^2}{1+\lambda
_i^2 q_i} \biggr) .
\end{equation}
We now recall the Efron--Stein inequality, as formulated in Theorem 9
of \cite{LugosiConcentration}: if $\alpha=f(X_1,\ldots,X_n)$, where
the $X_i$'s are independent, and $\alpha_i$ is a measurable function of
$(X_1,\ldots,X_{i-1},X_{i+1},\ldots,X_n)$, then
\[
\operatorname{var}(\alpha)\leq\sum_{i=1}^n \mathbf{E}\bigl((\alpha
-\alpha _i)^2\bigr) .
\]
In particular, for us, it means that
\[
\operatorname{var}(Z_{n,p}|\Lambda)\leq\sum_{i=1}^n \mathbf
{E}\biggl(\biggl(Z_{n,p}-Z_i-\frac {1}{n}\biggr)^2\Big|\Lambda\biggr) .
\]
If we now use equation \eqref{eq:ApproxZbyZi} and the fact that
$q_i\geq0$, we have
\[
n \biggl|Z_{n,p}-Z_i-\frac{1}{n} \biggr|=\frac{(1-\lambda_i
w_i)^2}{1+\lambda_i^2 q_i}\leq2(1+\lambda_i^2 w_i^2) .
\]
Moreover, conditional on $Y_{(-i)}=(Y_1,\ldots, Y_{i-1},Y_{i+1},\ldots
,Y_n)$ (and $\Lambda$ since all our arguments at this point are made
conditional on $\Lambda$), $w_i$ is ${\cal N}(0,\muHat_i\trsp{\cal
S}_i^{-2} \muHat_i)$ when the $Y$'s are ${\cal N}(0,\id_p)$,
because $w_i=\muHat_i\trsp{\cal S}_i^{-1} Y_i$. Therefore,
\[
\mathbf{E}(w_i^4|\Lambda)=3 \mathbf{E}((\muHat_i\trsp{\cal
S}_i^{-2} \muHat _i)^2|\Lambda) .
\]
Almost by definition, we have $\muHat_i\trsp{\cal S}_i^{-1} \muHat
_i\leq1$, since the vector $\ebold/\sqrt{n}$ has
norm 1 and $W_i(W_i\trsp W_i)^{-1}W_i\trsp$ is a projection matrix
(recall that ${\cal S}_i=W_i\trsp W_i/n$ and $\muHat_i\trsp=\ebold
\trsp W_i/n$). So we
would be done if we had uniform control on $\opnorm{{\cal S}_i^{-1}}$.
Let us now go around this difficulty.

$\bullet$ \textit{Regularization interlude.}

Let us consider,
for $t>0$, $Z(t)=\muHat\trsp({\cal S}+t\id_p)^{-1} \muHat$, where
$\muHat\trsp=\ebold\trsp W/n$. Clearly, $0\leq Z(t)\leq Z_{n,p}=Z(0)$,
because ${\cal S}+t\id_p\succeq{\cal S}\succeq0$ in the
positive-semidefinite ordering. In other respects, the decomposition
in equation \eqref{eq:ApproxZbyZi} is still valid if we replace $Z_i$
by $Z_i(t)$ and ${\cal S}_i$ by ${\cal S}_i(t)$ everywhere. However,
$\opnorm{({\cal S}_i(t))^{-1}} \leq1/t$. We therefore have
\[
\muHat_i\trsp{\cal S}_i(t)^{-2} \muHat_i\leq\opnorm{{\cal
S}_i^{-1}(t)}\| {\cal S}^{-1/2}_i\muHat_i\|_2^2\leq\frac{\muHat
_i\trsp{\cal S}_i^{-1}(t)\muHat_i}{t} \leq\frac{\muHat_i\trsp
{\cal S}_i^{-1}\muHat_i}{t} \leq\frac{1}{t}.
\]
So applying the previous analysis and using the fact that $\muHat
_i\trsp({\cal S}_i(t))^{-2}\muHat_i\leq1/t$, we conclude that
\[
\operatorname{var}(Z(t)|\Lambda)\leq\frac{8}{n^2}\sum_{i=1}^n \biggl(1+3
\frac
{\lambda_i^4}{t^2}\biggr) .
\]
So under our assumptions, $Z(t)$ can be approximated, in probability,
at least conditionally on $\Lambda$, by $\mathbf{E}(Z(t)|\Lambda)$.
If we
write the singular value decomposition of $W/\sqrt{n}=\sum_{i=1}^p
\sigma_i u_i v_i\trsp$, where $\sigma_1\geq\sigma_2\geq\cdots
\geq\sigma_p$, we have
$W{\cal S}^{-1}W\trsp/n=\sum_{i=1}^p u_iu_i\trsp$, $W({\cal
S}(t))^{-1}W\trsp/n=\sum_{i=1}^p \sigma_i^2/(\sigma_i^2+t)
u_iu_i\trsp$, and therefore
\begin{eqnarray*}
0&\leq& Z_{n,p}-Z(t)=\frac{t}{n} \sum_{i=1}^p \frac{1}{\sigma_i^2 +t}
(u_i\trsp\ebold)^2\\
 &\leq&\frac{t}{\sigma_p^2+t} \frac{1}{n}\sum
_{i=1}^p (u_i\trsp\ebold)^2 \leq\frac{t}{\sigma_p^2+t} \frac
{\| \ebold\|_2^2}{n}=\frac{t}{\sigma_p^2+t} .
\end{eqnarray*}
To get the inequality above, we used the fact that the $\{u_i\}
_{i=1}^p$ are orthonormal in $\mathbb{R}^n$, and can therefore be
completed to form an orthonormal basis of this vector space. The
quantities $u_i\trsp\ebold$ are
naturally the coefficients of $\ebold$ in this basis, and we know that
their sum of squares should be the squared norm of $\ebold$, which is $n$.

Let us now call ${\cal L}_{\eps_0,\delta}$ the set of matrices
$\Lambda$ such that $p/N<1-\eps_0$ and $C_0 (N-1)/(n-1)>\delta$.
Under our assumptions, for a $\delta_0$ bounded away
from 0 (e.g., $\delta_0=1/2 C_0 \eta_0$), $P(\Lambda\in{\cal
L}_{\eps_0,\delta_0})\tendsto1$. Let us pick such a $\delta_0$.
If $\Lambda\in{\cal L}_{\eps_0,\delta_0}$, according to Lemma \ref
{lemma:LowerBoundOnSmallestEig} and the proof of Theorem \ref
{thm:QuadFormsInverseEllipticalCase}, if $P_{\Lambda}$ denotes
probability conditional on $\Lambda$,
\[
P_{\Lambda} \bigl(\sigma_p\leq\sqrt{\delta_0} \bigl[\bigl(1-\sqrt
{1-\eps_0}\bigr)-t \bigr] \bigr)\leq\exp\bigl(-(n-1)\delta_0
t^2/C_0 \bigr) .
\]
Hence, when $\Lambda\in{\cal L}_{\eps_0,\delta_0}$, we can find,
for any $u>0$, an $\eta(u)>0$,
\[
P\bigl(|Z_{n,p}-Z(\eta(u))|>u\bigr)\leq f_n(\eps_0,\delta_0,\eta
(u),u)=f_n(u) ,
\]
where, $f_n(u)=f_n(\eps_0,\delta_0,\eta(u),u)\tendsto0$ as
$n\tendsto\infty$, for fixed $u$.

On the other hand, our conditional variance computations have
established that, for any $\eta>0$, $Z(\eta)-\mathbf{E}(Z(\eta
)|\Lambda )$ converges in probability
(conditional on $\Lambda$) to 0 if $\eta^{-2} \sum\lambda
_i^4/n^{2}$ tends to 0. We note that $0 \leq Z_{n,p}\leq1$ and that
the same is true for $\gamma_n(u)=\mathbf{E}(Z(\eta(u))|\Lambda)$.
Therefore, $|Z_{n,p}-\gamma_n(u)|\leq1$ and $\mathbf
{E}((Z_{n,p}-\gamma _n(u))^2|\Lambda)$ goes to zero, since
\begin{eqnarray*}
&&\mathbf{E}\bigl(\bigl(Z-\gamma_n(u)\bigr)^2|\Lambda\bigr)\\
&&\qquad\leq u^2 P\bigl(|Z_{n,p}-\gamma
_n(u)|\leq u|\Lambda\bigr)+ P\bigl(|Z_{n,p}-\gamma_n(u)|>u|\Lambda\bigr)\\
&&\qquad\leq u^2 +P\bigl(|Z_{n,p}-Z(\eta(u))|>u/2|\Lambda\bigr)+\frac
{4}{u^2}\operatorname{var}(Z(\eta(u))|\Lambda) .
\end{eqnarray*}
In other words, we also have, if $\Lambda\in{\cal L}_{\eps_0,\delta
_0}$, for any $u>0$,
\[
\operatorname{var}(Z_{n,p}|\Lambda)\leq u^2+f_n(u/2)+\frac
{32}{u^2}\frac
{1}{n^2}\sum_{i=1}^n \biggl(1+3\frac{\lambda_i^4}{\eta(u)^2}
\biggr) .
\]

Hence, if $\Lambda\in{\cal L}_{\eps_0,\delta_0}$ and $\sum_{i=1}^n
\lambda_i^4/n^2\tendsto0$, $\operatorname{var}(Z|\Lambda)$ goes to
zero as $n$
goes to infinity, and we conclude that, since $\mathbf{E}(Z|\Lambda)=p/n$,
\[
Z-\frac{p}{n} \tendsto0 \qquad \mbox{in probability, conditional on }
\Lambda.
\]

$\bullet$ \textit{Deconditioning on $\Lambda$.}

Let us call ${\cal L}^2_{\eps_0,\delta_0,t}$ the set of matrices such
that ${\cal L}^2_{\eps_0,\delta_0,t}={\cal L}_{\eps_0,\delta
_0}\cap\{(\frac{1}{n^2}\times \sum_{i=1}^n \lambda_i^4)\leq t\}$. Our
previous computations clearly show that we can find a function
$g_n(u)$, with $g_n(u)\tendsto0$ as $n\tendsto\infty$, such that,
for any $u>0$, when $\Lambda\in{\cal L}^2_{\eps_0,\delta_0,u^4\eta
(u)^2}\triangleq{\cal L}^2(u)$, $\operatorname{var}(Z_{n,p}|\Lambda
)\leq97 u^2
+g_n(u)$, and hence we have the ``uniform bound,'' if $\Lambda\in
{\cal
L}^2(u)$,
\[
P\biggl( \biggl|Z_{n,p}-\frac{p}{n} \biggr|>x|\Lambda\biggr)\leq\frac{97
u^2+g_n(u)}{x^2} .
\]
Now under our assumptions, $P(\Lambda\in{\cal L}^2(u))$ goes to 1 for
any given $u$, so we conclude, using the fact that
\begin{eqnarray*}
P(|Z_{n,p}-p/n|>x)&\leq& P[|Z_{n,p}-p/n|>x\ \&\ \Lambda\in{\cal
L}^2(u)]\\
&&{}+P[\Lambda\notin{\cal L}^2(u)] ,
\end{eqnarray*}
that
\[
Z_{n,p}-\frac{p}{n}\tendsto0 \qquad \mbox{in probability}.
\]
This last statement is now understood of course unconditionally on
$\Lambda$ and this proves the first part of the theorem.

$\bullet$ \textit{Proof of the second part of the theorem.}

We now focus on the $\hat{m}\trsp\SigmaHat^{-1}\hat{m}$ part of the
theorem. Let us call $\mathfrak{S}=\XTilde\trsp\XTilde/n$. Then,
$\frac{n-1}{n}\SigmaHat=\mathfrak{S}-\hat{m}\hat{m}\trsp$. Therefore,
\[
\frac{n}{n-1}\SigmaHat^{-1}=\mathfrak{S}^{-1}+\frac{\mathfrak
{S}^{-1}\hat{m}\hat{m}\trsp\mathfrak{S}^{-1}}{1-\hat{m}\trsp
\mathfrak
{S}^{-1}\hat{m}} .
\]
Hence,
\[
\frac{n}{n-1}\hat{m}\trsp\SigmaHat^{-1} \hat{m}=\frac{\hat
{m}\trsp
\mathfrak{S}^{-1}\hat{m}}{1-\hat{m}\trsp\mathfrak{S}^{-1}\hat
{m}}=\frac
{Z_{n,p}}{1-Z_{n,p}} .
\]
Since $Z_{n,p}\tendsto\rho$ in probability with $\rho\in(0,1)$, we
have the result announced in the theorem.
\end{pf*}

Now that we have proved Theorem \ref
{thm:simplequadformmuhatandsigmahat}, we need to turn to results that
will allow us to handle the case of nonzero population mean,
as well as questions such as the convergence of $\muHat\trsp\SigmaHat
^{-1} v$, for deterministic $v$.

\subsubsection{\texorpdfstring{On quantities of the type ${(\muHat-\mu)\trsp\SigmaHat^{-1}\mu}$}
{On quantities of the type (mu - mu)' Sigma -1 mu}}

Recall that the key quantity in the solution of problem \eqref
{eq:QP-eqc-Emp}, the problem of main interest in this paper, is of the
form $\VHat\trsp\SigmaHat^{-1} \VHat$.
Therefore, it is important for us to understand quantities of the type
\[
\zeta=\muHat\trsp\SigmaHat^{-1} v ,
\]
for a fixed vector $v$. At this point, we focus on the particular case
where $\mu=\mathbf{E}(X_i)=0$. To do so, we will need to study, if
${\cal
S}=Y\trsp\Lambda\trsp\Lambda Y/n$,
\[
\zeta=\frac{1}{n}\ebold\trsp\Lambda Y {\cal S}^{-1} v ,
\]
for a fixed vector $v$.
As it turns out, this random variable goes to zero in probability when
for instance $\| v\|_2=1$.
\begin{theorem}\label{thm:QuadFormsMuHatSigmaHatInvV}
Suppose $v$ is a deterministic vector, with $\| v\|_2=1$. Suppose
the assumptions stated in Theorem \ref
{thm:simplequadformmuhatandsigmahat} hold and also that
\renewcommand{\theequation}{Assumption-BLb}
\begin{equation}\label{eq:ConditionFiniteSecondEmpiricalMomentLambda}
\frac{1}{n}\sum_{i=1}^n \lambda_i^2 \mbox{ remains bounded with
probability going to 1.}\hspace*{-100pt}
\end{equation}
Consider
\[
\zeta=\frac{1}{n}\ebold\trsp\Lambda Y {\cal S}^{-1} v ,
\]
where ${\cal S}=\frac{1}{n}Y\trsp\Lambda^2 Y$.
Then
\[
\zeta\tendsto0 \qquad \mbox{in probability} .
\]
\end{theorem}

Before giving the proof, we note that if the $\lambda_i$'s are i.i.d.
and have a second moment, the ``extra'' condition on $\sum_{i=1}^n
\lambda_i^2/n$ introduced in this theorem (as compared to Theorem \ref
{thm:simplequadformmuhatandsigmahat}) is clearly satisfied by the law
of large numbers.

\begin{pf*}{Proof of Theorem~\ref{thm:QuadFormsMuHatSigmaHatInvV}}
The proof is quite similar to the proof of Theorem~\ref
{thm:simplequadformmuhatandsigmahat} above. We start by conditioning on
$\Lambda$.

Let us call $\zeta(t)$ the quantity obtained when we replace ${\cal
S}$ by ${\cal S}(t)={\cal S}+t\id$ in the definition of $\zeta$. Note
that since $Y$ is symmetric, $\zeta(t)\equalInLaw-\zeta(t)$,
conditionally on $\Lambda$, by arguments similar to those given in the
proof of Lemma~\ref{lemma:MeanRandomProjection}. Now $\zeta(t)$
clearly has an expectation (conditional on $\Lambda$), because
$\opnorm{S^{-1}(t)}\leq1/t$, for $t>0$, so $\mathbf{E}(\zeta
(t)|\Lambda )=0$. Now recall equation \eqref
{eq:IntermedExpansionQuadFormMatrixForm}: with the notations used there,
\[
\Lambda Y {\cal S}^{-1}=W_i{\cal S}_i^{-1}-\frac{\lambda_i^2}{n}
\frac{r_i Y_i\trsp{\cal S}_i^{-1}}{1+\lambda_i^2 q_i}+\lambda_i e_i
Y_i\trsp{\cal S}_i^{-1}-\lambda_i^3q_i \frac{e_iY_i\trsp{\cal
S}_i^{-1}}{1+\lambda_i^2 q_i} .
\]
Let us now call $q_i(t)=Y_i\trsp{\cal S}_i(t)^{-1} Y_i/n$,
$w_i(t)=\ebold\trsp W_i {\cal S}_i(t)^{-1} Y_i/n=\muHat_i\trsp{\cal
S}_i(t)^{-1} Y_i$ and $\theta_i(t)=Y_i\trsp{\cal S}_i(t)^{-1} v$.
Clearly, if $\zeta_i(t)$ is the random variable obtained by excluding
$Y_i$ from the computation of $\zeta(t)$ (e.g., by replacing $\lambda
_i$ by 0), we have
\begin{eqnarray*}
\zeta(t)&=&\zeta_i(t)-\frac{\lambda_i^2}{n}\frac{w_i(t)\theta
_i(t)}{1+\lambda_i^2 q_i(t)} +\frac{\lambda_i \theta_i(t)}{n}-\frac
{\theta_i(t)}{n}\frac{\lambda_i^3 q_i(t)}{1+\lambda_i^2 q_i(t)}\\
&=&\zeta_i(t)+\frac{1}{n} \biggl(\frac{\lambda_i \theta
_i(t)(1-\lambda_i w_i(t))}{1+\lambda_i^2 q_i(t)} \biggr).
\end{eqnarray*}
We remark that $\theta_i(t)|(Y_{(-i)},\Lambda)\sim{\cal N}(0,v\trsp
{\cal S}_i^{-2}(t) v)$ and recall that $w_i|(Y_{(-i)},\break\Lambda)\sim
{\cal N}(0,\muHat_i\trsp{\cal S}_i^{-2}\muHat_i)$. Using the fact
that $\| v\|_2=1$, $\opnorm{{\cal S}_i^{-2}(t)}\leq t^{-2}$ and
the remarks we made in the proof of Theorem \ref
{thm:simplequadformmuhatandsigmahat}, we get
that $\mathbf{E}([\theta_i(t)]^{2k}|(Y_{(-i)},\break\Lambda))\leq C_k
t^{-2k}$, $\mathbf{E}([w_i(t)]^{2k}|(Y_{(-i)},\Lambda))\leq C_k t^{-k}$,
where $C_1=1$ and $C_2=3$. We also have
\[
\bigl[\lambda_i \theta_i(t)\bigl(1-\lambda_i w_i(t)\bigr) \bigr]^2\leq
2 [\lambda_i^2\theta_i^2(t)+\lambda_i^4 \theta
_i^2(t)w_i^2(t) ].
\]
Hence, simply using the fact that $2(ab)^2\leq(a^4+b^4)$, we get
\[
\mathbf{E}\biggl(\biggl(\frac{\lambda_i \theta_i(t)(1-\lambda_i
w_i(t))}{1+\lambda_i^2 q_i(t)} \biggr)^2\Big|\Lambda\biggr)\leq2 \frac
{2\lambda_i^2}{t^2}+3 \lambda_i^4 \biggl(\frac{1}{t^2}+\frac
{1}{t^4} \biggr) .
\]

We conclude by the Efron--Stein inequality that, when $\Lambda$ is
such that $\sum_{i=1}^n \lambda_i^4/ n^2\tendsto0$, for any $t>0$,
\[
\zeta(t)\tendsto0 \qquad \mbox{in probability, conditionally on } \Lambda.
\]

As before, let us call ${\cal L}_{\eps_0,\delta}$ the set of matrices
$\Lambda$ such that $p/N<1-\eps_0$ and $C_0 (N-1)/(n-1)>\delta$.
Recall that under our assumptions, for $\delta_0$ bounded away from 0
(e.g., $\delta_0=C_0\eta_0/2$), $P(\Lambda\in{\cal L}_{\eps
_0,\delta_0})\tendsto1$.

As we saw before, when $\Lambda\in{\cal L}_{\eps_0,\delta_0}$,
$\opnorm{{\cal S}^{-1}}$ is bounded with high-probabil\-ity (conditional
on $\Lambda$), so we conclude that, for any $\eta>0$, we can find a
$t$ such that
\[
\opnorm{{\cal S}^{-1}-{\cal S}^{-1}(t)}<\eta\qquad\mbox{with
probability (conditional on $\Lambda$) going to }1 .
\]
We also notice that conditionally on $\Lambda$, $\muHat\sim{\cal
N} (0,\frac{\sum\lambda_i^2}{n^2} \id_p )$ and hence,
$\| \muHat\|_2^2 \sim\chi^2_p/n (\sum\lambda_i^2)/n$.
We recall that $\| v\|_2=1$, and since
\[
|\zeta-\zeta(t) |\leq\| \muHat\|_2\opnorm{{\cal
S}^{-1}-{\cal S}^{-1}(t)} \| v\|_2 ,
\]
we conclude that with high-probability (conditional on $\Lambda$), for
any $\eta>0$, $ |\zeta-\zeta(t) |\leq\eta$ and finally,
\[
\zeta\tendsto0 \qquad \mbox{in probability, conditionally on }\Lambda.
\]
Now along the same lines as what was done in the proof of Theorem \ref
{thm:simplequadformmuhatandsigmahat}, we can make all these probability
bounds uniform in $\Lambda$ when $\Lambda$ is in
a set of matrices such as ${\cal L}_{\eps_0,\delta_0}$ and when we
also have bounds on $\sum_{i=1}^n \lambda_i^4/n^2$ and $\sum_{i=1}^n
\lambda_i^2/n$. Under our assumptions, the set of $\Lambda$ for which
these conditions hold has measure going to 1, so we can finally
conclude---along the same lines (omitted here) as in the proof of
Theorem~\ref{thm:simplequadformmuhatandsigmahat}---that,
unconditionally on $\Lambda$,\looseness=-1
\[
\zeta\tendsto0 \qquad \mbox{in probability} .
\]\looseness=0
\upqed
\end{pf*}

After these preliminaries, we can finally state the theorem of main
interest. Recall that under the assumptions of Theorem \ref
{thm:QuadFormsInverseEllipticalCase}, if $v$ is deterministic,
\[
\frac{v\trsp\SigmaHat^{-1}v}{v\trsp\Sigma^{-1} v}\tendsto
\limitScaling\qquad \mbox{in probability} ,
\]
where $\limitScaling$ is defined in equation
\eqref{eq:defLimitScaling}.\vadjust{\goodbreak}
\begin{theorem}\label{thm:SummaryQuadFormsMuHatAndSigmaHat}
Suppose that $X_i=\mu+\lambda_i \Sigma^{1/2}Y_i$, where $Y_i$ are
i.i.d.\break ${\cal N}(0, \id_p)$ and $\{\lambda_i\}_{i=1}^n$ are random
variables, independent of $\{Y_i\}_{i=1}^n$. Let $v$ be a deterministic
vector. Suppose that $\rho_n=p/n$ has a finite nonzero limit, $\rho$
and that $\rho\in(0,1)$.

We call $\tau_i=\lambda_i^2$. We assume that $\tau_i\neq0$ for all
$i$ as well as
\renewcommand{\theequation}{Assumption-BL}
\begin{eqnarray}\label{eq:AssumptionSecondAndFourthMomentLambda}
\begin{tabular}{l}
\mbox{$\frac{1}{n^2}\sum_{i=1}^n \lambda_i^4\tendsto0$ in
probability and}\\[4pt]$\frac{1}{n}\sum_{i=1}^n \lambda_i^2 \mbox{
remains bounded in probability.}$\hspace*{-50pt}
\end{tabular}
\end{eqnarray}
If $\tau_{(i)}$ is the $i$th largest $\tau_k$, we assume that we can
find a random variable $N \in\mathbb{N}$ and positive real numbers
$\eps_0$ and $C_0$ such that
\renewcommand{\theequation}{Assumption-BB}
\begin{equation}
\cases{\displaystyle
P(p/N<1-\eps_0)\tendsto1,  \qquad\mbox{as } n\tendsto\infty ,\cr\displaystyle
P\bigl(\tau_{(N)}>C_0\bigr)\tendsto1 ,\cr\displaystyle
\exists\eta_0>0   \quad \mbox{such that}\quad   P(N/n>\eta_0)\tendsto1,\qquad \mbox{as } n\tendsto\infty.
}\hspace*{-100pt}
\end{equation}

We also assume that the empirical distribution of $\tau_i$'s converges
weakly in probability to a deterministic limit $G$.

We call $\Lambda$ the $n\times n$ diagonal matrix with $\Lambda
(i,i)=\lambda_i$, $Y$ the $n\times p$ matrix whose $i$th row is $Y_i$,
$W=\Lambda Y$ and ${\cal S}=W\trsp W/n=\sum_{k=1}^n \lambda_k^2
Y_kY_k\trsp/n$.
Finally, we use the notation $\omegahat=W\trsp\ebold/n$, $\muTilde
=\Sigma^{-1/2}\mu$.

Then, we have, for $\limitScaling$ defined as in equation \eqref
{eq:defLimitScaling},
\renewcommand{\theequation}{\arabic{equation}}\setcounter{equation}{8}
\begin{equation}\label{eq:HalfQuadFormMuHatSigmaHatv}\quad
\frac{\muHat\trsp\SigmaHat^{-1}v}{\sqrt{v\trsp\Sigma^{-1}
v}}=\frac{\mu\trsp\SigmaHat^{-1}v}{\sqrt{v\trsp\Sigma^{-1}
v}}+\lo_P(1)=\limitScaling\frac{\mu\trsp\Sigma^{-1}v}{\sqrt
{v\trsp\Sigma^{-1} v}}+\lo_P \biggl(1\vee\frac{\mu\trsp\Sigma
^{-1}v}{\sqrt{v\trsp\Sigma^{-1} v}} \biggr) ,
\end{equation}
the second statement holding if, for instance, $\mu$ and $v$ are such
that the first set of conditions in Lemma \ref
{lemma:DotProductsInvolvingSigmaInverse} are met.

Also,
\begin{equation}\label{eq:QuadFormMuHatSigmaHatMuHat}
\muHat\trsp\SigmaHat^{-1}\muHat=\mu\trsp\SigmaHat^{-1}\mu+\frac
{\rho_n}{1-\rho_n}+2\frac{n-1}{n}\frac{\omegahat\trsp{\cal
S}^{-1} \muTilde}{1-\omegahat\trsp{\cal S}^{-1}\omegahat} +\lo_P(1),
\end{equation}
and we recall that $\omegahat\trsp{\cal S}^{-1} \muTilde/\| \muTilde
\|=\lo_P(1)$ and $\omegahat\trsp{\cal S}^{-1} \omegahat
=p/n+\lo_P(1)$.
\end{theorem}

To be able to exploit equation \eqref{eq:QuadFormMuHatSigmaHatMuHat}
in practice, we make the following remarks. We can consider three
cases, having to do with the size of $\mu\trsp\Sigma^{-1} \mu=\|
\muTilde\|_2^2$:
\begin{enumerate}[3.]
\item If $\mu\trsp\Sigma^{-1} \mu\tendsto0$, then, $\muHat\trsp
\SigmaHat^{-1}\muHat=\frac{\rho_n}{1-\rho_n}+\lo_P(1)$.
\item If $\mu\trsp\Sigma^{-1} \mu\tendsto\infty$, then $\muHat
\trsp\SigmaHat^{-1}\muHat\sim\limitScaling\mu\trsp\Sigma^{-1}
\mu$.
\item Finally, if $\mu\trsp\Sigma^{-1} \mu$ stays bounded away from
0 and infinity,
\[
\muHat\trsp\SigmaHat^{-1}\muHat= \limitScaling\mu\trsp\Sigma
^{-1} \mu+\frac{\rho_n}{1-\rho_n}+\lo_P(1) .
\]
\end{enumerate}
A noticeable feature of these results is that the ``extra bias''
$\kappa
_n=\rho_n/(1-\rho_n)$, which comes essentially from mis-estimation of
$\mu$, is constant within the class of elliptical\vadjust{\goodbreak} distributions
considered here. This should be contrasted with the ``scaling,''
$\limitScaling$, which strongly depends on the empirical distribution
of the $\lambda_i^2$'s.

We now give a brief proof of Theorem \ref
{thm:SummaryQuadFormsMuHatAndSigmaHat}.
\begin{pf*}{Proof of Theorem
\ref{thm:SummaryQuadFormsMuHatAndSigmaHat}}
We first note that $\Sigma^{1/2}\omegahat=\hat{m}$ in the notation of
Theorem~\ref{thm:simplequadformmuhatandsigmahat}. Also, $\muHat=\mu
+\Sigma^{1/2}\omegahat=\mu+\hat{m}$. Finally,
\[
\frac{n-1}{n}\SigmaHat=\Sigma^{1/2}{\cal S} \Sigma^{1/2} - \hat{m}
\hat{m}\trsp=\Sigma^{1/2} ({\cal S}-\omegahat\omegahat\trsp
)\Sigma^{1/2} .
\]

\textit{Proof of equation \eqref{eq:QuadFormMuHatSigmaHatMuHat}}.
By writing $\muHat=\mu+\hat{m}$, we
clearly have
\[
\muHat\trsp\SigmaHat^{-1}\muHat=\mu\trsp\SigmaHat^{-1}\mu
+2\hat{m}\trsp\SigmaHat^{-1}\mu+\hat{m}\trsp\SigmaHat^{-1}\hat{m}.
\]
We have already seen in Theorem \ref
{thm:simplequadformmuhatandsigmahat} that the third term tends to
$\kappa=\rho/(1-\rho)$. On the other hand, half of the middle term
is equal to
\[
\frac{n}{n-1}\omegahat\trsp({\cal S}-\omegahat\omegahat\trsp
)^{-1}\muTilde.
\]
Since $ ({\cal S}-\omegahat\omegahat\trsp)^{-1}={\cal
S}^{-1}+{\cal S}^{-1}\omegahat\omegahat\trsp{\cal
S}^{-1}/(1-\omegahat\trsp{\cal S}^{-1}\omegahat)$, we have
\begin{eqnarray*}
\frac{n}{n-1}\hat{m}\trsp\SigmaHat^{-1}&=&\omegahat\trsp{\cal
S}^{-1} \biggl(1+\frac{\omegahat\trsp{\cal S}^{-1}\omegahat
}{1-\omegahat\trsp{\cal S}^{-1}\omegahat} \biggr)\Sigma
^{-1/2}\\
&=&\frac{1}{1-\omegahat\trsp{\cal S}^{-1}\omegahat}\omegahat
\trsp{\cal S}^{-1}\Sigma^{-1/2} ,
\end{eqnarray*}
and we deduce the result of equation \eqref
{eq:QuadFormMuHatSigmaHatMuHat}. We now remark that $\omegahat\trsp
{\cal S}^{-1}\omegahat$ is equal to the quantity $Z_{n,p}$ in Theorem
\ref{thm:simplequadformmuhatandsigmahat}. The fact that $\omegahat
\trsp{\cal S}^{-1} \muTilde/\| \muTilde\|=\lo_P(1)$ follows from
applying Theorem~\ref{thm:QuadFormsMuHatSigmaHatInvV} with $v=\muTilde
/\| \muTilde\|_2$.

\textit{Proof of equation \eqref
{eq:HalfQuadFormMuHatSigmaHatv}}.
 The proof of this result follows from
a decomposition similar to the one we just made. Clearly the only
question is whether $\hat{m}\trsp\SigmaHat^{-1}v/\sqrt{v\trsp
\Sigma
^{-1}v}$ goes to 0. As we just saw,
\[
\frac{n}{n-1}\hat{m}\trsp\SigmaHat^{-1}v=\frac{1}{1-\omegahat
\trsp
{\cal S}^{-1}\omegahat}\omegahat\trsp{\cal S}^{-1}\Sigma^{-1/2}v .
\]
The results of Theorem~\ref{thm:QuadFormsMuHatSigmaHatInvV} guarantee that
\[
\frac{\omegahat\trsp{\cal S}^{-1}\Sigma^{-1/2}v}{\| \Sigma
^{-1/2}v\|_2}\tendsto0 \qquad \mbox{in probability} .
\]
Since $\omegahat\trsp{\cal S}^{-1}\omegahat$ tends to $\rho<1$ and
$\| \Sigma^{-1/2}v\|_2^2=v\trsp\Sigma^{-1} v$, we have shown the
result stated in equation \eqref{eq:HalfQuadFormMuHatSigmaHatv}.
\end{pf*}

\renewcommand{\theequation}{\arabic{equation}}\setcounter{equation}{10}
\subsection{On the effect of correlation between observations}\label
{Sec:Elliptical:Subsec:CorrelObs}
It is clear that in financial practice and other applied settings, the
assumption that the returns (or observed data vectors) are independent
is often questionable. So for quadratic programs
with linear equality constraints (including the Markowitz problem but
also going beyond it), it is natural to ask what is the impact of
correlation in our observations on the empirical solution of the
problem. In our notation, this means that the vectors $X_i$\vadjust{\goodbreak} and $X_j$
are correlated; we refer to this situation as the correlated case
or as the case of temporal correlation.

Our work on the elliptical case comes in handy here and allows us to
also draw conclusions concerning the correlated case. We consider a
particular model, namely we assume that the $n \times p$ data matrix
$X$ is given by
\[
X=\ebold_n\mu\trsp+\Lambda Y \Sigma^{1/2} ,
\]
 where $\Lambda$ is a deterministic but \textit{not} necessarily a
diagonal matrix,
and $Y$ is a matrix with i.i.d. ${\cal N}(0,1)$ entries. We assume
throughout that $\Lambda$ is full rank. The model we consider now is
more general than the one we looked at before, since if $\Lambda=\id
_n$, we get the i.i.d. Gaussian case, and if $\Lambda$ is diagonal we
are back in an ``elliptical'' case (where the ellipticity parameters are
assumed to be deterministic, which amounts to doing computations
conditional on $\Lambda$). But when $\Lambda$ is not diagonal, $X_i$
and $X_j$ might be correlated. [In all the situations where $\Lambda$
is deterministic, the marginal distribution of $X_i$ is ${\cal N}(\mu,
s_i^2 \Sigma)$,
where $s_i$ is the norm of the $i$th row of $\Lambda$.]

Because we want to focus here on robustness questions arising when
going from independent Gaussian random variables to correlated ones, we
will assume throughout that $\Lambda$ is deterministic. (Allowing
$\Lambda$ to be random simply requires some minor technical
modifications but would make the exposition a bit less clear.) Our main
results in this subsection can be interpreted as saying that
that the Gaussian analysis of Section~\ref{sec:GaussianCase}, carried
out in the setting of independent observations, is not robust against
these independence assumptions. The results
change quite significantly when the vectors of observations are correlated.

In general, we write the singular value decomposition of the $n\times
n$ matrix $\Lambda$ as $\Lambda=ADB\trsp$ [see \cite{hj}, page
414], where $A$ and $B$ are orthogonal, and $D$ is diagonal. Therefore,
$AA\trsp=\id_n$, and
\[
\frac{1}{n}(X-\ebold_n\mu\trsp)\trsp(X-\ebold_n\mu\trsp)=\frac
{1}{n}\Sigma^{1/2} Y\trsp B D^2 B\trsp Y \Sigma^{1/2}\equalInLaw
\frac{1}{n}\Sigma^{1/2} Y\trsp D^2 Y \Sigma^{1/2} .
\]
So we are almost back in the elliptical case. The key difference now is
that what will matter in our analysis
are not the diagonal entries of $\Lambda\trsp\Lambda$, but rather
its eigenvalues (see Proposition \ref
{prop:QuadFormsInverseEllipticalCaseCorrelated}). Also, we will see (in
Proposition~\ref{proposition:QuadFormsMuHatSigmaInverseCorrCase}) that
the results change quite significantly when we look at quantities like
$\muHat\trsp\SigmaHat^{-1}\muHat$.

\subsubsection{\texorpdfstring{On quadratic forms involving ${\SigmaHat^{-1}}$}
{On quadratic forms involving Sigma -1}}

As a counterpart to Theorem~\ref{thm:QuadFormsInverseEllipticalCase},
we have the following proposition.
\begin{propositionNumThm}\label{prop:QuadFormsInverseEllipticalCaseCorrelated}
Suppose the $n\times p$ data matrix $X$ (whose $i$th row is the $i$th
vector of observations) can be written as
\[
X=\ebold_n\mu\trsp+\Lambda Y \Sigma^{1/2} ,
\]
  where
$\Lambda$  is a deterministic but \textit{not} necessarily
diagonal matrix. Suppose that the eigenvalues of $\Lambda\trsp\Lambda$ satisfy \eqref
{eq:AssumptionSmallestSingularValueBoundedBelow} with a deterministic
$N$ and that the spectral distribution of $\Lambda\trsp\Lambda$
converges weakly to a probability distribution $G$. Suppose also that
$p/n\tendsto\rho\in(0,1)$. Call $\SigmaHat$ the classical sample
covariance matrix, that is,
\[
\SigmaHat=\frac{1}{n}(X-\bar{X})\trsp(X-\bar{X}) .
\]
Then, if $v$ is a deterministic vector, we have
\[
\frac{v\trsp\SigmaHat^{-1} v}{v\trsp\Sigma^{-1} v}\tendsto
\limitScaling\qquad \mbox{in probability} ,
\]
where $\limitScaling$ satisfies, if $G$ is the limiting spectral
distribution $\Lambda\trsp\Lambda$
\[
\int\frac{dG(\tau)}{1+\rho\tau\limitScaling}=1-\rho.
\]
\end{propositionNumThm}

The proposition shows that Theorem \ref
{thm:QuadFormsInverseEllipticalCase} essentially applies again;
however, now what matters, unsurprisingly, are the singular
values of $\Lambda$ and not its diagonal entries. The proof of
Proposition~\ref{prop:QuadFormsInverseEllipticalCaseCorrelated}, or
rather the adjustments needed to make the proof of Theorem \ref
{thm:QuadFormsInverseEllipticalCase} go through, are given in the
\hyperref[appm]{Appendix}, Section~\ref{App:Sec:Generalizations:Subsec:QuadFormsV}.

\subsubsection{\texorpdfstring{On quadratic forms involving $\muHat$ and $\SigmaHat^{-1}$}
{On quadratic forms involving mu and Sigma -1}}

This is the situation where the results are most different from that of
the uncorrelated case. Once again, here we will be content to just
state the results; a detailed justification of our claims is in the
\hyperref[appm]{Appendix}, Section~\ref{App:Subsec:GeneralQuadFormsProjMat}.

As before, the most complicated aspect of the problem is to understand
quantities of the type $\muHat\trsp\SigmaHat^{-1}\muHat$, in the
situation where $\mu=0$. In this setting, we have the following result.
\begin{propositionNumThm}\label{proposition:QuadFormsMuHatSigmaInverseCorrCase}
Suppose the $n\times p$ data matrix $\XTilde$ is such that, for $Y$ an
$n\times p$ matrix with i.i.d. ${\cal N}(0,1)$ entries, and $\Lambda$
a deterministic matrix,
\[
\XTilde=\Lambda Y \Sigma^{1/2} .
\]
We assume that \eqref{eq:AssumptionSmallestSingularValueBoundedBelow}
holds for the eigenvalues of $\Lambda\trsp\Lambda$, for a
deterministic sequence $N(n)$.
We write the singular value decomposition of $\Lambda$ as $\Lambda
=ADB\trsp$.

We call ${\cal S}=\XTilde\trsp\XTilde/n$ and $\hat{m}=\Sigma
^{1/2}Y\trsp\Lambda\trsp\ebold/n$, that is, the sample mean of the
columns of $\XTilde$. We denote by $d_i$ the diagonal elements of $D$,
and $\widetilde{Y}=B\trsp Y\equalInLaw Y$. We also call
\[
F=\frac{1}{n}\sum_{i=1}^n d_i^2 \widetilde{Y}_i \widetilde
{Y}_i\trsp,\qquad F_i=F-\frac{1}{n} d_i^2 \widetilde{Y}_i \widetilde
{Y}_i\trsp,\qquad P=D\YTilde(\YTilde\trsp D^2 \YTilde
)^{-1}\YTilde\trsp D .
\]
If we call $\omega=A\trsp\ebold$, and $q_i=\YTilde_i\trsp F_i^{-1}
\YTilde_i/n$, we have, if $\| \omega\|_4^4/n^2$ and \mbox{$\| d\|
_4^4/n^2\tendsto0$},
\[
\hat{m}\trsp{\cal S}^{-1} \hat{m}-\kappa(n,p) \tendsto0 \qquad \mbox{in probability,}
\]
where
\[
\kappa(n,p)=\frac{1}{n}\sum_{i=1}^n \omega_i^2 \mathbf{E}(P(i,i))
\quad \mbox
{and}\quad  P(i,i)=1-\frac{1}{1+q_i d_i^2} .
\]
Further,
\[
\hat{m}\trsp\SigmaHat^{-1}\hat{m}-\frac{\kappa(n,p)}{1-\kappa
(n,p)}\tendsto0 \qquad \mbox{in probability} .
\]

Furthermore, under the above assumptions, if the spectral distribution
of $\Lambda'\Lambda$ converges to $G$ and $(\sum_{i=1}^n \omega_i^2
d_i^2)/n$ remains bounded, a result similar to Theorem~\ref
{thm:SummaryQuadFormsMuHatAndSigmaHat} holds, with $\limitScaling$
being computed by solving equation \eqref{eq:defLimitScaling} with the
corresponding $G$ and $\kappa(n,p)$ playing the role of $\rho
_n/(1-\rho_n)$.
\end{propositionNumThm}

Essentially the previous proposition tells us that when dealing with
correlated variables, the new $\kappa(n,p)$ replaces the old $\kappa
=\rho/(1-\rho)$. We note that
there are no inconsistencies with our previous results as $\sum_i
P(i,i)=\operatorname{trace}(P)=p$ and in the ``elliptical'' case (i.e.,
$\Lambda$
diagonal), $\omega_i^2=1$, so the previous proposition is consistent
with the results we have obtained in the elliptical case. We also
remark that $\| \omega\|=\sqrt{n}$, since $A$ is orthogonal.

Finally, in the case where the $d_i$'s have a limiting spectral
distribution and satisfy \eqref
{eq:AssumptionSmallestSingularValueBoundedBelow}, further
computations show that $q_i-\rho_n\limitScaling\tendsto0$. However,
this does not help (in general) in getting a simpler expression for
$\kappa(n,p)$.

\subsection{On the bootstrap}

An interesting aspect of the analysis of elliptical models is that it
also shed lights on the properties of the bootstrap in this context. As
a matter of fact, the nonparametric bootstrap yields covariance
matrices that have a structure similar to those computed from
elliptical distributions: if we call $D$ the diagonal matrix whose
$i$th diagonal entry is the number of times observation $X_i$ appears
in our bootstrap sample, we have, if $\SigmaHat^*$ is the bootstrapped
covariance matrix,
\[
\SigmaHat^*=\frac{1}{n-1}X\trsp D X-\frac{n}{n-1}\muHat^*(\muHat
^*)\trsp,
\]
where $X$ is our original data matrix, and $\muHat^*$ is the sample
mean of our bootstrap sample, which can also be written $\muHat
^*=X\trsp D \ebold/n$. Unless otherwise noted, we assume in the
discussion that follows that the population mean $\mu$ is 0. Since the
covariance matrix is shift-invariant, we can make this assumption
without loss of generality. We call
\[
\mathfrak{S}^*=\frac{1}{n}X\trsp D X   \quad \mbox{and}\quad  {\cal
S}^*=\Sigma^{-1/2}\mathfrak{S}^*\Sigma^{-1/2} .
\]
As we will see shortly, understanding the properties of $\SigmaHat^*$
boils down to understanding those of ${\cal S}^*$ so we will focus on
this slightly more convenient object in this short discussion.

We note that if $X$ is Gaussian, $\mathfrak{S}^*$ can be thought of as
a ``covariance matrix'' computed from the elliptical data $\XTilde
_i=d_i^{1/2}X_i$. The same remark applies when $X$ is elliptical, that is,
for us, $X_i=\lambda_i {\cal N}(0,\Sigma)$: all we need to do is
change the ``ellipticity parameter'' $\lambda_i$ to $\sqrt
{d_i}\lambda
_i$. The same remark is also applicable to the case of correlated
observations, that is, $X=\Lambda Y \Sigma^{1/2}$, where $\Lambda$ is
not diagonal anymore. Studying the bootstrap properties of such a model
is the same as studying that of the model where we replace $\Lambda$
by $\sqrt{D}\Lambda$. We therefore would like to apply directly all
the results we have obtained above in our study of elliptical models to
better understand the bootstrap. For quantities of the form $v\trsp
(\SigmaHat^*)^{-1} v$, we will see that we can essentially do it, but
differences will appear when dealing with $(\muHat^*)\trsp(\SigmaHat
^*)^{-1} \muHat^*$, which yields statistics that are not exactly
analogous to corresponding statistics appearing in the elliptical case.

Our focus will be on bias properties of bootstrapped replications, so
we will aim for convergence in probability results and not fluctuation
behavior. Our overall strategy here is to show convergence in
probability of the quantities we are interested in as functions of both
the $d_i$'s and $X_i$'s. We will derive the convergence properties of our
bootstrapped statistics by then conditioning on the data and arguing
that with high probability (over the $X_i$'s), this does not change the
results much. We first give some needed background on the bootstrap in
Sections~\ref{subsubsec:NeededConvPropertiesBoot}
and~\ref{subsubsec:PropertiesBootWeights}, then turn to properties of
quantities like $v\trsp(\SigmaHat^*)^{-1}v$ (in Section \ref
{subsubsec:InvCovMatBootData}) and finally study $(\muHat^*)
\trsp(\SigmaHat^*)^{-1}\muHat^*$ (in Section \ref
{subsubsec:BootResQuadFormsMuHatSigmaHat}), where we will see (in
Proposition~\ref{proposition:BootstrappingGaussianData}) some key
differences with the elliptical case.
We conclude this subsection with a brief discussion of the parametric
bootstrap and the conclusions that can be reached about it through our results.

\subsubsection{A remark on needed convergence properties}\label
{subsubsec:NeededConvPropertiesBoot}
Making statements about bootstrapped statistics requires us to make
statements that are conditional on the observed data. This is not a
trivial matter for the statistics we deal with since they cannot be
easily described in terms of simple formulas involving the original
observations. However, we can take a roundabout way: by showing joint
convergence in probability (joint here refers to the ``new'' data being
the vectors of bootstrapped weights and observations), we can obtain
interesting conclusions conditional on the data. Though this is not
difficult to show, we give full arguments here for the sake of completeness.

We will look at our statistics as functions of the number of times an
observation appears in the sample and also, of course, of our
observations. In other words, the original statistic, $T_n$ can be written
\[
T_n=f(1,\ldots,1,X_1,\ldots,X_n)
\]
and, the bootstrapped version $T_n^*$ is, if observation $X_i$ appears
$w_i^*$ times in the bootstrap sample,
\[
T_n^*=f(w^*_1,\ldots,w^*_n,X_1,\ldots,X_n) .
\]
The following simple proposition is used repeatedly in our bootstrap work.
\begin{propositionNumThm}\label
{prop:JoinConvInProbaImpliesBootstrapConvergence}
Let us consider a statistic $T_n=f(w_1,\ldots,w_n,X_1,\ldots,\break X_n)$,
where $w_i$ is the number of times $X_i$ appears in our sample. Suppose
that the vector of weights, $w$, is independent of the data matrix $X$.
Denote by ${\cal Q}_n$ the joint probability distribution of the
$w_i$'s, ${\cal P}_n$ the joint probability distribution of the $X_i$'s
and ${\cal R}_n={\cal Q}_n\times{\cal P}_n$ the probability
distribution of $(w_1,\ldots,w_n,X_1,\ldots,X_n)$.

Suppose we have established that $T_n$ tends in ${\cal
R}_n$-probability to $c$, a deterministic object, as $n\tendsto\infty$.

Then we have, with ${\cal P}_n$-probability going to 1 as $n\tendsto
\infty$,
\[
T_n|\{X_i\}_{i=1}^n\tendsto c \qquad\mbox{in } {\cal Q}_n\mbox
{-probability} .
\]
In other words, calling ${\cal X}_n=\{X_i\}_{i=1}^n$, for all $\eps,
\eta>0$, if $Q_n(\eps)={\cal Q}_n(|T_n-c|>\eps|{\cal X}_n)$, ${\cal
P}_n(Q_n(\eps)>\eta)\tendsto0$ as $n$ tends to infinity.
\end{propositionNumThm}

In the case where the weights $w_i$ are obtained by standard
bootstrapping, $Q_n$ is multinomial($1/n,\ldots,1/n,n$). Then,
$T_n|X_n$ has the distribution of the usual
bootstrap quantity $T_n^*$. We will focus on this case more
specifically later.
\begin{pf*}{Proof of Proposition~\ref{prop:JoinConvInProbaImpliesBootstrapConvergence}}
The proof and the statement are almost obvious but we include them for
the sake of completeness. Let us call $\tau_n=|T_n-c|$ and ${\cal
X}_n=\{X_1,\ldots,X_n\}$. By assumption, $\tau_n\tendsto0$ in ${\cal
R}_n$ probability. Hence,
\[
\myExp_{{\cal R}_n} (1_{\tau_n>\eps} )=\myExp_{{\cal
P}_n} (\myExp_{{\cal Q}_n} [1_{\tau_n>\eps}|{\cal
X}_n ] )\tendsto0 .
\]
Let us call $Q_n(\eps)={\cal Q}_n(|T_n-c|>\eps|{\cal X}_n)$. Clearly,
$0\leq Q_n(\eps)\leq1$ and\break  $\myExp_{{\cal P}_n} (Q_n(\eps
) ) \tendsto0$, so for any $\eta>0$,
\[
{\cal P}_n\bigl(Q_n(\eps)>\eta\bigr) \tendsto0 .
\]
\upqed
\end{pf*}

We now investigate the case of the classical bootstrap, that is, the
situation in which ${\cal Q}_n$ is multinomial$(\frac{1}{n},\ldots
,\frac{1}{n},n)$.

\subsubsection{Empirical distribution of bootstrap weights}\label
{subsubsec:PropertiesBootWeights}
As we saw in Theorem~\ref{thm:QuadFormsInverseEllipticalCase}, the
empirical distribution of the ellipticity parameters affect crucially
statistics of the type $v\trsp\SigmaHat^{-1}v$, so to understand the
effect of bootstrapping, we need to understand the empirical
distribution of the bootstrap weights. This question has surely been
investigated, but we did not find a good reference, so we provide the
result and a simple proof for the convenience of the reader.\vadjust{\goodbreak}

\begin{propositionNumThm}\label{prop:EmpDistBootWeights}
Let the vector $w$ be distributed according to a multi\-nomial$(\frac
{1}{n},\ldots,\frac{1}{n},n)$ distribution. Call $F_n$ the empirical
distribution of the vector $w$. Then
\[
F_n \WeakCv\Poisson(1) \qquad \mbox{in probability} ,
\]
where $\Poisson(1)$ is the Poisson distribution with parameter 1.

\end{propositionNumThm}

\begin{pf}
Let us first start by an elementary remark: suppose $\pi_1,\ldots,\pi
_n$ are i.i.d. with distribution $\Poisson(1)$. Call $\Pi_n=\sum
_{i=1}^n \pi_i$. Then
\[
(\pi_1,\ldots,\pi_n )| \{\Pi_n=n \} \sim
\operatorname{multinomial}\biggl(\frac{1}{n},\ldots,\frac{1}{n},n\biggr) .
\]
This result is a simple application of Bayes's rule and the fact that
$\Pi_n\sim\Poisson(n)$.

Let us now show that if $f$ is bounded and continuous, and if $W\sim
\Poisson(1)$,
\[
\myExp_{F_n}(f)=\frac{1}{n}\sum_{i=1}^n f(w_i) \tendsto\mathbf{E}(f(W))
\qquad \mbox{in probability} .
\]
To do so, we note that $w_i\sim\operatorname{binomial}(n,1/n)$ and therefore
its marginal distribution is asymptotically $\Poisson(1)$. Therefore,
\[
\mathbf{E}(\myExp_{F_n}(f))\tendsto\mathbf{E}(f(W)) .
\]
Now all we need to do is therefore to show that $\operatorname
{var}(\myExp_{F_n}(f))$ goes to zero. Clearly, by independence of the
$\pi_i$'s,
\[
\operatorname{var}\Biggl(\frac{1}{n}\sum_{i=1}^n f(\pi_i)\Biggr)=\frac
{1}{n}\operatorname{var}(f(W))=\gO
\biggl(\frac{1}{n} \biggr) ,
\]
because $f$ is bounded. But our first remark implies that
\[
\operatorname{var}(\myExp_{F_n}(f))=\operatorname{var}\Biggl( \frac
{1}{n}\sum_{i=1}^n f(\pi_i)\Big |\Pi_n=n\Biggr).
\]
Now,
\begin{eqnarray*}
\operatorname{var}\Biggl(\frac{1}{n}\sum_{i=1}^n f(\pi_i)\Biggr)&=&\mathbf
{E}\Biggl(\operatorname{var}\Biggl( \frac{1}{n}\sum_{i=1}^n f(\pi_i)\Big |\Pi
_n\Biggr)\Biggr)+\operatorname{var}\Biggl( \mathbf{E}\Biggl(\frac{1}{n}\sum_{i=1}^n f(\pi
_i)\Big |\Pi_n\Biggr)\Biggr)\\
&\geq&\operatorname{var}\Biggl( \frac{1}{n}\sum_{i=1}^n f(\pi
_i)\Big |\Pi_n=n\Biggr)
P(\Pi_n=n) .
\end{eqnarray*}
Since $\Pi_n$ has $\Poisson(n)$ distribution, $P(\Pi_n=n)\sim
1/\sqrt{2\pi n}$. Hence,
\[
\operatorname{var}(\myExp_{F_n}(f))=\operatorname{var}\Biggl(\frac
{1}{n}\sum_{i=1}^n f(\pi_i)|\Pi_n=n\Biggr) =\gO(n^{-1/2})\tendsto0 ,
\]
and the result is established.\vadjust{\goodbreak}
\end{pf}

We will also need later to use on the following (coarse) fact:
\begin{factNumThm}\label{fact:BoundMaximumBootstrapSample}
Let the vector $w$ be distributed according to a multino\-mial$(\frac
{1}{n},\ldots, \frac{1}{n},n)$ distribution. Then
\[
P\Bigl(\max_{i=1,\ldots,n} w_i>(\log n)\Bigr)=\gO\biggl(\frac{n^{3/2}}{(\log
n)!} \biggr) .
\]
In particular, this probability goes to 0 faster than any $n^{-a}$, $a>0$.
\end{factNumThm}

The proof of the fact is elementary, and relies on the representation
used above for the vector $w$, a simple union bound, the fact that
$P(\Poisson(n)=n)\sim n^{-1/2}$ and the fact that
$P(\Poisson(1)\geq M)\leq(M!)^{-1} M/(M-1)$ which is easy to see by
writing explicitly the probability we are trying to compute.

With these preliminaries behind us, we are now ready to tackle the
question of understanding the (first-order) bootstrap properties of the
statistics appearing in the study of quadratic programs with linear
equality constraints.

\subsubsection{On inverse covariance matrices computed from
bootstrapped data}\label{subsubsec:InvCovMatBootData}
Our aim in this subsubsection and the next is to find analogs to
Theorems~\ref{thm:QuadFormsInverseEllipticalCase} and Theorems \ref
{thm:SummaryQuadFormsMuHatAndSigmaHat}. Our first result along these
lines is an analog of Theorem~\ref{thm:QuadFormsInverseEllipticalCase}.

We present the result in the case of Gaussian data, where we can get a
somewhat explicit expression for the quantity we care about, and
discuss possible extensions below.
\begin{theorem}\label{thm:bootQuadFormsInverseEllipticalCaseNormalCase}
Suppose we observe $n$ i.i.d. observations $X_i$, where $X_i$ are
i.i.d. in $\mathbb{R}^p$ with distribution ${\cal N}(\mu,\Sigma_p)$.
Call $\rho_n=p/n$ and assume that $\rho_n\tendsto\rho\in(0,1-\me
^{-1})$. Call $\SigmaHat^*$ the covariance matrix computed after
bootstrapping the $X_i$'s. Call ${\cal P}_n$ the joint distribution of
the $X_i$'s.

If $v$ is a (sequence of) deterministic vectors, then conditional on $\{
X_i\}_{i=1}^n$, with high ${\cal P}_n$ probability,
\[
\frac{v\trsp(\SigmaHat^*)^{-1} v}{v\trsp\Sigma^{-1} v}\tendsto
\limitScaling\qquad \mbox{in probability} ,
\]
where $\limitScaling$ satisfies, if $G$ is a $\Poisson(1)$ distribution
\begin{equation}\label{eq:defLimitScalingBootstrap}
\int\frac{dG(\tau)}{1+\rho\tau\limitScaling}=1-\rho.
\end{equation}
\end{theorem}

\begin{pf}
As before, we call ${\cal Q}_n$ the law of the bootstrap weights [i.e.,
multinomial$ (\frac{1}{n},\ldots,\frac{1}{n},n )$] and
${\cal R}_n={\cal Q}_n \times{\cal P}_n$.
Without loss of generality, we can assume that $\mu=0$. Let us call
$D$ the diagonal matrix\vadjust{\goodbreak} containing the bootstrap weights. We have
$\muHat^*=X\trsp D \ebold/n$. Also, it is true that
\[
\SigmaHat^*=\frac{1}{n-1} \biggl(X-\frac{\ebold(\muHat^*)\trsp
}{n}\biggr)\trsp D\biggl (X-\frac{\ebold(\muHat^*)\trsp}{n}\biggr) .
\]
Since $\ebold\trsp D \ebold=n$, we also have
\[
(n-1)\SigmaHat^*=X\trsp D \biggl(\id-\frac{\ebold\ebold\trsp
D}{n}\biggr)X=X\trsp D^{1/2}\biggl (\id-\frac{1}{n} D^{1/2}\ebold\ebold
\trsp D^{1/2} \biggr)D^{1/2}X .
\]
Because $X$ is of the form $X=Y\Sigma^{1/2}$ under our assumptions, we
see that
\begin{eqnarray}
\hspace*{-80pt}\SigmaHat^*=\Sigma^{1/2}{\cal S}^*\Sigma^{1/2}  \nonumber\\
\eqntext{\displaystyle\mbox{where }{\cal S}^*=\frac{1}{n-1}Y\trsp D^{1/2} \biggl(\id-\frac{1}{n}
D^{1/2}\ebold\ebold\trsp D^{1/2} \biggr)D^{1/2} Y .}
\end{eqnarray}
If we call $\delta=D^{1/2}\ebold$, we have $\| \delta\|_2^2=n$
because the sum of the bootstrap weights is $n$. Therefore, $H_{\delta
}=\id_n-\delta\delta\trsp/n\succeq0$. Also, $H_{\delta}$ (like
$H$) is a projection matrix and a rank 1 perturbation of $\id_n$.

The situation is therefore very similar to the question we studied in
Theorem~\ref{thm:QuadFormsInverseEllipticalCase}, except that $H=\id
_n-\ebold\ebold\trsp/n$ is replaced by $H_{\delta}=\id_n-\delta
\delta\trsp/n$. All the arguments given there hold provided we can
show that \eqref{eq:AssumptionSmallestSingularValueBoundedBelow} is
satisfied for the bootstrap weights in
the situation we have here.

Now let us call $N$ the number of nonzero bootstrap weights. In the
notation of Theorem~\ref{thm:QuadFormsInverseEllipticalCase}, $\lambda
_i=\sqrt{d_i}$ and $\tau_i=d_i$. So clearly, $\tau_{(N)}\geq1$. So
$C_0=1/2$ is a possibility. Also, $N/n\tendsto1-1/e$ in probability,
so $p/N$ has a limit in probability and this limit is bounded away from
1 because of our assumption that $\rho_n\tendsto\rho\in(0,1-1/e)$.
Finally, we can pick $\eta_0=(1-1/e)/2$.

So the proof of Theorem~\ref{thm:QuadFormsInverseEllipticalCase}
applies [it is easy to see here that the assumption that $\tau_i\neq
0$ can be dispensed of, because we know that the nonzero $\tau_i$'s
are large enough for our arguments to go through, and there are enough
of them that we do not have problems (at least in probability) with
$\SigmaHat^{-1}$ not being defined], and we have the announced result.
\end{pf}

The previous theorems settled the question of understanding the impact
of the nonparametric bootstrap on statistics of the form $v\trsp
\SigmaHat^{-1}v$ in the situation where the original data were
Gaussian. A similar analysis could be carried out in the case of
elliptical data, when we assume that the ``ellipticity'' parameters,
$\lambda_i$, are such (\ref
{eq:AssumptionSmallestSingularValueBoundedBelow}) is satisfied for the
``new weights'' $\tau_i=\lambda_i^2 w_i$. The result would then depend
on the limiting distribution of $\lambda_i^2 w_i$ (if it exists),
where $w_i$ is the bootstrap weight given to observation $i$.

\subsubsection{\texorpdfstring{Bootstrap analogs of Theorems \protect\ref{thm:QuadFormsMuHatSigmaHatInvV} and \protect\ref{thm:SummaryQuadFormsMuHatAndSigmaHat}}
{Bootstrap analogs of Theorems 4.5 and 4.6}}\label{subsubsec:BootResQuadFormsMuHatSigmaHat}

An important piece of our analysis of quadratic programs with linear
equality constraints when the data are elliptically distributed was the
study of quadratic forms of the type $\muHat\trsp\SigmaHat
^{-1}\muHat$. It is natural to ask what happens to them when we
bootstrap the data.\vadjust{\goodbreak} In the elliptical case, we saw that the key
statistic was of the form,
when $\mu=0$ and ${\mathfrak S}=\Sigma^{1/2}Y\trsp\Lambda^2 Y\Sigma
^{1/2}/n$,
\[
\muHat\trsp{\mathfrak S}^{-1} \muHat=\frac{1}{n}\ebold\trsp
\Lambda Y(Y\trsp\Lambda^2 Y)^{-1} Y\trsp\Lambda\ebold.
\]
However, in the bootstrap case, if $\Lambda$ is the diagonal matrix
containing the bootstrap weights, we have ${\mathfrak S}^*=\Sigma
^{1/2}Y\trsp\Lambda Y \Sigma^{1/2}/n$, but $\muHat^*=Y\trsp\Lambda
\ebold/n$, so the key statistic is going to be of the form
\[
(\muHat^*)\trsp({\mathfrak S}^*)^{-1}(\muHat^*)=\frac{1}{n}\ebold
\trsp\Lambda Y(Y\trsp\Lambda Y)^{-1} Y\trsp\Lambda\ebold.
\]
This creates complications because the matrix $\Lambda Y(Y\trsp\Lambda
Y)^{-1} Y\trsp\Lambda$ is not a projection matrix, and hence some of
our previous analysis cannot be applied directly. However, this
statistic can be rewritten, if we denote $w=\Lambda^{1/2} \ebold$, as
\[
\frac{1}{n}w\trsp\Lambda^{1/2} Y(Y\trsp\Lambda Y)^{-1} Y\trsp
\Lambda^{1/2} w=\frac{1}{n}w\trsp P_{\Lambda^{1/2}} w ,
\]
where $P_{\Lambda^{1/2}}$ is now a projection matrix. As before its
off-diagonal elements have mean 0 (conditional on $\Lambda$), but now
we also need to understand\break $\sum_{i=1}^n w_i P_{i,i}/n$ and not only
$\sum_{i=1}^n P_{i,i}/n$. A detailed analysis of the former quantity
is done in Appendix~\ref{App:Subsec:BootSpecResults}.

We naturally now assumes that $p/n$ has a finite limit, $\rho$ in
$(0,1-1/e)$. As explained in Appendix \ref
{App:Subsec:BootSpecResults}, $\sum_{i=1}^n w_i P_{i,i}/n\tendsto
(\limitScaling-1)/\limitScaling$ in ${\cal Q}_n$-probability, with
${\cal P}_n$ probability going to 1, where $\limitScaling$ is computed
by solving equation~\eqref{eq:defLimitScalingBootstrap} [i.e., using
$\Poisson(1)$ for $G$ in that equation].

Similarly, it is explained there, that with ${\cal P}_n$ probability
going to 1, when $X_i$ have mean 0,
\[
(\muHat^*)\trsp(\SigmaHat^*)^{-1}\muHat^* \tendsto\limitScaling
-1\geq\frac{\rho}{1-\rho}  \qquad \mbox{in } {\cal Q}_n \mbox{-probability}.
\]

Finally, an analog of Theorem~\ref{thm:QuadFormsMuHatSigmaHatInvV}
holds, so we have an analog of Theorem~\ref
{thm:SummaryQuadFormsMuHatAndSigmaHat}, where $\limitScaling$ is as
defined above, and $\rho_n/1-\rho_n$ needs to be replaced by
$\limitScaling-1$.

In summary, we have the following proposition.
\begin{propositionNumThm}\label{proposition:BootstrappingGaussianData}
Call $\limitScaling$ the quantity defined by equation \eqref
{eq:defLimitScalingBootstrap}.

Suppose the data $X_1,\ldots,X_n$ is i.i.d. ${\cal N}(\mu,\Sigma)$,
and call ${\cal P}_n$ the corresponding probability distribution.
Suppose $v$ is a given deterministic sequence of vectors. Under the
assumptions of Theorem \ref
{thm:bootQuadFormsInverseEllipticalCaseNormalCase}, we have, when
bootstrapping the data, with ${\cal P}_n$ probability going to 1
\begin{eqnarray*}
\frac{v\trsp(\SigmaHat^*)^{-1}v}{v\trsp\Sigma^{-1} v}&\tendsto&
\limitScaling\qquad\mbox{in } {\cal Q}_n \mbox{-probability,}\\
\frac{(\muHat^*)\trsp(\SigmaHat^*)^{-1}v}{\sqrt{v\trsp\Sigma
^{-1} v}}&\tendsto&0 \qquad\mbox{in } {\cal Q}_n\mbox{-probability, when }
\mu=0 ,\\
(\muHat^*)\trsp(\SigmaHat^*)^{-1}\muHat^*&\simeq&\limitScaling\mu
\trsp\Sigma^{-1}\mu+(\limitScaling-1)+\lo_{{\cal Q}_n} \bigl(\sqrt
{\mu\trsp\Sigma^{-1}\mu},1 \bigr).
\end{eqnarray*}
\end{propositionNumThm}

We note that our techniques could yield generalizations of the previous
fact for the case where the data is elliptically distributed. However,
in the case where $X_i$ have mean 0, the quantity $(\muHat^*)\trsp
(\SigmaHat^*)^{-1}\muHat^*$ does not seem to have a limiting value
that is writable in compact form, so we do not dwell on this question further.

Naturally, the motivation behind the previous proposition is practical
and the results are interesting from that standpoint. They show that
the bootstrap yields inconsistent estimators of the population quantities,
something that is not completely unexpected when we understand the
random matrix aspects of these questions. Perhaps even more interesting
is that bootstrap estimates of bias are themselves inconsistent: as a
matter of
fact, the key quantity that measures bias in the Gaussian case is
$1/(1-p/n)$; when bootstrapping it is replaced by $\limitScaling$, as
defined in equation \eqref{eq:defLimitScalingBootstrap}.
These results therefore cast some doubts on the practical relevance of
the bootstrap for the high-dimensional problems we are considering, at
least when the bootstrap is used in ``classical'' ways.

\subsubsection{On the parametric bootstrap}
In the settings considered here, it is also natural to ask how the
parametric bootstrap would behave. For instance, if we assumed
Gaussianity of the data, we could just estimate $\Sigma$ and $\mu$
(e.g., naively, by $\SigmaHat$ and $\muHat$) and use a parametric
bootstrap to get at the quantities we are interested in.

Naturally, the analysis of such a scheme is similar to the analysis of
the Gaussian case carried out in Section~\ref{sec:GaussianCase}, where
the population parameters $\Sigma$ and $\mu$ need to be replaced by
the estimators we use in our parametric bootstrap. The same would be
true if we were to do a parametric bootstrap for elliptical data, but
we would have to use the results of Section~\ref{sec:EllipticalCase} instead.

Our computations show that the parametric bootstrap could be used in
the problems under study to estimate the bias of various plug-in
estimators: we would for instance recover the correct $\limitScaling$
by considering $v\trsp(\Sigma^*_{\mathrm{parametric}})^{-1} v/\break v\trsp
\SigmaHat^{-1} v$. We note, however, that
our analyses, and the estimation work we carry out in Section \ref
{sec:ComparisonGaussianElliptical} could do this too, at a cheaper
numerical cost.

Finally and very interestingly, we see that a naive use of the
parametric bootstrap to estimate the bias in the empirical efficient
frontier---a perhaps reasonable idea at first glance---would yield
inconsistent estimates of bias.

\section{Robustness, bias  and improved estimation}\label
{sec:ComparisonGaussianElliptical}
We now go back to our original problem, which was to understand the
relationship between the solution of problem \eqref{eq:QP-eqc-Emp} and
the solution of problem \eqref{eq:QP-eqc-Pop} (see page \pageref
{eq:QP-eqc-Emp} for definitions).

It is naturally important to understand the effect of making the
assumption that the data is normally distributed as compared to, say,
an assumption of elliptical distribution for the data. The following
discussion fleshes out some of our theoretical results and what their
significance is when solving\vadjust{\goodbreak} quadratic programs with linear
constraints. The discussion is an application of the work done in
Sections~\ref{Sec:GeneralQP}--\ref{sec:EllipticalCase}. It might
appear to be mainly heuristic, but precise statements can be easily
deduced from the precise statements of the theorems
given in the corresponding technical sections.

We discuss here only the case of i.i.d. data. As we have shown above,
the bootstrap case and the case of correlated observations are more
complicated to handle, and the formulas are not as explicit in those
cases as they are in the case of i.i.d. data. But for certain cases,
one could plug-in our earlier results for those situations to obtain
explicit results about efficient frontiers and weight vectors in those
cases too.

As a matter of notation, all of our approximation statements hold with
high-probability asymptotically, unless otherwise noted. We will carry
out our work under the model put forward in Theorem \ref
{thm:QuadFormsInverseEllipticalCase}, assuming that the $\lambda_i$'s
are i.i.d. and the following assumptions:\vspace*{-2pt}
\begin{enumerate}
\item[Assumption A1:]  for all $i \in\{1,\ldots,k\}$, $v_i\trsp\Sigma
^{-1} v_i$ stays bounded away from 0. $v_k$ is assumed to be equal to
$\mu$.
\item[Assumption A2:]  the smallest eigenvalue of $M=V\trsp\Sigma
^{-1}V$ stays bounded away from 0 and the condition number of $M$
remains bounded.
\item[Assumption A3:] if $\eps=\pm1$, $(v_i+\eps v_j)\trsp\Sigma
^{-1} (v_i+\eps v_j)$ stays bounded away from infinity.
\item[Assumption A4:]  \eqref
{eq:AssumptionSmallestSingularValueBoundedBelow} and \eqref
{eq:AssumptionSecondAndFourthMomentLambda} hold. (See Theorem~\ref
{thm:SummaryQuadFormsMuHatAndSigmaHat} for definitions.)
\item[Assumption A5:]  we have, for some\vspace*{1pt} $\eps>0$, if $u_n= (2\log
(n)+(\log n)^{\eps} )^{1/2}+\sqrt{2\pi}$, $\frac{\opnorm
{\Sigma}}{\operatorname{trace}(\Sigma)}u_n^2\tendsto0$, where $\opnorm
{\Sigma}$
is the largest eigenvalue of $\Sigma$.\vspace*{-2pt}
\end{enumerate}
These assumptions guarantee that the noise terms involving $\muHat$ do
not overwhelm the signal terms involving $\mu$, and also that we can
safely take inverses of our approximations to get approximations of
their inverses. Also, all the key results we obtained in Sections \ref
{sec:GaussianCase} and~\ref{sec:EllipticalCase} are applicable, and
our conclusions will of course heavily rely on them.

We will use the notation $\rho_n=p/n$. We recall that in the Gaussian
case, the quantity $\limitScaling$ appearing below is approximately
equal to $1/(1-\rho_n)$ and in the elliptical case, it is always
greater than $1/(1-\rho_n)$, as we explained after the proof of
Theorem~\ref{thm:QuadFormsInverseEllipticalCase}.

We start by investigating the case of equality constraints. We discuss
inequality constraints in
Section~\ref{Subsec:InequalityConstraints}.\vspace*{-4pt}

\subsection{Relative positions of efficient frontiers: Gaussian vs.
elliptical case}
When assumptions (A1--A4) hold, it is clear that
\begin{equation}\label{eq:ApproxMHat}
\MHat=\VHat\trsp\SigmaHat^{-1}\VHat\simeq\limitScaling V\trsp
\Sigma^{-1} V + \frac{\rho_n}{1-\rho_n} e_k e_k\trsp.\vspace*{-2pt}
\end{equation}
Now recall that in the elliptical case, $\limitScaling\geq
1/(1-p/n)=\limitScaling^{G}$, that is, the ``$\limitScaling$''
corresponding to the Gaussian case. Calling $\widehat{M}_E$ the
empirical estimator of $M$ we get in the\vadjust{\goodbreak} elliptical case and $\widehat
{M}_G$ its analog in the Gaussian case, we have, when A1--A4 are
satisfied, with high-probability,
\[
\widehat{M}_E\succeq\widehat{M}_G ,
\]
at least asymptotically.

We now call $\fEmp^{(E)}$ and $\fEmp^{(G)}$ the ``efficient
frontiers'' obtained by solving problem \eqref{eq:QP-eqc-Emp} when the
data is respectively elliptical and Gaussian. Recall that under our
assumptions, $\mu$ and $\Sigma$ are the same for the two problems, so
the population version corresponding to the two problems is the same.
We call the population solution, that is, the efficient frontier
computed with the population parameters, $\fTheo$. Naturally, this is
the quantity we are fundamentally interested in estimating.

Using the fact that $\fEmp=U\trsp\widehat{M}^{-1} U$, the following
important results.
\begin{theorem}\label{thm:PositionFrontierEllipticalAndGaussianCase}
When assumptions \textup{A1--A4} are satisfied, we have with high-probability and
asymptotically,
\[
\fEmp^{(E)}\leq\fEmp^{(G)}\leq\fTheo.
\]
In other words, risk underestimation in the empirical quadratic program
with linear equality constraints is least severe (within the class of
elliptical models) in the Gaussian case.

In other respects, we have, asymptotically, with high-probability, if
$\kappa_n=\rho_n/(1-\rho_n)$,
\begin{equation}\label{eq:PositionEmpFrontier}
\fEmp^{(E)}\simeq\frac{1}{\limitScaling} \biggl(\fTheo-\frac
{\kappa_n}{\limitScaling} \frac{ (e_k\trsp M^{-1} U
)^2}{1+ ({\kappa_n/\limitScaling}) e_k\trsp M^{-1} e_k} \biggr) .
\end{equation}

\end{theorem}

Another way of phrasing this result is the fact that the Gaussian
analysis gives the most optimistic view of risk underestimation within
the class of elliptical models considered here.

Practically, it means that users of Markowitz-type optimization should
be wary of the empirical solution they get, and even of the correction
that Gaussian results suggest. If the data is elliptical, they will
underestimate the risk of their portfolio even more than the Gaussian
results suggest.

Let us now give a proof of Theorem \ref
{thm:PositionFrontierEllipticalAndGaussianCase}.
\begin{pf*}{Proof of Theorem~\ref{thm:PositionFrontierEllipticalAndGaussianCase}}
Under the assumptions of the theorem, we can use the approximation in
equation \eqref{eq:ApproxMHat}. The first part of the theorem has been
argued before, so we do not need to do anything else to obtain it.

The second part follows directly from a rank one perturbation argument.
We have
\begin{eqnarray*}
\fEmp^{(E)}&\simeq& U\trsp\biggl(\limitScaling V\trsp\Sigma^{-1} V +
\frac{\rho_n}{1-\rho_n} e_k e_k\trsp\biggr)^{-1} U\\
&=&\frac
{1}{\limitScaling} U\trsp\biggl(M+\frac{\kappa_n}{\limitScaling
}e_ke_k\trsp\biggr)^{-1}U .
\end{eqnarray*}
Using the classic result $(M+\nu\nu\trsp)^{-1}=M^{-1}-M^{-1}\nu\nu
\trsp M^{-1}/(1+\nu\trsp M^{-1} \nu)$, we conclude that
\[
U\trsp\biggl(M+\frac{\kappa_n}{\limitScaling}e_ke_k\trsp
\biggr)^{-1}U=U\trsp M^{-1} U -\frac{\kappa_n}{\limitScaling}\frac
{(U\trsp M^{-1}e_k)^2}{1+ ({\kappa_n/\limitScaling})e_k\trsp
M^{-1} e_k} .
\]
We now recall from Section~\ref{Sec:GeneralQP} that $\fTheo=U\trsp
M^{-1} U$, and we have the announced result.
\end{pf*}

Equation \eqref{eq:PositionEmpFrontier} naturally suggests better ways
of estimating $\fTheo$ than using $\fEmp$. We postpone a discussion
of this issue to Section~\ref{subsec:EstimationEffFrontier} because
it requires somewhat lengthy preliminaries.

\subsection{Issues concerning the weights of the portfolio}\label
{Sec:GaussEllip:Subsec:WeightsPortfolio}
Besides problems in the location of the efficient frontiers, our
analysis reveals another very interesting phenomenon: problems with
estimating $\wTheo$, the optimal vector of weights. In particular, one
can show that the mean return of the portfolio is poorly estimated and
the weight given to each asset is biased.

\begin{theorem}[(Bias in weights)]\label{thm:BiasInWeights}
Suppose assumptions \textup{A1--A4} hold. We have, asymptotically and with
high-probability,
\begin{equation}\label{eq:BiasInWeights}
\wEmp\simeq\wTheo-\zeta(\limitScaling)\frac{\kappa
_n}{\limitScaling} w_b ,
\end{equation}
where
\[
\zeta(\limitScaling)=\frac{e_k\trsp M^{-1} U}{1+ ({\kappa
_n/\limitScaling})e_k\trsp M^{-1} e_k} ,\qquad w_b=\Sigma^{-1} VM^{-1}
e_k .
\]
This approximation is valid when looking at linear combinations of the
vector of weights: if $\gamma\in\mathbb{R}^n$ is deterministic and
assumption \textup{A3} extended to include this vector holds,
\[
\gamma\trsp\wEmp= \gamma\trsp\biggl(\wTheo-\zeta(\limitScaling
)\frac{\kappa_n}{\limitScaling} w_b \biggr)+\lo_P(1) .
\]
\end{theorem}

We note that the last assertion of the theorem does not necessarily
immediately follow from equation \eqref{eq:BiasInWeights} in
high-dimension, but it is true in the setting we consider. A
particularly interesting corollary is the following statement
concerning inconsistent estimation of the returns.
\begin{corollaryNumThm}[(Poor estimation of returns)]
Recall that with our notations, $\wTheo\trsp\mu=u_k=\mu_P$. In
practical terms, $\mu_P$ corresponds to the desired expected returns
we wish to have for our ``portfolio.''
Under the same assumptions as that of Theorem~\ref{thm:BiasInWeights},
we have
\[
\mu\trsp\wEmp\simeq\mu_P\frac{1}{1+ ({\kappa
_n/\limitScaling})e_k\trsp M^{-1} e_k}-\frac{\kappa_n}{\limitScaling
}\frac{\sum_{i<k} u_i e_k\trsp M^{-1} e_i}{1+ ({\kappa
_n/\limitScaling})e_k\trsp M^{-1} e_k} .
\]
\end{corollaryNumThm}

The previous corollary is a statement about poor estimation of returns
for the following reason: $\muHat\trsp\wEmp=\mu_P$ by construction,
so one might naively hope that, for a new observation $X_{n+1}$,
independent of $X_1,\ldots,X_n$ and with the same distribution as them,
$\mathbf{E}(\wEmp\trsp X_{n+1}|X_1,\ldots,X_n)=\wEmp\trsp\mu
\simeq\mu
_P$. However, as the previous corollary shows, this is not satisfied.
We note that the factor affecting $\mu_P$ is a shrinkage factor,
always smaller than 1 because $M$ is positive semi-definite. The other
term could have either sign, so its effect on return estimation is less
interpretable. For large $\mu_P$, it is nonetheless clear that the
previous corollary shows that the returns are overestimated: the
realized returns are (asymptotically and with high-probability) less
than $\mu_P$. Hence, our result can be seen as a generalization of the
overestimation of returns result first found in \citep
{JobsonKorkie80}, in the low-dimensional Gaussian case.

We now prove these two results. The proof of the corollary is at the
end of the proof of the theorem.

\begin{pf*}{Proof of Theorem~\ref{thm:BiasInWeights}}
Under the assumptions of the theorem we have
\[
\MHat\simeq\limitScaling M+\kappa_n e_ke_k\trsp,
\]
and our assumptions guarantee that we can take inverses and still have
valid approximations. Hence, using the classic formula for inversion of
a rank one perturbation of a matrix [see \cite{hj}, page 19], we have
\[
\MHat^{-1}\simeq\frac{1}{\limitScaling} \biggl(M^{-1}-\frac{\kappa
_n}{\limitScaling}\frac{M^{-1}e_k e_k\trsp M^{-1}}{1+ ({\kappa
_n/\limitScaling})e_k\trsp M^{-1} e_k} \biggr) .
\]
Now recall that $\wEmp=\SigmaHat^{-1} \VHat\MHat^{-1}U$ and $\wTheo
=\Sigma^{-1} V M^{-1}U$. For a deterministic $\gamma$, our work in
Section~\ref{sec:EllipticalCase} indicates that $\gamma\trsp
\SigmaHat^{-1} \VHat= \limitScaling\gamma\trsp\Sigma^{-1} V +\lo_P(1)$.
So we conclude that
\[
\gamma\trsp\wEmp= \limitScaling\gamma\trsp\Sigma^{-1} V \frac
{1}{\limitScaling} \biggl(M^{-1}-\frac{\kappa_n}{\limitScaling
}\frac{M^{-1}e_k e_k\trsp M^{-1}}{1+ ({\kappa_n/\limitScaling
})e_k\trsp M^{-1} e_k} \biggr)U+\lo_P(1) .
\]
In other words, we have
\[
\gamma\trsp\wEmp= \gamma\trsp\Sigma^{-1} V M^{-1} U - \frac
{\kappa_n}{\limitScaling}\frac{\gamma\trsp\Sigma^{-1} V M^{-1}e_k
e_k\trsp M^{-1}U}{1+({\kappa_n/\limitScaling
})e_k\trsp M^{-1}
e_k}+\lo_P(1) ,
\]
or, as announced,
\[
\gamma\trsp\wEmp= \gamma\trsp\wTheo-\frac{\kappa
_n}{\limitScaling}\gamma\trsp w_b \zeta(\limitScaling)+\lo
_P(1) .
\]
It seems difficult to say more, because $w_b$ and $\zeta$ are
population parameters and their properties and values may vary from
problem to problem.

$\bullet$ \textit{Proof of the corollary.}

We now assume that $\gamma=\mu$. We remark that $\mu=Ve_k$, by
construction of~$V$. Therefore,
\[
\mu\trsp w_b=e_k\trsp V\trsp\Sigma^{-1} V M^{-1} e_k=e_k\trsp M
M^{-1} e_k=1 .
\]
Further,
\[
e_k\trsp M^{-1} U=\sum_{i=1}^k u_i e_k\trsp M^{-1} e_i=\sum_{i<k}u_i
e_k\trsp M^{-1} e_i+\mu_P e_k\trsp M^{-1} e_k .
\]
These two remarks and the result of Theorem~\ref{thm:BiasInWeights}
give the conclusion of the corollary.
\end{pf*}

\subsection{Bias correction for the weights}
An important question now that we have identified possible problems
with the empirical weights is to try and correct them. We propose such
a scheme, suggested by our computations.

Our investigations will rely on the following asymptotic result,
discussed in Theorem~\ref{thm:BiasInWeights}: in the notations of this theorem,
\[
\gamma\trsp\wEmp= \gamma\trsp\wTheo-\frac{\kappa
_n}{\limitScaling}\gamma\trsp w_b \zeta(\limitScaling)+\lo
_P(1) .
\]

Our efforts will focus on trying to estimate $w_b/\limitScaling$ and
$\zeta(\limitScaling)$, as $\kappa_n=\rho_n/(1-\rho_n)$ is known
and computable from the data.

Recall that we assumed that $v_k=\mu$ and let us call
\[
\MTilde=\MHat-\kappa_n e_k e_k\trsp.
\]
Under the assumptions underlying the previous computations, we have
\[
\MTilde\simeq\limitScaling M .
\]
In practice, we wish $\MTilde$ to be a positive semi-definite matrix---something that
is guaranteed asymptotically, but might require checking
and potentially corrections in practice.

We propose to use:
\begin{enumerate}[2.]
\item As an estimator of $w_b$,
\[
\widehat{w}_b=\SigmaHat^{-1}\VHat\MTilde^{-1} e_k.
\]
\item As an estimator of $\zeta(\limitScaling)/\limitScaling$,
\[
\widehat{z}=\frac{e_k\trsp\MTilde^{-1}U}{1+\kappa_n e_k\trsp
\MTilde^{-1} e_k} .
\]
\end{enumerate}
For any deterministic $\gamma$ (such that the assumptions\vspace*{1pt} of Theorem
\ref{thm:BiasInWeights} hold), $\gamma\trsp\widehat{w}_b\simeq
\gamma\trsp w$, because $\gamma\trsp\SigmaHat^{-1} \VHat\simeq
\limitScaling\gamma\trsp\Sigma^{-1} V$ and $\MTilde^{-1}U\simeq
M^{-1} U/\limitScaling$. Also,\break $e_k\trsp\MTilde^{-1} U\simeq
\limitScaling^{-1} e_k\trsp M^{-1} U$, and $e_k \trsp\MTilde^{-1}
e_k \simeq\limitScaling^{-1} e_k\trsp M^{-1} e_k$, so $\widehat
{z}\simeq\zeta(\limitScaling)/\limitScaling$. Hence,
\[
\gamma\trsp\widehat{w}_b\widehat{z} \simeq\gamma\trsp w \frac
{\zeta(\limitScaling)}{\limitScaling} .\vadjust{\goodbreak}
\]
In other words, we have found an asymptotically consistent way of
estimating the quantities of interest.
Hence, the estimator we propose to use is
\begin{equation}\label{eq:FirstEstimator}
\widehat{\wTheo}=(\wEmp+\kappa_n \widehat{z} \widehat
{w}_b)=\SigmaHat^{-1}\VHat\MTilde^{-1} U .
\end{equation}
Interestingly, this proposal does not require us to estimate
$\limitScaling$. Furthermore, because we have consistency of the
estimator in the whole class of elliptical distributions, this
estimator is fairly robust to distributional assumptions about the
data. Finally, the estimator is consistent in the sense that all
(deterministic and given) linear combinations of $\widehat{\wTheo}$
are consistent for the corresponding linear combinations of $\wTheo$
(provided these linear combinations are such that the assumptions of
Theorem~\ref{thm:BiasInWeights} apply to them). (Naturally, we cannot
take a supremum over too large a class of $\gamma$'s.)

\subsubsection*{The estimator satisfies the constraints} It is nonetheless
natural to raise the following question: does the proposed estimator
satisfy the constraints of the problem? If not, our proposal would be
problematic, but it is indeed the case that our estimator satisfies the
constraints $\widehat{\wTheo}{}\trsp v_i=u_i$ for all $i \in\{1,\ldots
,k-1\}$. Naturally, the last constraint (i.e., $\widehat{\wTheo}{}\trsp
\mu=u_k=\mu_P$) is difficult to satisfy exactly because $\mu$ is
unknown, so it is also less of a concern.

Let us now briefly justify our claim concerning the satisfaction of the
equality constraints. By construction, $\wEmp$ satisfies the
constraints $\wEmp\trsp v_i=u_i$, $1\leq i \leq k-1$, so all we have
to show is that the $k\times1$ vector $\VHat\trsp\widehat{w}_b$ is
proportional to $e_k$. We recall that $\MTilde=\MHat-\kappa_n e_k
e_k\trsp$, so
\[
\widehat{w}_b=\SigmaHat^{-1}\VHat(\MHat-\kappa_n e_k
e_k\trsp)^{-1} e_k .
\]
Using the standard formula for the inverse of a rank-1 perturbation of
a matrix, we therefore get
\begin{eqnarray*}
\widehat{w}_b&=&\SigmaHat^{-1}\VHat\biggl(\MHat^{-1}+\kappa_n \frac
{\MHat^{-1}e_k e_k\trsp\MHat^{-1}}{1-\kappa_n e_k\trsp\MHat^{-1}
e_k} \biggr) e_k\\
&=&\SigmaHat^{-1}\VHat\MHat^{-1} e_k+\kappa_n \SigmaHat^{-1}\VHat
\MHat^{-1} e_k \frac{e_k\trsp\MHat^{-1}e_k}{1-\kappa_n e_k\trsp
\MHat^{-1} e_k}\\
&=&\frac{1}{1-\kappa_n e_k\trsp\MHat^{-1} e_k} \SigmaHat^{-1}\VHat
\MHat^{-1} e_k.
\end{eqnarray*}
Once we recall that $\MHat=\VHat\trsp\SigmaHat^{-1} \VHat$, we
immediately get the equality
\[
\VHat\trsp\widehat{w}_b=\frac{1}{1-\kappa_n e_k\trsp\MHat^{-1}
e_k} e_k ,
\]
which shows that $v_i\trsp\widehat{w}_b=0$ for $1\leq i \leq k-1$, as
announced.

Finally, from a practical point of view, one might be worried that the
estimator proposed in equation \eqref{eq:FirstEstimator} ``puts too
much weight on the theory and not enough on the data'' and that better
practical performance might be achieved by tuning more finely our
corrections to the data. For instance, one might propose, we think
reasonably, to use, instead of $\MTilde$ the matrix $\MTilde(\lambda
_1)=\MHat-\lambda_1 \kappa_n e_k e_k\trsp$, where $\lambda_1$
would be picked by some form of cross-validation based on the new
estimator $\widehat{\wTheo}(\lambda_1)=\wEmp+\kappa_n \widehat
{z}(\lambda_1)\widehat{w}_b(\lambda_1)$. We do not discuss this
issue any further in this paper as we plan to address it in another,
more applied, article. We do, however, show the performance of our
estimator in limited simulations in Section~\ref{subsec:NumericalResults}.

\subsection{Improved estimation of the frontier}\label
{subsec:EstimationEffFrontier}
We now discuss the question of improved estimation of the efficient
frontier. This is naturally an important quantity in the problem, and,
as we hope to have shown, a difficult one to estimate by naive methods.
One aspect of its importance is that it gives us a benchmark of
performance for optimal portfolios. We therefore think that in a
financial context, it might be of great interest in particular to regulators.
\subsubsection{Estimation of ${\limitScaling}$}

Though we have seen that we could devise a scheme to improve the
estimation of the weights without having to estimate~$\limitScaling$,
this latter quantity is still an important one to estimate if we want
to better understand the pitfalls we might be facing.

In the elliptical case, where $X_i=\mu+\lambda_i \Sigma^{1/2} Y_i$,
we wish to estimate $\lambda_i^2$, as we have seen that $\limitScaling
$ is ``driven'' by this quantity. We now describe heuristics that
suggest how to estimate $\limitScaling$; more detailed consistency
arguments follow in Proposition~\ref
{prop:ConsistencyEstimateLimitScaling}. To estimate $\limitScaling$,
we recall that standard concentration of measure results (see below)
say that with very high probability, if the largest eigenvalue of
$\Sigma$ does not grow too fast,
\[
\frac{\| \Sigma^{1/2}Y_i\|^2_2}{p}\simeq\frac{\mathrm
{trace}(\Sigma)}{p} .
\]
Hence, in this setting, the concentration of measure phenomenon can be
used for practical purposes.
Now, note that $\| \mu-\muHat\|_2^2\simeq\frac{\mathrm
{trace}(\Sigma)}{n}$, because under our assumptions A1--A4 and the
assumption of
independence of the $\lambda_i$'s, $\sum_{i=1}^n \lambda
_i^2/n\tendsto1$ and A5 implies that the previous approximation holds. Hence,
\[
\frac{\| X_i-\muHat\|_2^2}{p}\simeq\lambda_i^2 \frac{\mathrm
{trace}(\Sigma)}{p} .
\]
We now propose the following estimator for $\lambda_i^2$:
\[
\widehat{\lambda_i^2}=\frac{\| X_i-\muHat\|_2^2}{\sum_{i=1}^n
\| X_i-\muHat\|_2^2/n}=\frac{\| X_i-\muHat\|_2^2}{\mathrm
{trace}(\SigmaHat)} .
\]
If we denote $\rho_n=p/n$, we then propose to estimate $\limitScaling
$ using the positive solution of
\[
\hat{g}_n(x)=1-\rho_n  \qquad
\mbox{where } \hat{g}_n(x)=\frac{1}{n}\sum_{i=1}^n \frac
{1}{1+x\widehat{\lambda_i^2} \rho_n } .
\]
We note that this is just the discretized version of the equation
characterizing~$\limitScaling$. ($\hat{g}_n$ is clearly a continuous
convex decreasing function of $x$ on $[0,
\infty)$, so the existence and uniqueness of a solution to
$g(x)=1-\rho_n$ is clear.)
\subsubsection{Estimation of the efficient frontier}
We recall an important result from Theorem \ref
{thm:PositionFrontierEllipticalAndGaussianCase}: under the assumptions
made in this section,
\[
\fEmp\simeq\frac{1}{\limitScaling}\fTheo-\kappa_n \frac{
(e_k\trsp M^{-1} U/\limitScaling)^2}{1+ ({\kappa
_n/\limitScaling}) e_k\trsp M^{-1} e_k} .
\]
Now recall that we have a consistent estimator of $e_k\trsp
M^{-1}/\limitScaling$, that is, $e_k\trsp\MTilde^{-1}$, and we just
discussed how to estimate $\limitScaling$.

As an estimator of the efficient frontier we therefore propose
\[
\widehat{\fTheo}=\widehat{\limitScaling} \biggl(\fEmp+\kappa_n
\frac{(e_k\trsp\MTilde^{-1}U)^2}{1+\kappa_n e_k\trsp\MTilde
^{-1}e_k} \biggr) .
\]

We also note that $\MTilde$ could be replaced by $\MTilde(\lambda
_1)$ described above with a similar cross-validation
scheme.\subsubsection{Consistency of the estimator of $\limitScaling
$}\label{subsubsec:ConsistencyEstimateLimitScaling}
Let us now show that our proposed estimator of $\limitScaling$ is
consistent. We place ourselves in the setting where $\lambda_i$'s are
i.i.d. with a second moment and $\mathbf{E}(\lambda_i^2)=1$. Recall also
that the $Y_i$'s that appear below are such that $Y_i\sim{\cal
N}(0,\id_p)$.

We have the following proposition.
\begin{propositionNumThm}\label{prop:ConsistencyEstimateLimitScaling}
Let us call $u_n= (2\log(n)+(\log n)^{\eps} )^{1/2}+\sqrt
{2\pi}$ and $\opnorm{\Sigma}$ the largest eigenvalue of $\Sigma$.
Then we have, with probability going to 1,
\[
\max_{1\leq i \leq n} \biggl|\frac{\| \Sigma^{1/2} Y_i\|
_2^2}{p}-\frac{\operatorname{trace}(\Sigma)}{p} \biggr|\leq
\frac{\opnorm{\Sigma}}{p}(4+u_n^2)+2u_n\sqrt{\frac{\opnorm{\Sigma
}}{p}}\sqrt{\frac{\operatorname{trace}(\Sigma)}{p}} .
\]
Further, if $\hat{\limitScaling}_n$ is the solution of $\hat
{g}_n(x)=1-\rho_n$,
\[
\hat{\limitScaling}_n\tendsto\limitScaling\qquad \mbox{in probability} ,
\]
as soon as, for some $\eps>0$, $\frac{\opnorm{\Sigma}}{\mathrm
{trace}(\Sigma)}u_n^2\tendsto0$ as $n\tendsto\infty$.
\end{propositionNumThm}

\begin{pf}
Let us consider the function $F(Y)=\| \Sigma^{1/2}Y\|_2/\sqrt{p}$.
Clearly this function is $\| \Sigma^{1/2}\|_2/\sqrt{p}$-Lipschitz
with respect to Euclidian norm in $\mathbb{R}^p$.

Now suppose that $Y_0\sim{\cal N}(0,\id_p)$. Let us call $m_F$ a
median of $F(Y_0)$. Using standard results on the concentration
properties of Gaussian random\vadjust{\goodbreak} variables [see \cite{ledoux2001},
Chapter 1 and Theorem 2.6], we have
\[
P \biggl( \biggl|\frac{\| \Sigma^{1/2}Y_0\|_2}{\sqrt{p}}-m_F
\biggr|>t \biggr)\leq2 \exp\biggl(-\frac{p t^2}{2\opnorm{\Sigma}}\biggr) .
\]
Hence, using a simple union bound argument, we have, after some
algebra, if $t_n=\sqrt{\opnorm{\Sigma}/p}(2\log(n)+\log(n)^{\eps})^{1/2}$,
\[
P \biggl(\max_{1\leq i \leq n} \biggl|\frac{\| \Sigma^{1/2}Y_i\|
_2}{\sqrt{p}}-m_F \biggr|>t_n \biggr)\leq2 \exp\bigl(-(\log
(n))^{\eps}/2\bigr) .
\]
So with large probability,
\[
\max_{1\leq i \leq n} \biggl|\frac{\| \Sigma^{1/2}Y_i\|_2}{\sqrt
{p}}-m_F \biggr|\leq t_n .
\]
Now, if we call $\mu_F=\mathbf{E}(F(Y_0))$, we have, using
Proposition 1.9
in \cite{ledoux2001},
\begin{eqnarray*}
|m_F-\mu_F|
&\leq&\sqrt{\pi}\sqrt{\frac{2\opnorm{\Sigma}}{p}} ,\\
0&\leq&\frac{\operatorname{trace}(\Sigma)}{p}-\mu_F^2
\leq4 \frac
{\opnorm
{\Sigma}}{p} .
\end{eqnarray*}
Now, using the fact that $\max_{1\leq i \leq n}|a_i^2-b^2|\leq\max
_{1\leq i \leq n}|a_i-b|(2b+\break \max_{1\leq i \leq n}|a_i-b|)$, and the
fact that $\mu_F \leq\sqrt{\operatorname{trace}(\Sigma)/p}$, we have
\begin{eqnarray*}
&&\max_{1\leq i \leq n} \biggl|\frac{\| \Sigma^{1/2}Y_i\|
_2^2}{p}-\frac{\operatorname{trace}(\Sigma)}{p} \biggr|\\
&&\qquad\leq\max_{1\leq i \leq n}\biggl |\frac{\| \Sigma^{1/2}Y_i\|_2}{\sqrt
{p}}-\mu_F \biggr| \biggl(2\sqrt{\operatorname{trace}(\Sigma)/p}+\max
_{1\leq i
\leq n}\biggl |\frac{\| \Sigma^{1/2}Y_i\|_2}{\sqrt{p}}-\mu
_F \biggr| \biggr) .
\end{eqnarray*}
Our previous results imply that with large probability,
\[
\max_{1\leq i \leq n} \biggl|\frac{\| \Sigma^{1/2}Y_i\|_2}{\sqrt
{p}}-\mu_F \biggr|\leq\sqrt{\frac{\opnorm{\Sigma}}{p}} u_n ,
\]
and therefore, with large probability,
\begin{eqnarray*}
&&\max_{1\leq i \leq n}\biggl |\frac{\| \Sigma^{1/2}Y_i\|
_2^2}{p}-\frac{\operatorname{trace}(\Sigma)}{p}\biggr |\\
&&\qquad\leq4 \frac
{\opnorm{\Sigma}}{p}+
\sqrt{\frac{\opnorm{\Sigma}}{p}} u_n \Biggl(2\sqrt{\frac{\mathrm
{trace}(\Sigma)}{p}}+\sqrt{\frac{\opnorm{\Sigma}}{p}} u_n \Biggr)
\\
&&\qquad=\frac{\opnorm{\Sigma}}{p}(4+u_n^2)+2u_n\sqrt{\frac{\opnorm
{\Sigma}}{p}}\sqrt{\frac{\operatorname{trace}(\Sigma)}{p}} ,
\end{eqnarray*}
as announced in the proposition.
As a consequence, we have, if we call $v_n=\sqrt{\opnorm{\Sigma
}/\operatorname{trace}(\Sigma)}u_n$, with large probability,
\[
\max_{1\leq i \leq n} \biggl|\frac{\| \Sigma^{1/2}Y_i\|
_2^2}{\operatorname{trace}(\Sigma)}-1 \biggr|\leq
2v_n+v_n^2+\frac{\opnorm{\Sigma}}{\operatorname{trace}(\Sigma)} .
\]
Hence, when for some $\eps>0$, $v_n$ goes to zero, which implies that
$\frac{\opnorm{\Sigma}}{\operatorname{trace}(\Sigma)}\tendsto0$, we
have with
high probability,
\[
\max_{1\leq i \leq n}\biggl |\frac{\| \Sigma^{1/2}Y_i\|
_2^2}{\operatorname{trace}(\Sigma)}-1 \biggr|\tendsto0 .
\]

$\bullet$ \textit{Consistency of ${\hat{\limitScaling}_n}$.}

We now assume that $v_n$ goes to zero for some $\eps>0$ and turn to
showing the consistency of $\hat{\limitScaling}_n$.
First, let us note that
\[
\muHat-\mu|\{\lambda_i\}_{i=1}^n\sim\sqrt{\frac{1}{n}\sum
_{i=1}^n \lambda_i^2} \frac{1}{\sqrt{n}}{\cal N}(0,\Sigma) .
\]
Hence, by the same concentration arguments we just used, and using the
fact that the $\lambda_i$'s are i.i.d. with $\mathbf{E}(\lambda
_i^2)=1$, we have
\[
\| \muHat-\mu\|_2^2\tendsto\sqrt{p/n} \operatorname{trace}(\Sigma)/p .
\]
Now,
\begin{eqnarray*}
\biggl|\frac{\| X_i-\muHat\|_2^2}{\operatorname{trace}(\Sigma)}-\lambda
_i^2 \biggr|&\leq&\biggl|\frac{\| X_i-\mu\|_2^2}{\mathrm
{trace}(\Sigma)}-\lambda_i^2 \biggr|+\frac{\| \mu-\muHat\|
_2^2}{\operatorname{trace}(\Sigma)} \\
&=&\lambda_i^2 \biggl|\frac{\| \Sigma^{1/2}Y_i\|_2^2}{\mathrm
{trace}(\Sigma)}-1 \biggr|+\frac{\| \mu-\muHat\|_2^2}{\mathrm
{trace}(\Sigma)} .
\end{eqnarray*}
Also, the law of large numbers (for triangular arrays) imply that with
probability 1
\[
\frac{\operatorname{trace}(\SigmaHat)}{\operatorname{trace}(\Sigma)}\tendsto
1 .
\]
So we can write
\[
|\widehat{\lambda_i^2}-\lambda_i^2 |\leq\frac{\mathrm
{trace}(\Sigma)}{\operatorname{trace}(\SigmaHat)} \biggl|\frac{\|
X_i-\muHat\|_2^2}{\operatorname{trace}(\Sigma)}-\lambda_i^2
\biggr|+\lambda_i^2
\biggl|\frac{\operatorname{trace}(\Sigma)}{\operatorname{trace}(\SigmaHat
)}-1 \biggr| ,
\]
and we have now all the terms on the right-hand side under control.

In particular, it is clear that when $v_n\tendsto0$,
\[
\frac{1}{n}\sum_{i=1}^n |\widehat{\lambda_i^2}-\lambda
_i^2 |\tendsto0 .
\]

With all these preliminaries behind us, let us now turn to the final
part of the proof. Let us call
\[
g_n(x)=\frac{1}{n}\sum_{i=1}^n \frac{1}{1+x\lambda_i^2\rho_n} .
\]
With a slight abuse of notation, we note that $\limitScaling$ is the
solution of $g_{\infty}(x)=1-\rho$.
If $x_n$ is the solution of $g_n(x_n)=1-\rho_n$, it is clear that
$x_n$ is consistent for $\limitScaling$: we can just use the fact that
$g_n$ is decreasing and evaluate it at $y_1$ and $y_2$ which are on
either sides of $\limitScaling$. Clearly $g_n(y_1)$ is consistent for
$g_{\infty}(y_1)$, and similarly for $y_2$, so with high probability,
$x_n$ needs to be in $[y_1,y_2]$ asymptotically.

Recall that the roots we are looking for are positive. So we have
\[
\hat{g}_n(x)-g_n(x)=\frac{1}{n}\sum_{i=1}^n\frac{x\rho_n (\lambda
_i^2-\widehat{\lambda_i^2})}{(1+x\lambda_i^2 \rho_n)(1+x\widehat
{\lambda_i^2} \rho_n)} ,
\]
and therefore, for $x>0$,
\[
|\hat{g}_n(x)-g_n(x) |\leq x\rho_n \frac{1}{n}\sum
_{i=1}^n |\lambda_i^2-\widehat{\lambda_i^2} | .
\]
By noting that $\hat{g}_n(\hat{\limitScaling}_n)=1-\rho
_n=g_n(x_n)$, we have
\[
|\hat{g}_n(x_n)-\hat{g}_n(\hat{\limitScaling}_n) |\leq
x_n\rho_n \frac{1}{n}\sum_{i=1}^n |\lambda_i^2-\widehat
{\lambda_i^2} |\tendsto0 ,
\]
since $x_n$ is bounded above.

Now since $\hat{g}_n(x)$ is decreasing and is pointwise consistent for
$g_{\infty}(x)$, it is clear that we can find $y_3$, deterministic and
bounded away from $\infty$, such that asymptotically, $\hat
{\limitScaling}_n<y_3$, with high probability. Also, $\hat{g}_n(x)$
is convex, so this guarantees that $|\hat{g}'_n|$ can be bounded below
(uniformly in $n$ with high probability) on $[0,\max(y_3,y_2)]$ by a
quantity that is strictly greater than 0 with high probability. Note
that this latter interval contains both $x_n$ and $\hat{\limitScaling
}_n$ asymptotically. Using the mean value theorem, the fact that we
have a lower bound (different from 0) on $|\hat{g}'_n(x)|$ on $[0,\max
(y_3,y_2)]$, and the equation in the previous display, we can finally
conclude that
\[
x_n-\hat{\limitScaling}_n\tendsto0
\]
with high probability, and since $x_n$ is consistent for $\limitScaling
$, so is $\hat{\limitScaling}_n$.
\end{pf}

%
\subsubsection{On robust estimates of scatter}
We just saw that we could take advantage of the high-dimensionality of
the problem to essentially estimate $\lambda_i^2$, by using
concentration of measure arguments. This also allows us to propose
estimates of scatter that are tailored for high-dimensional problems.

In low-dimension, estimation of individual $\lambda_i^2$ is not
possible and a classic proposal for estimating the scatter matrix
$\Sigma$ is Tyler's estimator\vadjust{\goodbreak} [see \cite
{TylerRobustEstimateScatterAoS87}], which is the solution $V_n$
(defined up to scaling), of the equation
\[
p \sum_{i=1}^n \frac{(X_i-\mu)(X_i-\mu)\trsp}{(X_i-\mu)\trsp
V_n^{-1}(X_i-\mu)} =V_n .
\]
It has been observed in a random matrix context [see \cite
{FrahmJaekel05} and \cite{BiroliBouchaudPottersStudentEnsemble07}]
that when using Tyler's estimator in connection with elliptically
distributed data, one seemed to recover a spectrum that looked similar
to predictions of the Mar\v{c}enko--Pastur law, at least in the case
of $\id_p$
scatter. At this point, the evidence is mostly based on simulations
though a rigorous proof seems feasible with a little bit of effort (the
argument given in \citep{BiroliBouchaudPottersStudentEnsemble07} is
interesting though it falls short of a ``full proof,'' which is
acknowledged in that paper). We do not try to give a proof here because
this is quite far from being the topic of this paper.

As a high-dimensional alternative to Tyler's estimator, we could use
\[
\widetilde{V}_n=\frac{1}{n}\sum_{i=1}^n \frac{(X_i-\mu)(X_i-\mu
)\trsp}{\widehat{\lambda_i^2}} .
\]
One potential advantage of this proposal over Tyler's estimator is that
Tyler's estimator is a priori not-defined when $p>n$, because it
becomes impossible to invert $V_n$. Also, this estimator is rather
quick to compute and does not require multiple inversions of $p\times
p$ matrices, where $p$ is large [Tyler's estimator is generally found
through an iterating procedure---see \cite{FrahmJaekel05} and
references therein]. The spectral properties of $\widetilde{V}_n$ are
also quite easy to analyze in light of the detailed work we carried out
concerning consistency of our estimator of $\limitScaling$. For
instance in the simple case where $\mu$ is known, it is easy to see
that under some conditions on $\Sigma$ and the $\lambda_i$'s, the
limiting spectral distribution of $\widetilde{V}_n$ will satisfy a
Mar\v{c}enko--Pastur-type equation. (Because this is really
tangential to our main
points in the paper, we do not give further details.)

Note that these estimates of scatter essentially make the influence of
the $\lambda_i$'s on the problem disappear, at least as far as
covariance (or really scatter) is concerned. So to answer a question
asked by an insightful referee, it is reasonable to think that another
approach might be to turn the problem back to an essentially Gaussian
problem by using an estimate of scatter instead of an estimate of
covariance---if we ignore problems due to mean estimation. Since in the
Gaussian case, $\limitScaling=1/(1-\rho)$, corrections are relatively
easy then. However, the impact of mean estimation needs to be
investigated and furthermore, at this point there are no rigorous
results that we know of (only very limited simulations) concerning the
spectral properties of Tyler's estimator in high-dimension. So we leave
further investigations of the properties of these estimates of scatter
to future work, as they are not a primary concern in this already long
paper (after all we have a provably consistent estimator that takes
care of all the problems and is fast to compute).

Let us however note that using estimates of scatter (instead of
covariance) would likely yield a serious improvement in terms of the
realized risk of portfolios which is discussed in the paper \citep
{nekMarkoRealizedRisk}. However, these questions touch more on the
issue of allocation, whereas we are concerned in this paper with
estimating the efficient frontier and have shown that we can do this
well (at least asymptotically and theoretically) independently of
allocation issues, a fact that is potentially useful for, for instance,
creating benchmarks.

\subsection{Numerical results and practical considerations}\label
{subsec:NumericalResults}
This subsection gives some numerical results to assess the quality of
the proposed estimators for both weights and ``efficient frontier.'' The
simulation analysis is done in an a priori quite favorable case---the
question being whether even then the theory could be useful in practice.

Our aim was to investigate among other things the improvement in the
quality of our approximations as $n$ and $p$ grew to infinity. Hence,
we present the results of two simulation setups: one where $n=250,
p=100$ and one where $n=2500,p=1000$. We chose to work with simulations
where we picked both $\Sigma$ and $\mu$ so that we could guarantee,
for instance, that the efficient frontier was basically the same for
both simulations.

More specifically, we chose $\Sigma$ to be a $p\times p$ Toeplitz
matrix, with $\Sigma(i,j)=\alpha^{|i-j|}$, where $\alpha=0.4$. In the
smaller dimensional simulation, that is, $p=100$, we picked $v_1$ to be
the eigenvector associated with the 90th smallest eigenvalue of $\Sigma
$. Calling $\beta_2$ the eigenvector associated with the 15th smallest
eigenvalue of $\Sigma$, we picked $v_2=\mu$ to be $\sqrt
{0.3}v_1+\sqrt{0.7}\beta_2$. In the larger dimensional simulation, we
used for $v_1$ the eigenvector associated with the 900th smallest
eigenvalue of $\Sigma$, while $\beta_2$ was now associated with the
150th smallest eigenvalue of $\Sigma$. $\mu=v_2$ was computed in the
same fashion in both simulations.

The simulations are here to illustrate ``how large is large,'' that is,
when the asymptotics kick-in and our theoretical predictions become
accurate. The parameters were chosen so that we would be close to
satisfying assumptions A1--A5. Also, the choice of $v_1$ and $v_2$
guarantees that the off-diagonal elements of $M$ are not zero, which we
thought might make the problem easier and lead to overoptimistic
pictures. (This choice of parameters is not motivated by a particular
problem in Finance. We also note that if we knew that
the covariance matrix were Toeplitz, we could resort to regularization
methods to better solve the problem. However, if we applied the same
random rotation to $\Sigma$, $v_1$ and $v_2$, it becomes less clear
how one could use other approaches than the ones presented here for estimation.)

We did simulations both in the Gaussian case and in the case of an
elliptical distribution as described above, that is, $X_i=\mu+\lambda
_i \Sigma^{1/2} Z_i$,\vadjust{\goodbreak} where $\lambda_i$ was proportional to a
$t$-distributed random variables with 6 degrees of freedom and scaled
to have variance 1. We picked 6 degrees of freedom to have simulations
with relatively heavy tails and capture visually the corresponding
effects. It was also naturally a way to investigate the practical
robustness of our estimators and compare with the Gaussian case. We
call below the set of simulations involving the $t$-distribution the
``$t_6$'' case because of its similarity with multivariate $t$-distributions.

We repeated 1000 times the simulations in all the cases considered. We
chose $u_1=1$ and $u_2$ (the ``target returns'' in a financial context)
ranging from 0.1 to 5.

We note that our estimators require taking inverses of matrices which
naturally raises the question of how well conditioned those matrices
are. This is particularly the case when we deal with $M$ and $\MTilde
$: if $M$ is poorly conditioned, even though $\MTilde$ is a good
estimator of $\limitScaling M$, it can turn out that $\MTilde^{-1}$ is
a relatively poor estimator of $M^{-1}\limitScaling^{-1}$. In our
simulations, both $M$ and $\Sigma$ were well conditioned but in
practice, one should be aware of potential difficulties that may arise
if, for instance, $\MTilde$ indicates that $M$ may be ill conditioned.
When this is the case, it is actually quite easy to make the estimators
perform poorly (but of course this violates assumptions
A1--A5).\vspace*{-2pt}

%
\begin{figure}

\includegraphics{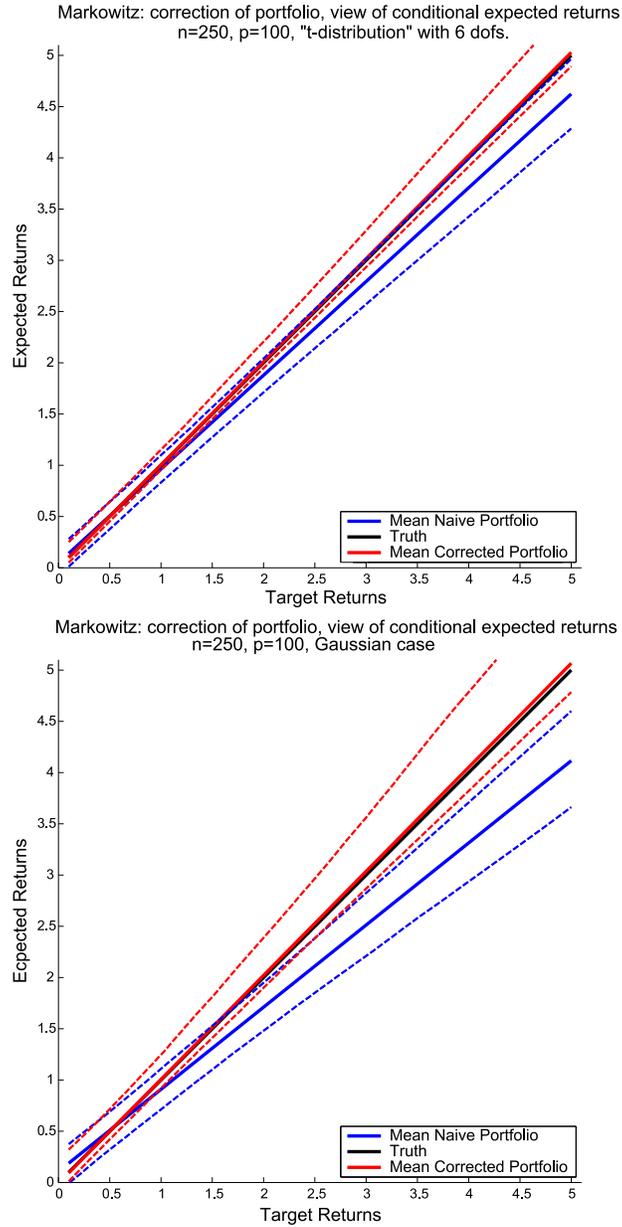}

\caption{Performance of naive and
corrected portfolios, for scaled ``$t_6$'' (top picture) and Gaussian
returns. Here $n=250, p=100$ and the number of simulations is 1000. The
dashed lines represent 95\% confidence bands. The $x$-axis represents
the returns an investor expects. The $y$-axis represents what she would
actually get on average (i.e., $\mu\trsp\widehat{w}$). The plots show
both the bias in the naive solution (blue solid lines) and the fact
that our estimator is nearly unbiased (red solid lines). They also
illustrate the robustness of our corrections. The black line is very
close to the red line, showing a very good correction (on average) in
this setting where assumptions \textup{A1--A5} are satisfied.}\label
{fig:ReturnsCorrectionSmall}
\end{figure}

%
\begin{figure}

\includegraphics{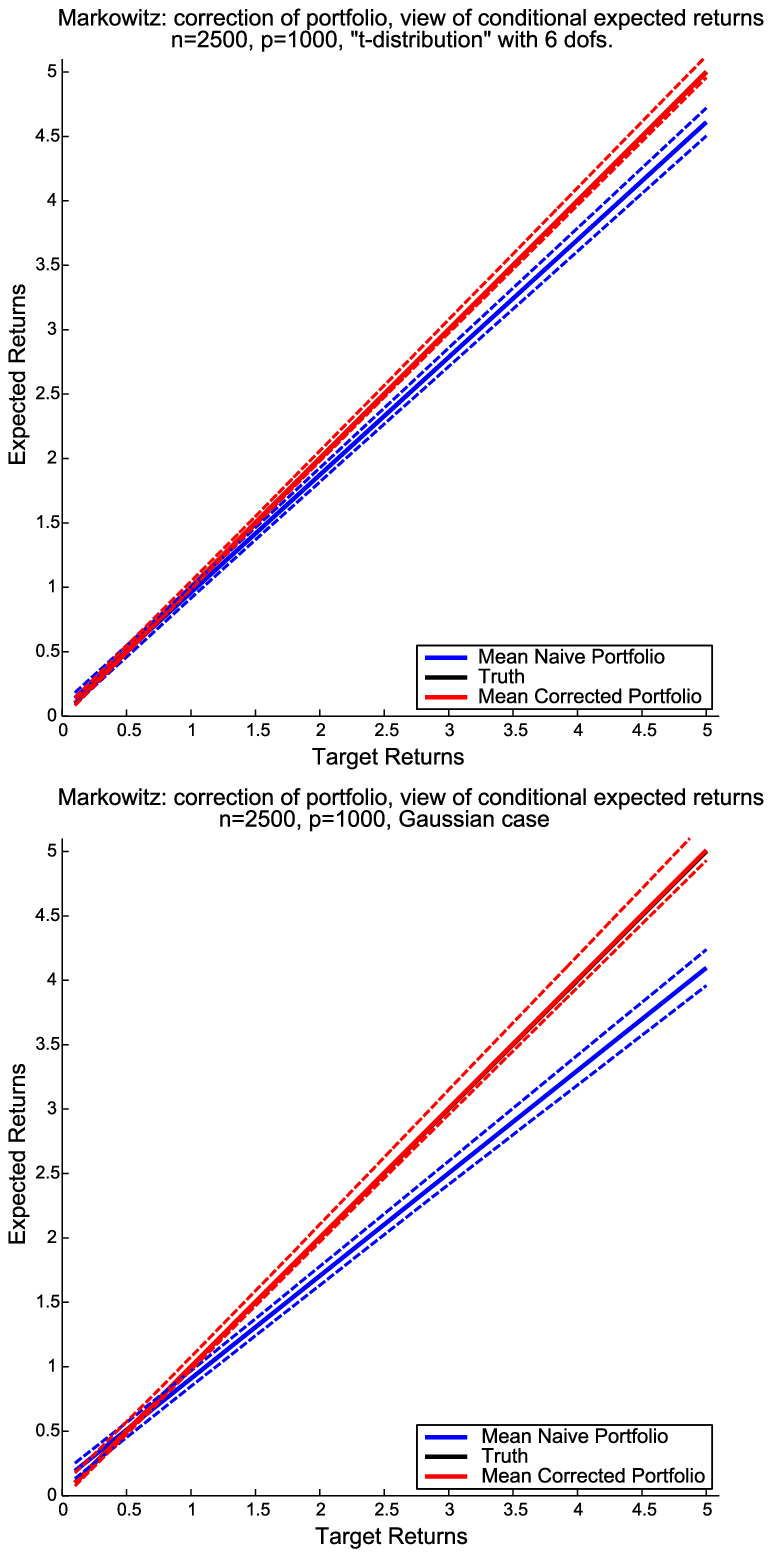}

\caption{Performance of naive and
corrected portfolios, for scaled ``$t_6$'' (left picture) and Gaussian
returns. Here $n=2500, p=1000$ and the number of simulations is 1000.
The dashed lines represent 95\% confidence bands. The $x$-axis
represents the returns an investor expects. The $y$-axis represents
what she would actually get on average (i.e., $\mu\trsp\widehat
{w}$). The plots show both the bias in the naive solution (blue solid
lines) and the fact that our estimator is nearly unbiased (red solid
lines). They also illustrate the robustness of our corrections. Note
the narrower confidence bands as compared to Figure \protect\ref
{fig:ReturnsCorrectionSmall}. The black line is essentially hidden
under the red line, showing a near perfect correction (on average) in
this setting where assumptions \textup{A1--A5} are satisfied.}\label
{fig:ReturnsCorrection}
\end{figure}

%
\begin{figure}[t!]
\centering
\begin{tabular}{@{}cc@{}}

\includegraphics{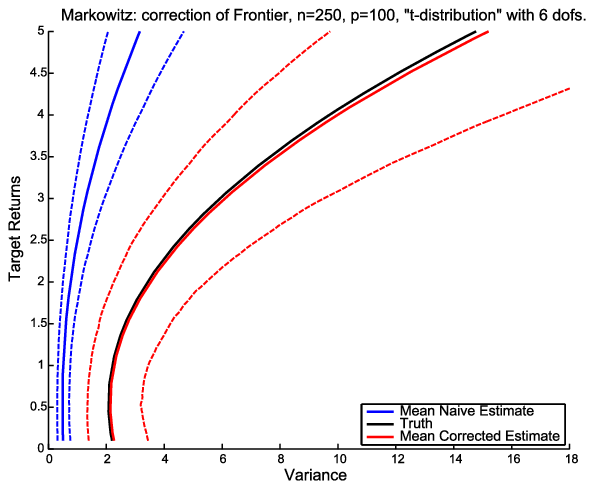}
&\includegraphics{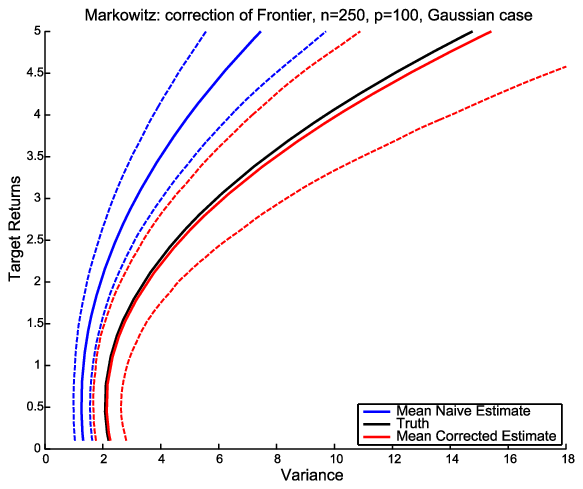}\\
\footnotesize{(a)}&\footnotesize{(b)}\\[8pt]

\includegraphics{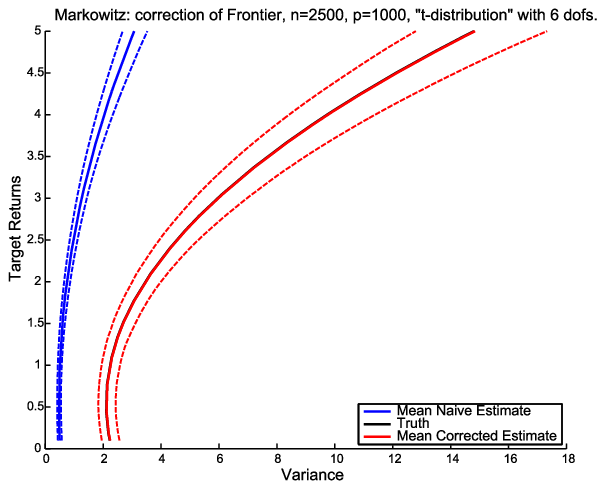}
&\includegraphics{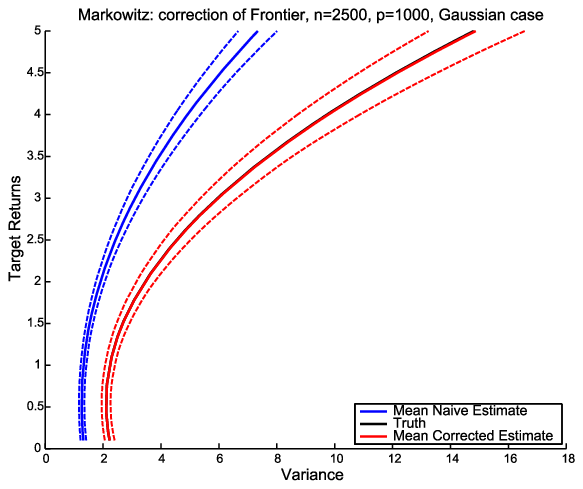}\\
\footnotesize{(c)}&\footnotesize{(d)}
\end{tabular}
\caption{Performance of naive and
corrected frontiers, for scaled ``$t_6$'' [\textup{(a)} and \textup{(c)}]
and Gaussian returns [\textup{(b)} and \textup{(d)}]. Here, in the left column $n=250$
and $p=100$. In the right column, $n=2500, p=1000$. The number of
simulations is 1000 in all pictures. The dashed lines represent
(empirical) 95\% confidence bands. (The confidence bands corresponds
are computed for a fixed $y$.) The $x$-axis represents our estimate of
variance of the optimal portfolio. The $y$-axis represents the target
returns for the portfolio. The plots show both the bias in the naive
solution (blue solid curves) and the fact that our estimator is nearly
unbiased (red solid curves near, or covering the black curve, the
population solution). They also illustrate the robustness of our
corrections. Another striking feature is the lack of robustness of
Gaussian computations, since the ``efficient frontiers'' computed with
``$t_6$'' returns are different from the Gaussian ones. The fact that,
as our theoretical work predicts, Gaussian computations underestimate
risk-underestimation in the class of elliptical distributions
considered in the paper is illustrated by the fact that the ``$t_6$''
curves are to the left of the Gaussian curves. Note the narrower
confidence bands in the larger dimensional simulations [\textup{(c)} and \textup{(d)}].
The black line is essentially hidden under the red line in \textup{(c)} and \textup{(d)},
showing a near perfect correction (on average) in this setting where
assumptions \textup{A1--A5} are satisfied.}\label{fig:FrontiersCorrection}\vspace*{12pt}
\end{figure}

\subsubsection{Estimation of portfolio weights}
As we have seen earlier, the ``naive'' weights obtained by plugging-in
the sample mean and the sample covariance matrix in our quadratic
program with linear equality constraints are biased, in the sense that
their projection in any given direction will generally be biased.

Here we show the performance of our estimator as measured by its
projection on $v_k=\mu$. It is a natural direction to consider since,
for instance, in a financial context and under our modeling
assumptions, it gives us the expected returns of our portfolio
(conditional on $X_1,\ldots,X_n$).

As our limited simulations indicate, our estimator appears to be
practically unbiased here (even in the ``lower-dimensional'' case),
which means in a financial context that the corresponding investment
strategy will yield the returns that the investor expected. (We note
that from a mean--variance point of view, we do not claim that our
estimator is optimal. Work is under way to find better performing
portfolios---but it requires a new set of theoretical investigations
whose results are postponed to another paper. In limited simulations,
it appeared that our ``debiased'' portfolio performed similarly to the
naive one from a mean--variance point of view, its main advantage being
that it delivers the returns that the investor expects.)

We present two pictures, Figure~\ref{fig:ReturnsCorrectionSmall},
  page \pageref{fig:ReturnsCorrectionSmall}
and Figure~\ref{fig:ReturnsCorrection}, page \pageref{fig:ReturnsCorrection} to give a sense to the reader of
the impact of the size of $n$ and $p$ on the estimators we proposed
[the ``larger-dimensional'' case gives quite significantly better
results, with narrower confidence bands, though (empirical)
near-unbiasedness is present in both cases].\vadjust{\goodbreak}

\subsubsection{Correction to the frontier}
We now turn to the issue of estimating the ``efficient frontier,'' that
is, the curve that represents the minima of our convex optimization
problem (\ref{eq:GeneralQP}), on page \pageref{eq:GeneralQP}. The
pictures we present on Figure~\ref{fig:ReturnsCorrection} (see page \pageref{fig:FrontiersCorrection}) were
obtained from the simulations we described above. We chose to plot the
variance (i.e., $\min w\trsp\Sigma w$) on the $x$-axis and the ``target
returns'' [i.e., the $u_k$'s in the notation of equation \eqref
{eq:GeneralQP}] on the $y$-axis as this is the convention in financial
applications.

As the reader can see, our estimator turns out to be essentially
unbiased, even in the ``lower-dimensional'' case. We note too that the
variance can be quite large but that the confidence bands obtained from
our corrections were always to the right of the confidence bands
obtained from the naive estimator, meaning that if one is concerned
with risk estimation that in (essentially) the worst case for our
estimator, we still obtained a better performing estimator than in
(essentially) the best case for the naive estimator. (We do not claim
that this is always the case and it might be an artifact of the
simulation setup chosen here.)

Finally, for graphical purposes and to help comparisons, we chose to
put all the graphs on the same scale. Some of the information on our
original graphs (for the ``lower-dimensional'' case) was therefore left
out but can be inferred by ``naturally'' extrapolating the curves shown
on our graphs which are essentially parabolas.

\subsection{Remarks on inequality constraints}\label
{Subsec:InequalityConstraints}
Our work has mostly been concerned with obtaining results for the case
of a quadratic program with linear equality constraints. We now explain
that our results can also be used to obtain approximation results
concerning the case of a quadratic program with linear inequality constraints.

In this subsection we therefore consider the problem
\renewcommand{\theequation}{QP-ineqc-Pop}
\begin{equation}\label{eq:GeneralQPIneqConstraintsPop}
\fTheo(Q)=\cases{\displaystyle
\inf_{w \in\mathbb{R}^p} w\trsp\Sigma w, \cr\displaystyle
V\trsp w \in Q .
}
\end{equation}
Here $Q$ is a subset of $\mathbb{R}^k$, and $V$ is a $p\times k$
matrix. We naturally want to relate the solution of the above problem
to that of the empirical version of the problem:
\renewcommand{\theequation}{QP-ineqc-Emp}
\begin{equation}\label{eq:GeneralQPIneqConstraintsEmp}
\fEmp(Q)=\cases{\displaystyle
\inf_{w \in\mathbb{R}^p} w\trsp\SigmaHat w, \cr\displaystyle
\VHat\trsp w \in Q .
}
\end{equation}
When $Q$ is a product of intervals, we obtain a quadratic program with
linear inequality constraints. But our formulation allows us to deal
with even more complicated constraint structures.
We note that if $G(U)$ is the solution of problem \eqref
{eq:GeneralQPIneqConstraintsPop} with $Q=\{U\}$ (i.e., a singleton),
where $U$ is a vector in $\mathbb{R}^k$, we are back in the case of
the equality constrained problem that we worked with for most of this\vadjust{\goodbreak}
paper. Let us call $\widehat{G}(U)$ the solution of problem \eqref
{eq:GeneralQPIneqConstraintsEmp} with $Q=\{U\}$. We now make the simple
following observation: note that
\begin{eqnarray*}
\fTheo(Q)&=&\inf_{U\in Q} G(U) ,\\
\fEmp(Q)&=&\inf_{U\in Q} \widehat{G}(U) .
\end{eqnarray*}

The main idea here is that we can find a deterministic equivalent to
$\fEmp(Q)$ and we can relate this deterministic equivalent to $\fTheo(Q)$.

Recall from Section~\ref{Sec:GeneralQP} that $G(U)=U\trsp M^{-1} U$
and $\widehat{G}(U)=U\trsp\MHat^{-1}U$. Recall also that under the
assumptions A1--A5 made at the beginning of this section, we have found
a deterministic equivalent to $\MHat^{-1}$: we have shown that $\MHat
\simeq\limitScaling M+\kappa_n e_k e_k\trsp=M_0(\limitScaling
,\kappa_n)$ in probability. The previous result is valid entry-wise,
and since we assume that $k$ stays bounded in the asymptotics we are
considering, it is also valid in operator norm. Now, $\MHat$ is
invertible with probability one under our assumptions, so we have,
using the first resolvent identity, that is, $A^{-1}-B^{-1}=A^{-1}(B-A)B^{-1}$,
\[
\opnorm{\MHat^{-1}-M^{-1}_0(\limitScaling,\kappa)}\leq\opnorm
{\MHat^{-1}}\opnorm{M_0^{-1}(\limitScaling,\kappa_n)}\opnorm{\MHat
-M_0(\limitScaling,\kappa_n)} .
\]
Hence, since $\opnorm{M^{-1}}$ remains bounded under our assumptions,
\[
\opnorm{\MHat^{-1}-M^{-1}_0(\limitScaling,\kappa)}\tendsto0\qquad \mbox
{in probability}.
\]

Under our assumptions, we also know that the smallest eigenvalue of
$\MHat$ and $M_0(\limitScaling,\kappa)$ stay bounded away from 0.
Therefore, for any $\delta>0$, we know that asymptotically, and with
probability 1,
\[
\forall U \in\mathbb{R}^k\qquad |U\trsp\MHat^{-1}U-U\trsp
M^{-1}_0(\limitScaling,\kappa_n)U |\leq\delta\| U\|_2^2 .
\]
Furthermore, let us note that assumption A2 guarantees that $\opnorm
{M}$ remains bounded and hence so do $\opnorm{\MHat}$ and $\opnorm
{M_0(\limitScaling,\kappa_n)}$.

We have the following theorem.

\begin{theorem}\label{thm:UniformApproxConvexMaps}
Suppose $G_0$ and $\widehat{G}_n$ are maps from $\mathbb{R}^k$ to
$\mathbb{R}^+$ such that:
\begin{enumerate}[3.]
\item$G_0$ is deterministic and $\widehat{G}_n$ is possibly random.
\item$G_0(0)=\widehat{G}_n(0)=0$.
\item$\exists c_0>0$ such that, $\forall U$, $G_0(U)\geq c_0 \| U\|
_2^2$. Similarly, $\exists\hat{c}_n>0$ such that, $\forall U$,
$\widehat{G}_n(U)\geq\hat{c}_n \| U\|_2^2$. Furthermore, $\hat
{c}_n\tendsto c_0$ with probability 1.
\item$\exists\delta_n$ such that $\delta_n\tendsto0$ in
probability and $\forall U$, $ |\widehat{G}_n(U)-G_0(U)
|\leq\delta_n \| U\|^2_2$.
\end{enumerate}
Assume that $k$ is fixed as $n\tendsto\infty$. Suppose $Q$ is a
(nonempty) subset of $\mathbb{R}^k$ and that we can find $U_0\in Q$
such that $G_0(U_0)<\infty$ and $U_0\neq0$. Then,
\[
\lim_{n\tendsto\infty} \inf_{U\in Q} \widehat{G}_n(U)=\inf_{U\in
Q} G_0(U) \qquad \mbox{in probability.}
\]
\end{theorem}

We have the following corollary:
\begin{corollary}
When assumptions \textup{A1--A5} are satisfied
\[
\fEmp(Q)\tendsto\inf_{U\in Q} U\trsp M^{-1}_0(\limitScaling,\kappa
_n)U \qquad \mbox{in probability}.
\]
\end{corollary}

Hence, we have found a deterministic equivalent to $\fEmp(Q)$. It
should also be noted that because $U\trsp M^{-1}_0(\limitScaling
,\kappa_n)U\leq\frac{1}{\limitScaling} M^{-1}$, we also have
\[
\fEmp(Q)\leq\frac{1}{\limitScaling} \inf_{U\in Q} U\trsp M^{-1}U
=\frac{1}{\limitScaling} \fTheo(Q)\qquad\mbox{with high probability}.
\]
Hence, our results on risk underestimation remain valid, even with
these more general (nonequality) linear constraints. The comparison
theorems between Gaussian and elliptical assumptions remain also valid,
because of similar comparison theorems for their deterministic
equivalents. [Note also that when $\MHat\simeq\limitScaling M$ (i.e.,
when the sample mean does not appear in $\VHat$), the previous
inequalities become equalities.] Finally, our corrections also give a
way to get a consistent estimator of $\fTheo(Q)$: one can simply solve
the optimization problem over $Q$ with $\MHat$ replaced by $1/\widehat
{\limitScaling}(\MHat-\kappa_n e_k e_k\trsp)$ in the definition of
$\widehat{G}$.

Note that the corollary follows immediately from Theorem \ref
{thm:UniformApproxConvexMaps} because of our remarks on the operator
norm of $\MHat$ and $\MHat^{-1}$ and their deterministic equivalents.

Let us now prove Theorem~\ref{thm:UniformApproxConvexMaps}.
\begin{pf*}{Proof of Theorem~\ref{thm:UniformApproxConvexMaps}}
Let us pick $U_0$ in $Q$. We can do so because $Q$ is nonempty. We
assume without loss of generality that $0\notin Q$, for otherwise the
problem is trivial, since $0$ is the global minimizer of both
(deterministic and stochastic) problems.

Let us pick $r_0=\sqrt{2 G_0(U_0)/c_0}$, with $U_0\neq0$. Suppose
that $U\notin B(0,r_0)$, where $B(0,r_0)$ is the closed ball of radius
$r_0$ with center $0$. Then, our assumptions on $G_0$ guarantee that
\[
G_0(U)\geq c_0 \| U\|_2^2>c_0 \frac{2 G_0(U_0)}{c_0}=2 G_0(U_0) .
\]
So $U_0\in B(0,r_0)$. Also, if we call $Q(r_0)=Q\cap B(0,r_0)$,
$Q(r_0)$ is nonempty and
\[
\inf_{U\in Q} G_0(U)=\inf_{U\in Q(r_0)} G_0(U) ,
\]
because if $U$ is outside of $B(0,r_0)$, $G_0(U)>G_0(U_0)$. Now,
suppose that $\{\alpha_t\}_{t\in T}$ and $\{\beta_t\}_{t\in T}$ are
two sets of real numbers. We have
\[
|\inf\alpha_j-\inf\beta_k |\leq\sup|\alpha_i-\beta
_i| .
\]
As a matter of fact, for any $j$,
\[
(\inf\alpha_k) -\beta_j\leq\alpha_j-\beta_j\leq|\alpha_j-\beta
_j|\leq\sup_k |\alpha_k-\beta_k|.
\]
Now $\sup_j[(\inf\alpha_k) -\beta_j]=(\inf\alpha_k)-(\inf\beta
_k)$. And the previous display guarantees that $\sup_j [(\inf\alpha
_k) -\beta_j]\leq\sup_k |\alpha_k-\beta_k|$. By symmetry of the
role of $\alpha$ and $\beta$, we therefore have
\[
|\inf\alpha_j-\inf\beta_k |\leq\sup|\alpha_i-\beta
_i| .
\]
Hence, we can conclude that
\[
\Bigl|\inf_{U\in Q(r_0)} G_0(U)-\inf_{U\in Q(r_0)} \widehat
{G}_n(U) \Bigr|\leq\sup_{U\in Q(r_0)} |G_0(U)-\widehat
{G}_n(U) |\leq\delta_n r_0^2 ,
\]
by our assumptions, and the fact that $\| U\|_2^2\leq r_0^2$ in $Q(r_0)$.
Hence, since $r_0$ stays fixed as $n\tendsto\infty$,
\[
\Bigl|\inf_{U\in Q} G_0(U)-\inf_{U\in Q(r_0)} \widehat
{G}_n(U) \Bigr|\tendsto0 \qquad \mbox{in probability.}
\]
If we can show that with high-probability,
\[
\inf_{U\in Q(r_0)} \widehat{G}_n(U)=\inf_{U\in Q} \widehat{G}_n(U) ,
\]
the result will be shown. First we note that if $U \notin B(0,r_0)$,
\[
\widehat{G}_n(U)\geq\hat{c}_n \frac{2 G_0(U_0)}{c_0}>(1+\delta)
\widehat{G}_n(U_0)
\]
for some $\delta>0$ with high probability under our assumptions. Let
us call $E_{\delta}$ the event $E_{\delta}=\{2\hat
{c}_n/c_0>(1+\delta)\widehat{G}_n(U_0)/G_0(U_0)\}$. Of course,
$P(E_{\delta})\tendsto1$ under our assumptions, since $2\hat
{c}_n/c_0\tendsto2$ in probability and $\widehat
{G}_n(U_0)/G_0(U_0)\tendsto1$ in probability. When $E_{\delta}$ is
true, we have
\[
\inf_{U\in B^c(0,r_0)\cap Q} \widehat{G}_n(U)\geq(1+\delta)
\widehat{G}_n(U_0)> \widehat{G}_n(U_0)\geq\inf_{U\in Q(r_0)}
\widehat{G}_n (U) .
\]
So when $E_{\delta}$ is true, and hence with high-probability,
\[
\inf_{U\in Q} \widehat{G}_n(U)=\inf_{U\in Q(r_0)} \widehat{G}_n
(U) .
\]
We can finally conclude that
\[
\inf_{U\in Q} \widehat{G}_n(U)\tendsto\inf_{U\in Q} G_0(U) \qquad \mbox{in probability,}
\]
and the theorem is proved.
\end{pf*}

\section{Conclusion}
This study of quadratic programs with linear constraints whose
parameters are estimated from data has highlighted the difficulties
created by the high-dimensionality of the data. In particular, we have
shown that the fact that $n$ (the number of observations used to
estimate the parameters) and\vadjust{\goodbreak} $p$ both grew to infinity lead to a
systematic underestimation of the minimal ``risk'' one exposed oneself
to when approaching the optimization problem \eqref{eq:QP-eqc-Pop} by
solving its naive proxy \eqref{eq:QP-eqc-Emp}.

Our study produced exact distributional results in the Gaussian case
(Section~\ref{sec:GaussianCase}) and convergence results in
probability in the elliptical case (Section~\ref{sec:EllipticalCase}),
which also allowed us to reach conclusions for the bootstrap and the
case of nonindependent data (in particular, it covers the case of
Gaussian data correlated in time). As explained in Section~\ref
{sec:ComparisonGaussianElliptical}, the study of the Gaussian case
gives an over-optimistic assessment of risk underestimation in the
context we study: in the class of elliptical distributions we consider,
risk is minimally underestimated in the Gaussian case, and the
situation is more dire for other elliptical distributions. Our study
also highlights the fact that standard bootstrap estimates of bias will
be inconsistent. It also suggests that in the case of correlated
Gaussian observations, risk underestimation is likely to be more severe
than in the i.i.d. case.


Another benefit of our analysis is that it sheds light on what is
creating those difficulties and allows us to propose robust corrections
to these problems. As shown in the theoretical part of the paper and
illustrated in our limited simulation work, they are robust in the
class of elliptical distributions we consider. They also appear to work
reasonably well in practice (when the underlying assumptions hold), as
our (somewhat limited) simulation work seems to indicate.

Perhaps surprisingly, we did not need to make very strong assumptions
about the covariance matrix at stake or the mean, whereas recent
statistical work focused on estimation of covariance matrices [see
\cite{nekSparseMatrices} or \cite{bickellevina06}] tends to do so.
This is in part because our theoretical analysis clearly showed what
functionals of these two parameters one needed to estimate, and hence
we were able to bypass stronger requirements by focusing on those
particular functionals and correcting the first order errors that
appeared. In other words, even though our aim was to estimate a
complicated function of the population covariance matrix and of the
population mean, for which we do not have good estimators in
high-dimension in general, we were able to use poor estimators of both
(and our theoretical analysis) to get an accurate estimator of the
functional of interest. This is an interesting result in the context of
high-dimensional statistics more generally, as it suggests that we
might be able to estimate certain functions of high-dimensional
parameters without having to accurately estimate the parameters
themselves [and hence we might be able to bypass in some situations
sparsity (or other similar requirements) for the population quantities].

Beside the interesting statistical and mathematical questions this
study raised, we hope that it might also be helpful to, for instance,
financial regulators by perhaps providing them with more realistic
benchmarks for the performance of optimal portfolios and that it sheds
light on how the high-dimensionality of the data affects the proper
assessment of risk of large portfolios obtained by solving
high-dimensional optimization problems.

\renewcommand{\theequation}{A.\arabic{equation}}
\setcounter{equation}{0} 
\begin{appendix}\label{appm}
%
%


\section{Classical results of linear algebra}
\subsection{On inverses of partitioned matrices}\label
{subsec:ClassicalLinearAlgebra}

In our study of the Gaussian case, and in particular in connection with
properties of Wishart matrices, we relied several times on properties
of the inverse of a partitioned matrix. Here is a detailed statement of
what we needed.

Let $A$ be a generic matrix, and let us decompose it by blocks
\[
A=
\pmatrix{
A_{11} &A_{12}\cr
A_{21}& A_{22}
}.
\]

Let us call $A^{-1}$ the inverse of $A$. We assume that all inverses we
take are well defined. Let us write
\[
A^{-1}=
\pmatrix{
A^{11} &A^{12}\cr
A^{21}& A^{22}
}.
\]
Then, it is well known that [see, e.g., \cite{mardiakentbibby}, pages
458--459, or \cite{boydvandenberghe04}, page 650]
\begin{eqnarray}
A^{11}&=& (A_{11}-A_{12}A_{22}^{-1}A_{21} )^{-1} ,\\
A^{22}&=& (A_{22}-A_{21}A_{11}^{-1}A_{12} )^{-1} ,\label
{eq:Rep22BlockInInverse}\\
A^{12}&=&-A_{11}^{-1}A_{12}A^{22} ,\\\label{eq:Rep21BlockInInverse}
A^{21}&=&-A^{22}A_{21}A_{11}^{-1} .
\end{eqnarray}

\section{Random matrix results}
\subsection{Lower bounds on smallest eigenvalue}
In many proofs in the course of the paper we needed to have
quantitative bounds on the behavior of the smallest eigenvalue of a
number of matrices and made repeated use of the following lemma.

\begin{lemma}\label{lemma:LowerBoundOnSmallestEig}
Suppose $Y$ is a $n\times p$ matrix, with i.i.d. ${\cal N}(0,1)$
entries, with $p/n\tendsto\rho$, and $0<\rho<1$.

Suppose $\Lambda$ is an $n\times n$ diagonal and deterministic matrix
and that we can find $N$, $C>0$ and $\eps>0$ such that, if $\tau_i$
is the $i$th largest eigenvalue of $\Lambda\trsp\Lambda$, $\tau
_{N}>C$, for some fixed $C>0$. $N$ is such that, for $p$ and $n$ large,
$p/N<1-\eps$ and $N/n$ stays bounded away from 0. Finally, we assume
that all the diagonal entries of $\Lambda$ are different from 0.

Call $H=\id-\delta\delta\trsp/n$, where $\| \delta\|_2^2=n$.
Then $\lambda_p$, the smallest eigenvalue of $Y\trsp\Lambda\trsp
H\Lambda Y/n-1$, is bounded away from 0 with high-probability.

In particular, when $p/N<1-\eps$, if $\mathfrak{C}_n=C \frac{N-1}{n-1}$,
\[
P \bigl(\sqrt{\lambda_p}\leq\sqrt{\mathfrak{C}_n} \bigl[\bigl(1-\sqrt
{1-\eps}\bigr)-t \bigr] \bigr)\leq\exp\bigl(-(N-1)t^2 \bigr) .
\]
\end{lemma}

The following proof makes clear that the result holds also when some of
the diagonal entries of $\Lambda$ are equal to zero if we make the
following modification: $n$ should now denote the number of nonzero
entries on the diagonal of $\Lambda$, and the corresponding
assumptions about $p$ and $N$ should then hold. We also point out that
under our assumptions $H$ is an orthogonal projection matrix.

\begin{pf*}{Proof of Lemma~\ref{lemma:LowerBoundOnSmallestEig}}
Before we start the proof per se, we need some notations: we call
$\lambda_k$ the $k$th largest eigenvalue of a symmetric matrix. In
other words, the eigenvalues are decreasingly ordered and $\lambda
_1\geq\lambda_2 \geq\cdots.$

The result is known if $\Lambda= \id_n$, since
\[
\frac{1}{n-1}Y\trsp HY\equalInLaw\frac{1}{n-1}\Wishart_p(\id
_p,n-1) .
\]
Using \cite{DavidsonSzarek01}, Theorem II.13, we have the following
result: the smallest eigenvalue of a matrix with distribution $\Wishart
(\id_p,n_0)/n_0$ is strongly concentrated around $(1-\sqrt{p/n_0})^2$
when $p<n_0$, and
\[
P \bigl(\sqrt{\lambda_{p}}<\bigl(1-\sqrt{p/n_0}\bigr)-t \bigr)\leq\exp
(-n_0t^2) .
\]

This gives our result in the case where $\Lambda=\id_n$. Let us now
investigate what happens when $\Lambda$ is not $\id_n$.

The matrix $M=\Lambda\trsp H \Lambda$ is a rank-1 perturbation of
$\Lambda\trsp\Lambda$ and is positive semi-definite, because $H$ is.
Therefore, for any $k\geq2$, $\lambda_{k-1}(\Lambda\trsp H \Lambda
)=\lambda_{k-1}(M)\geq\lambda_k(\Lambda\trsp\Lambda)$, by the
interlacing Theorem 4.3.4 in \cite{hj}. $M$ has rank $n-1$ matrix
since, $M\Lambda^{-1}\delta=0$ and $\textrm{rank}(M)\geq\textrm
{rank}(\Lambda\trsp\Lambda)-1=n-1$.

We can diagonalize $M=ODO\trsp$, where $D$ has $(n-1)$ nonzero
coefficients, and because $O\trsp Y\equalInLaw Y$, we have
\[
Y\trsp M Y = Y\trsp\Lambda\trsp H \Lambda Y\equalInLaw Y\trsp D
Y=\sum_{i=1}^{n-1} d_i Y_i Y_i\trsp,
\]
where $d_i$ are the nonzero diagonal entries of $D$. Because $M$ is
positive semi-definite, we have $d_i\geq0$ for all $i$. In other
respects, because for all $k\leq n-1$, $d_{k}\geq\lambda
_{k+1}(\Lambda\trsp\Lambda)=\tau_{k+1}$ by our remark on
interlacing inequalities. Hence, we have, if $\succeq$ denotes
positive-semidefinite ordering,
\[
\sum_{i=1}^{n-1} d_i Y_i Y_i\trsp\succeq\sum_{i=1}^{N-1} d_i Y_i
Y_i\trsp\succeq\tau_N \sum_{i=1}^{N-1} Y_i Y_i\trsp= \tau_N
\Wishart_p(\id_p,N-1) .
\]
Therefore, we have in law,
\[
\frac{1}{n-1}Y\trsp\Lambda\trsp H\Lambda Y\succeq C \frac
{N-1}{n-1} \frac{1}{N-1}\Wishart_p(\id_p,N-1).
\]
As we recalled above, the smallest eigenvalue of $\Wishart_p(\id
_p,N-1)/(N-1)$ remains bounded away from 0 with high-probability in our
setting because $p/N$ remains bounded away from 1 by assumption. We
also assumed that $N/n$ and $C$ were bounded away from 0. If we call
$\mathfrak{C}=\liminf_{n\tendsto\infty} C \frac{N-1}{n-1}$, we
have $\mathfrak{C}>0$, and, for any $\eta>0$,
%
according to the result of \cite{DavidsonSzarek01} we have, for
$\lambda_p=\lambda_p (\frac{1}{n-1}Y\trsp\Lambda\trsp
H\Lambda Y )$ and $\mathfrak{C}_n=C (N-1)/(n-1)$,
\[
P \bigl(\sqrt{\lambda_p}\leq\sqrt{\mathfrak{C}_n} \bigl[\bigl(1-\sqrt
{p/(N-1)}\bigr)-t \bigr] \bigr)\leq\exp\bigl(-(N-1)t^2 \bigr) .
\]
In particular, when $p/N$ is such that $p/N\leq1-\eps$,
\[
P \bigl(\sqrt{\lambda_p}\leq\sqrt{\mathfrak{C}_n} \bigl[\bigl(1-\sqrt
{1-\eps}\bigr)-t \bigr] \bigr)\leq\exp\bigl(-(N-1)t^2 \bigr) .
\]

Interestingly, this bound is ``quite uniform'' in $\Lambda$, in the
sense that the only characteristics of $\Lambda$ that matter are
$\mathfrak{C}_n=C \frac{N-1}{n-1}$ and $N$.
\end{pf*}

\section{\texorpdfstring{Generalizations of the proof of Theorem~\lowercase{\protect\ref{thm:simplequadformmuhatandsigmahat}}}
{Generalizations of the proof of Theorem 4.3}}\label{appendix:SectionMuHatSigmaHat}

This part of the Appendix explains how to appropriately modify the
proofs of Theorems~\ref{thm:QuadFormsInverseEllipticalCase} and
~\ref{thm:SummaryQuadFormsMuHatAndSigmaHat} to
obtain the results we need in the case of correlated observations
(Section~\ref{Sec:Elliptical:Subsec:CorrelObs}) and the bootstrap.

\subsection{\texorpdfstring{On ${v\trsp\SigmaHat^{-1} v}$ when the observations are correlated}
{On v' Sigma -1 v when the observations are correlated}}\label{App:Sec:Generalizations:Subsec:QuadFormsV}

We explain in this subsection how to modify the proof of Theorem \ref
{thm:QuadFormsInverseEllipticalCase} in the case where the vectors of
observations $X_i$ and $X_j$ are potentially correlated. The data was
assumed to have the following representation, in matrix form:
\[
X=\ebold\mu\trsp+\Lambda Y \Sigma^{1/2} ,
\]
where $\Lambda$ is $n\times n$, deterministic but not necessarily
diagonal and $Y$ has i.i.d. ${\cal N}(0,1)$ entries. We also wrote the
SVD of $\Lambda$ as $\Lambda=ADB\trsp$, where $A$ and $B$ are orthogonal.

If we call $H=\id_n-\ebold\ebold\trsp/n$, we have, of course,
\[
\SigmaHat=\frac{1}{n-1} X\trsp H X=\frac{1}{n-1} \Sigma^{1/2}
Y\trsp\Lambda\trsp H \Lambda Y \Sigma^{1/2} .
\]
The orthogonality of $B$ implies $BY\equalInLaw Y$, and we have
\[
\SigmaHat\equalInLaw\frac{1}{n-1} \Sigma^{1/2}Y\trsp D (A\trsp H A)
D Y\Sigma^{1/2} .
\]
If we now call $\delta=A\trsp\ebold$, we see that $\| \delta\|
_2^2=n$, because $A$ is orthogonal. It can also easily be seen that
$A\trsp H A=\id_n-\delta\delta\trsp/n=H_{\delta}$. Because of the
remark we just made on\vadjust{\goodbreak} the norm of $\delta$, $H_{\delta}$ is clearly
an orthogonal projection matrix. So we have to understand
\[
\SigmaHat\equalInLaw\frac{1}{n-1} \Sigma^{1/2}Y\trsp D H_{\delta}
D Y\Sigma^{1/2} ,
\]
which is extremely close to the situation of Theorem \ref
{thm:QuadFormsInverseEllipticalCase}, where we had to work with
\[
\SigmaHat\equalInLaw\frac{1}{n-1} \Sigma^{1/2}Y\trsp D H_{\ebold}
D Y\Sigma^{1/2} .
\]
$D$ now plays the role $\Lambda$ played in Theorem \ref
{thm:QuadFormsInverseEllipticalCase} and the main modification is that
$H=H_{\ebold}$ is now replaced by $H_{\delta}$.

An examination of the proof of Theorem \ref
{thm:QuadFormsInverseEllipticalCase} shows that we never relied on the
fact that we used specifically $H_{\ebold}$ (instead of $H_{\delta}$)
in that proof. All we used was the fact that our $H$ there was a rank-1
perturbation of $\id_n$ and an orthogonal projection matrix.
Similarly, Lemma~\ref{lemma:LowerBoundOnSmallestEig}, on which we
relied in the course of the proof of Theorem \ref
{thm:QuadFormsInverseEllipticalCase}, handles $H_{\delta}$ for general
$\delta$ with squared norm $n$ without any problems, so it is still
usable in the course of the current study.

Because we know that the squared singular values of $\Lambda$ (and
hence the eigenvalues of $D$) satisfy \eqref
{eq:AssumptionSmallestSingularValueBoundedBelow}, the proof of Theorem
\ref{thm:QuadFormsInverseEllipticalCase} goes through without further
modifications and Proposition \ref
{prop:QuadFormsInverseEllipticalCaseCorrelated} holds.

\subsection{On quadratic forms involving random projection
matrices}\label{App:Subsec:GeneralQuadFormsProjMat}
A recurrent issue in the questions we addressed was the understanding
of statistics of the form
\[
\frac{1}{n} u\trsp P u ,
\]
where $P$ is a random projection matrix and $u$ a (generally
deterministic) vector of dimension $n$. In particular, the projection
matrices we dealt with were of the form
\[
P=\Lambda Y (Y\trsp\Lambda^2 Y)^{-1}Y\trsp\Lambda,
\]
for $\Lambda$ a (possibly random) $n\times n$ diagonal matrix and $Y$
an $n\times p$ matrix with i.i.d. ${\cal N}(0,1)$ entries. We also
assume that $\Lambda$ and $u$ are independent of~$Y$. Finally, we
assume that $\| u\|_2/\sqrt{n}=1$.

In the course of the text, we carried out successfully computations
when $u=\ebold$, but relied to do so on properties of $\operatorname{trace}(P)$.
The case of general $u$ is more involved and is treated here.

\begin{lemma}\label{lemma:GeneralQuadFormsRandomProj}
Assume that $\Lambda$ and u (which is deterministic) are such that
\[
\frac{1}{n^2} \sum_{i=1}^n u_i^4 \tendsto0 \quad  \mbox{and} \quad \frac
{1}{n^2} \sum_{i=1}^n \lambda_i^4 \tendsto0
\]
and that \eqref{eq:AssumptionSmallestSingularValueBoundedBelow} holds
for $\Lambda$ for a certain sequence $N(n)$.\vadjust{\goodbreak}

Under the preceding assumptions, we have, if $Z(u)=\frac{1}{n} u\trsp
P u $,
\[
Z(u)-\frac{1}{n}\sum_{i=1}^n u_i^2 \mathbf{E}(P(i,i)|\Lambda
)\tendsto0
\qquad \mbox{in probability}
\]
conditionally on $\Lambda$.
\end{lemma}

\begin{pf}
We simply sketch the modifications to the proof given after the
statement of Theorem~\ref{thm:simplequadformmuhatandsigmahat}. As
noted in Lemma~\ref{lemma:MeanRandomProjection}, the off-diagonal
elements of $P$ have mean 0 conditionally on $\Lambda$. Now, using the
same notations as in Theorem~\ref{thm:simplequadformmuhatandsigmahat},
we have, using equation \eqref{eq:keyExpansionProjectionMatrix} there,
if $Z_i(u)$ is the quantity obtained by replacing $\lambda_i$ by $0$
in $Z$, $r_i=W_i {\cal S}_i^{-1} Y_i$, $w_i=r_i\trsp u/n$ and $u_i$ is
the $i$th coordinate of $u$,
\[
Z(u)=Z_i(u)+\frac{1}{n} \frac{1}{1+\lambda_i^2 q_i} (-\lambda
_i^2 w_i^2+2\lambda_i u_i w_i +\lambda_i^2 u_i w_i ) .
\]
The expression between the parentheses is easily seen to be equal to
$(1+\lambda_i^2 q_i)u_i^2-(\lambda_i w_i -u_i)^2$. We get an analog
of equation \eqref{eq:ApproxZbyZi}
\[
Z(u)=Z_i(u)+\frac{u_i^2}{n} -\frac{1}{n}\frac{(\lambda_i
w_i-u_i)^2}{1+\lambda_i^2 q_i}.
\]
Clearly, from the definition of $w_i$, $w_i| \{Y_{(-i)},\Lambda
\}\sim{\cal N}(0,u\trsp W_i{\cal S}_i^{-2} W_i\trsp u/n^2).$
Since by assumption $\| u\|_2=\sqrt{n}$, we have
\[
0\leq u\trsp W_i{\cal S}_i^{-1} W_i\trsp u/n^2=u\trsp W_i(W_i\trsp
W_i)^{-1}W_i\trsp u/n \leq1
\]
because $W_i(W_i\trsp W_i)^{-1}W_i\trsp$ is an orthogonal projection
matrix (hence its eigenvalues are only 0 and 1) and $\| u/\sqrt{n}\|_2=1$.

So we are exactly in the situation we were in during the proof of
Theorem~\ref{thm:simplequadformmuhatandsigmahat}, except for a term in
$u_i^4$ that now appears in our bound on the variance. Hence, with our
extra assumption on $\| u\|_4^4/n^2$, we conclude similarly (after a
regularization step) that $Z(u)$ converges in probability, conditional
on $\Lambda$ to its conditional mean which is simply
\[
\frac{1}{n}u_i^2 \mathbf{E}(P(i,i)|\Lambda) .
\]
\upqed
\end{pf}

We remark that to get an analog of Theorem \ref
{thm:QuadFormsMuHatSigmaHatInvV}, where now
\[
\zeta=\frac{1}{n}u\trsp\Lambda Y{\cal S}^{-1} v ,
\]
one just needs to go through the proof and replace the $w_i$ appearing
there by the ``new''
$w_i=u\trsp W_i {\cal S}_i^{-1} Y_i/n$. Exactly the
same arguments go through when $\sum_{i=1}^n u_i^2\lambda_i^2/n$
remains bounded. So under this condition, $\zeta$ tends to zero in probability.

With the help of the previous lemma, we can now prove the gist of
Proposition~\ref{proposition:QuadFormsMuHatSigmaInverseCorrCase}.
\begin{fact}\label{fact:PropCorrCaseHolds}
Proposition~\ref{proposition:QuadFormsMuHatSigmaInverseCorrCase} holds.
\end{fact}

\begin{pf}
We note that Proposition \ref
{proposition:QuadFormsMuHatSigmaInverseCorrCase} is essentially an
application of the previous lemma, with appropriate change of notation.
Recall the notations from the proposition. We have $\XTilde=\Lambda Y
\Sigma^{1/2}$ and $\Lambda$, which is $n\times n$, has singular value
decomposition $ADB\trsp$. Also, ${\cal S}=\XTilde\trsp\XTilde/n$,
$\YTilde=B\trsp Y$, $F=\YTilde\trsp D^2 \YTilde/n$. Hence, in the
language of the proposition,
\[
\hat{m}\trsp{\cal S}^{-1}\hat{m}=\frac{1}{n^2}\ebold\trsp
AD\YTilde
F^{-1}\YTilde\trsp D A \ebold= \omega\trsp P \omega,
\]
where $P=D\YTilde(\YTilde\trsp D^2 \YTilde)^{-1}\YTilde
\trsp D$ and $\omega=A\trsp\ebold$. When the assumptions of the
proposition are in force, $\Lambda$ is deterministic and Lemma \ref
{lemma:GeneralQuadFormsRandomProj} applies; from which we conclude
\[
\hat{m}\trsp{\cal S}^{-1}\hat{m}-\frac{1}{n}\sum_{i=1}^n \omega_i^2
\mathbf{E}(P(i,i)) \tendsto0 \qquad \mbox{in probability}.
\]
This gives us the analog of Theorem~\ref{thm:simplequadformmuhatandsigmahat}.

To get the analog of Theorem~\ref{thm:QuadFormsMuHatSigmaHatInvV}, we
just need $\sum_{i=1}^n \omega_i^2 d_i^2/n$ to remain bounded, which
is an assumption stated in Proposition \ref
{proposition:QuadFormsMuHatSigmaInverseCorrCase}.
\end{pf}

\subsection{Bootstrap specific results}\label{App:Subsec:BootSpecResults}

\subsubsection*{Bootstrapping mean 0 Gaussian data}

Our analysis of the bootstrap problem requires an analysis similar to
the one we performed in the previous subsection. In particular, there
we have $u=\Lambda^{1/2} \ebold$, where $\Lambda$ contains the
bootstrap weights. Since those add-up to $n$, the assumption $\| u\|
_2^2=n$ was clearly satisfied. Also, in the situation where
$p/n\tendsto\rho\in(0,1-1/e)$, we are guaranteed that
\[
P^*=\Lambda^{1/2} Y (Y\trsp\Lambda Y)^{-1}Y\trsp\Lambda^{1/2}
\]
is well defined with high-probability. When conditioning on $\Lambda$,
we see that we can work only with the submatrix $\Lambda^*$ (of size
$n^*$) whose diagonal entries are nonzero. This submatrix has its
diagonal entries bounded away from 0 as they are at least equal to 1.
Also, using arguments similar to those given in the proof of Lemma~\ref
{lemma:LowerBoundOnSmallestEig}, we see that we can get a uniform (in
$\Lambda$) lower bound on the smallest singular value of $\Lambda Y$,
which holds with probability exponentially [in $(n^*-p)$] close to 1.

So now we assume that we are dealing with $\Lambda$ such that $n^*-p$
tends to $\infty$, the empirical distribution of $\Lambda$ goes to
$\Poisson(1)$ and $\sum\lambda_i^2/n^2 \tendsto0$. We also assume
that \eqref{eq:AssumptionSmallestSingularValueBoundedBelow} are
satisfied for this $\Lambda$. Finally, we assume that $\{\sum_{i=1}^n
\lambda_i^2/n\leq10\}$. We call the corresponding set of matrices
${\cal G}_{B_n}$.
When the diagonal entries of $\Lambda$ are drawn from a
multinomial$(\frac{1}{n},\ldots,\frac{1}{n},n)$ it is clear that
these conditions are satisfied with probability going to 1. The only
thing that might require an explanation\vadjust{\goodbreak} is why the condition $\{\sum
_{i=1}^n \lambda_i^2/n\leq10\}$ holds with probability going to 1.
The mean of $\sum_{i=1}^n \lambda_i^2/n$ clearly goes to 2, using the
marginal distribution of $\lambda_i$. On the other hand, the arguments
we gave in Proposition~\ref{prop:EmpDistBootWeights} show that its
variance goes to 0, so this quantity goes to 2 in probability and
therefore is less than 10 with probability going to 1.

The main question that we still have to address is that of the behavior of
\[
\frac{1}{n}\sum_{i=1}^n u_i^2 \mathbf{E}(P^*(i,i)|\Lambda)
\]
when $u_i^2=\lambda_i$. By definition,
\[
P^*(i,i)=\frac{1}{n} \lambda_i Y_i\trsp\Biggl(\frac{1}{n}\sum
_{i=1}^n \lambda_i Y_iY_i\trsp\Biggr)^{-1}Y_i=1-\frac{1}{1+\lambda
_i Y_i\trsp{\cal S}_i^{-1} Y_i/n} ,
\]
where ${\cal S}_i=\frac{1}{n} \sum_{j\neq i} \lambda_j Y_j Y_j\trsp
$. Now concentration arguments (see, e.g., Section~\ref
{Subsec:InequalityConstraints}) show that, if $\sigma_p({\cal S}_i)$
is the smallest singular value of ${\cal S}_i$,
\[
P \biggl( \biggl|\frac{Y_i\trsp{\cal S}_i^{-1} Y_i}{p}-\frac
{\operatorname{trace}({\cal S}_i^{-1})}{p} \biggr|>t \big|{\cal
S}_i^{-1}
\biggr)=\gO\bigl(\exp\bigl(-pt^2\sigma^2_p({\cal S}_i)/2\bigr)\bigr) .
\]
We also know that with overwhelming probability (measured over
$Y_{(-i)}= \{Y_1,\ldots,Y_{i-1},Y_{i+1},\ldots,Y_n \}$),
$\sigma_p({\cal S}_i)$ is bounded away from 0, conditionally on
$\Lambda$, when $\Lambda$ is such that \eqref
{eq:AssumptionSmallestSingularValueBoundedBelow} holds. (Note for
instance that ${\cal S}_i \succeq\sum_{i\neq j} Y_jY_j\trsp/n$ and
use Lemma~\ref{lemma:LowerBoundOnSmallestEig}.) Hence, we conclude that
\[
\frac{Y_i\trsp{\cal S}_i^{-1} Y_i}{p}\simeq\frac{\mathrm
{trace}({\cal S}_i^{-1})}{p} .
\]
Hence, conditionally on $\Lambda$,
\[
P^*(i,i)\simeq1-\frac{1}{1+\lambda_i  ({p/n})  ({\mathrm
{trace}({\cal S}_i^{-1})/p})} ,
\]
with very high-probability, that is, the probability that the
difference between the two is greater than $\lambda_i (p/n) t$ is $\gO
(\exp(-C(n^*-p)t^2))$ for a fixed $C$ (by arguments similar to those
given in Lemma~\ref{lemma:LowerBoundOnSmallestEig}).
In other respects, we note that rank-1 perturbation arguments give, if
${\cal S}=\frac{1}{n}Y\trsp\Lambda Y$,
\[
\operatorname{trace}({\cal S}_i^{-1})-\operatorname{trace}({\cal S}^{-1})=\frac
{\lambda_i}{n}
\frac{Y_i\trsp{\cal S}_i^{-2}Y_i}{1+\lambda_i Y_i\trsp{\cal
S}_i^{-1}Y_i/n} .
\]
In particular, when $\Lambda$ is such that \eqref
{eq:AssumptionSmallestSingularValueBoundedBelow} holds, by using a
union bound argument,
\[
P \biggl(\max_{i=1,\ldots,n} \biggl|\frac{\mathrm
{trace}({\cal S}_i^{-1})-\operatorname{trace}({\cal S}^{-1})}{p}
\biggr|>\eps\big|\Lambda
\biggr)\tendsto0.
\]
We also note that $\operatorname{trace}({\cal S}^{-1})/p\tendsto
\limitScaling$
conditionally on $\Lambda$, if $\Lambda$ is such that its empirical
distribution goes to $\Poisson(1)$, $\Lambda\in{\cal G}_{B_n}$ and
$p/n\tendsto\rho$.\vadjust{\goodbreak}

Therefore, we also have by a simple union bound argument, conditional
on $\Lambda$, and assuming that $\Lambda$ is such that its empirical
distribution goes to $\Poisson(1)$, $\Lambda\in{\cal G}_{B_n}$ and
hence $\sum\lambda_i^2/n$ is less than 10,
\[
\frac{1}{n}\sum_{i=1}^n \lambda_i P^*(i,i)\simeq1-\frac{1}{n}\sum
_{i=1}^n \frac{\lambda_i}{1+\lambda_i \rho_n \limitScaling} .
\]
Now when $\Lambda\WeakCv\Poisson(1)$, which we write $G$, and $\rho
_n\tendsto\rho$,
\[
\frac{1}{n}\sum_{i=1}^n \frac{\lambda_i}{1+\lambda_i \rho_n
\limitScaling}\tendsto\int\frac{\tau \,dG(\tau)}{1+\tau\rho
\limitScaling} .
\]
But in light of the Mar\v{c}enko--Pastur equation, we have, under
these circumstances,
\[
\frac{1}{n}\sum_{i=1}^n \lambda_i P^*(i,i)\tendsto1-\frac
{1}{\limitScaling}=\frac{\limitScaling-1}{\limitScaling} .
\]

We finally conclude that conditional on $\Lambda$ being in ${\cal
G}_{B_n}$ (whose probability goes to 1),
\[
(\muHat^*)\trsp(\SigmaHat^*)^{-1}\muHat^* \tendsto\frac{
 (\limitScaling-1)/\limitScaling}{1- ({\limitScaling
-1})/{\limitScaling}}=\limitScaling-1\geq\frac{\rho}{1-\rho} ,
\]
since we know that $\limitScaling\geq1/(1-\rho)$ when $G$ is
$\Poisson(1)$, since its mean is 1.

Similar arguments as the ones used in the proofs in the main body of
the paper show that the same convergence in probability result holds
unconditionally on $\Lambda$, the problem being to get bounds that are
uniform in $\Lambda$, when $\Lambda\in{\cal G}_{B_n}$.

Hence, an analog of Theorem~\ref{thm:simplequadformmuhatandsigmahat}
follows (with ${\cal P}_n$ probability going to 1), where the ratio
$\rho/(1-\rho)$ is replaced by $\limitScaling-1$. The analog of
Theorem~\ref{thm:QuadFormsMuHatSigmaHatInvV} follows from the
arguments given in Appendix~\ref{App:Subsec:GeneralQuadFormsProjMat},
if we can show, in the notation used there that $\sum_{i=1}^n (u_i
d_i)^2/n$ remains bounded with probability going to 1. Note that
$u_i=d_i=\sqrt{\lambda_i}$ here, where $\lambda_i$ are the bootstrap
weights, so we just need to show that $\sum_{i=1}^n \lambda_i^2/n$
remains bounded. But we did this when describing ${\cal G}_{B_n}$. 

We therefore have an analog of Theorem \ref
{thm:QuadFormsMuHatSigmaHatInvV} and also of Theorem \ref
{thm:SummaryQuadFormsMuHatAndSigmaHat} when bootstrapping Gaussian data.

\subsubsection*{Bootstrapping elliptically distributed data}
Finally, let us say a few words about what would happen if we replaced
the normality assumption for the $X_i$'s by an elliptical distribution
assumption. We focus on the case where $X_i=\lambda_i \Sigma^{1/2}
Y_i$, that is, the mean of the $X_i$'s is 0. The previous analyses make
clear that the key questions concern $v\trsp(\SigmaHat^*)^{-1} v$ and
$(\muHat^*)\trsp(\SigmaHat^*)^{-1}\muHat^*$.

The questions concerning $v\trsp(\SigmaHat^*)^{-1} v$ fall pretty
much directly under the study we have made of elliptical distributions,
since we\vadjust{\goodbreak} know, according to the proof of Theorem \ref
{thm:bootQuadFormsInverseEllipticalCaseNormalCase}, that
\[
\SigmaHat^*=\frac{1}{n-1} \Sigma^{1/2}Y\trsp\Lambda\trsp
D^{1/2}(\id_n-\delta\delta\trsp/n)D^{1/2}\Lambda Y \Sigma^{1/2} ,
\]
where $D$ is the diagonal matrix containing the bootstrap weights and
$\delta=D^{1/2} \ebold$. So, as long as $D^{1/2}\Lambda$ satisfies
\eqref{eq:AssumptionSmallestSingularValueBoundedBelow}, results
similar to Theorem \ref
{thm:bootQuadFormsInverseEllipticalCaseNormalCase} will hold.

The questions dealing with $(\muHat^*)\trsp(\SigmaHat^*)^{-1}\muHat
^*$ are more involved. Analyses similar to the ones performed above
show that the key quantity to understand is now
\[
\frac{1}{n}\ebold\trsp D\Lambda Y(Y\trsp\Lambda\trsp D \Lambda
Y)^{-1}Y\trsp\Lambda\trsp D \ebold=\frac{1}{n}u\trsp
P_{D^{1/2}\Lambda,Y} u ,
\]
where $P_{D^{1/2}\Lambda,Y}=D^{1/2}\Lambda Y (Y\trsp\Lambda\trsp D
\Lambda Y)^{-1} Y\trsp\Lambda D^{1/2}$ and $u=D^{1/2} \ebold$. The
analysis of this quadratic form can be carried out just like we did
above in the Gaussian case, that is, $\Lambda=\id_n$. However, the
remarks we made to get simplified expressions for the limit do not seem
to apply anymore: quantities of the type
\[
\frac{1}{n}\sum_{i=1}^n \frac{d_i}{1+\lambda_i^2 d_i \rho
\limitScaling} ,
\]
appear, where $\limitScaling$ is the solution of equation \eqref
{eq:defLimitScaling} with $G$ being the limit (if it exists) of the
empirical distribution of the random variables $\lambda_i^2 d_i$.
These quantities do not appear to simplify any further to yield a
clearer and more exploitable expression.
\end{appendix}

\section*{Acknowledgments}

I am very grateful to Nizar Touzi and Nicole El Karoui for several very
interesting discussions at the beginning of this project and for their
interest in it. I would also like to thank two anonymous referees for
their constructive comments and insightful questions.


\printaddresses

\end{document}